\documentclass[12pt, a4paper, twoside, notitlepage]{article}
\usepackage{amssymb,amssymb}
\usepackage{enumerate}
\usepackage{amsthm}
\usepackage{framed}
\usepackage{comment,url}
\usepackage[centertags]{amsmath}
\usepackage[english]{babel}
\usepackage{multimedia}
\usepackage[english]{babel}
\usepackage[latin1]{inputenc}
\usepackage{times}
\usepackage[T1]{fontenc}
\usepackage{color}
\usepackage{graphicx}
\usepackage{graphics}
\usepackage{pgf}
\usepackage{amssymb}
\usepackage{bbm}
\usepackage[centertags]{amsmath}
\usepackage[english]{babel}
\usepackage{multimedia}
\usepackage[english]{babel}
\usepackage[latin1]{inputenc}
\usepackage{times}
\usepackage[T1]{fontenc}
\usepackage{picinpar}    %for picture in paragraphs with text around
\usepackage{color}
\usepackage{graphicx}
\usepackage{graphics}
\usepackage{pgf}
\usepackage{todonotes}
\newcommand{\bb}{\mathbb}

\newcommand{\abs}[1]{\lvert#1\rvert}

\newcommand {\U}{{\bf u}}

\newcommand {\Div}{{ div}}

\newcommand {\Deriv}{\mathcal{D}^\alpha}

\newcommand{\dxran}{{\rm d}\mathbf{x}_t}
\newcommand{\dXran}{{\color{blue}{\rm d}\mathbf{X}_t}}

\newcommand{\rmd}{{\rm d}}

\def\bs#1{\boldsymbol{#1}}
\def\wh#1{\widehat{#1}}

\newcommand{\dede}[2]{\frac{\delta #1}{\delta #2}}

\newtheorem{theorem}{Theorem}[section]

\newtheorem{corollary}[theorem]{Corollary}
\newtheorem{definition}[theorem]{Definition}

\newenvironment{remark}[1][Remark]{\begin{trivlist}
\item[\hskip \labelsep {\bfseries #1}]}{\end{trivlist}}

\begin{document}

%%%%%%%%%%%%%%%%%%%%%%%%%%%%%%%%%%%%%%
%\begin{frontmatter}
\begin{titlepage}
%%%%%%%%%%%%%%%%%%%%%%%%%%%%%%%%%%%%%%%%%%%%%%%%%%%%%%%%%%%%%%%%%%%%%%%%%%%%%%%%%%%%%%%%%%%%%%%%
%\title{A Deterministic Idealized Climate Model and its Stochastic Extension}
\title{Hasselmann's Paradigm for Stochastic Climate Modelling based on Stochastic Lie Transport}
\author{D. Crisan, D. D. Holm, P. Korn}

%%%%%%%%%%%%%%%%%%%%%%%%%%%%%%%%%%%%%%%%%%%%%%%%%%%%%%%%%%%%%%%%%%%%%%%%%%%%%%%%%%%
%%%%%%%%%%%%%%

%%%\begin{frontmatter}
%
%
%
%
\maketitle
\begin{center}
{\it Dedicated to the memory of Charlie Doering}
\end{center}
\vskip0.5cm
%\begin{framed}
\begin{abstract}
A generic approach to stochastic climate modelling is developed for the example of an idealized Atmosphere-Ocean model that rests upon Hasselmann's paradigm for stochastic climate models. Namely, stochasticity is incorporated into the fast moving atmospheric component of an idealized coupled model by means of stochastic Lie transport, while the slow moving ocean model remains deterministic. More specifically the stochastic model SALT (stochastic advection by Lie transport) is constructed by introducing stochastic transport into the material loop in Kelvin's circulation theorem. The resulting stochastic model preserves circulation, as does the underlying deterministic climate model. A variant of SALT called LA-SALT (Lagrangian-Averaged SALT) is introduced in this paper. In LA-SALT, we replace the drift velocity of the stochastic vector field by its expected value. The remarkable property of LA-SALT is that the evolution of its higher moments are governed by linear deterministic equations.     
Our modelling approach is substantiated by establishing local existence results, first, for the deterministic climate model that couples compressible atmospheric equations to incompressible ocean equation, and second, for the two stochastic SALT and LA-SALT models. 
\end{abstract}
%\end{framed}

\newpage

\tableofcontents
%\setcounter{page}{1}
%\vskip2cm
%
%\vskip2cm
\end{titlepage}

%\end{frontmatter}
%\newpage
%% %\vspace{8cm}!single spacing
%% %\vspace{2cm}%!double spacing
%% \begin{frontmatter}
%% \vspace{0.5cm}
%% \newpage
% \vspace*{-1cm}
% \tableofcontents
% \begin{frontmatter}
%% \addtocontents{toc}{\protect\vspace*{1cm}}
% \end{frontmatter}
%
%% %%%%%%%%%%%%%%%%%%%%%%%%%%%%%%%%%%%%%%
% \newpage

%%%%%%%%%%%%%%%%%%%%%%%%%%%%%%%%%%%%%%%%%%%%%%%%
%%%%%%%%%%%%%%%%%%%%%%%%%%%%%%%%%%%%%%%%%%%%%%%%
 \section{Introduction}\label{sect_INTRO}
%%%%%%%%%%%%%%%%%%%%%%%%%%%%%%%%%%%%%%%%%%%%%%%%%%%%%%%%
%\todo[inline,color=pink]{PK: Below is a suggestion of a new introduction that tries to take into account some of the concerns Dan has formulated in our last meeting.}

%\color{blue}
Prediction of climate dynamics is one of the great societal and intellectual challenges of our time. The complexity of this task has prompted the formulation of {\it idealized} climate models that target the representation of selected spatio-temporal characteristics, instead of representing the full bandwidth of physical processes ranging from seconds to millennia and from centimeters to thousands of kilometers.  A climate model of full complexity would couple, for example, an atmospheric model, described by the compressible three-dimensional Navier-Stokes equations and a set of 
advection-diffusion equations for temperature and humidity, to an oceanic model, given by the three-dimensional incompressible Navier-Stokes equations and advection-diffusion equations for temperature and salinity. Each model would completed by thermodynamic relationships and physical parametrizations to account for non-resolved processes, and the boundary conditions would represent the physics of the air-sea interface and the ocean's mix layer. In contrast, idealized models tend to simplify these equations by reducing the number of state variables, terms in the equations, or spatial dimensions. The amount of simplification needed in each climate model is dictated by the climate processes under investigation and is often rationalized by heuristic scaling considerations.

In this paper we formulate a framework for deriving stochastic idealized climate models.
The deterministic version of the two types of stochastic climate model we derive here belongs to a class of idealized climate models that target the study of El-Nino-Southern Oscillation (ENSO). ENSO is an instability of the coupled atmosphere-ocean system that occurs with quasi-periodic frequency of 5-7 years. The fundamental instability mechanism for ENSO can and has been investigated with idealized models that consist of {\it two-dimensional} coupled atmosphere-ocean-equations (see e.g. \cite{ENSO_THEORY}). We provide more details on this class of deterministic models in Section \ref{subsect_DetModel}.

A conceptual picture of the integration of stochasticity into a climate model was formulated by Hasselmann \cite{Hasselmann1976}. In Hasselmann's paradigm, the atmosphere acts with high frequency on short time scales, represented as a stochastic white-noise forcing of the ocean. The integration of the atmosphere's stochastic white-noise forcing over long time scales produces a low-frequency response in the ocean. As a result of the back-reaction, a red spectrum of the atmosphere's climate fluctuation is produced which complies with a variety of observations of the internal variability of the climate system  \cite{Hasselmann1976}. For a description of Hasselmanns program in probabilistic terms we refer to \cite{Arnold2001}).  
 It is common practice in climate modelling to incorporate Hasselmann's paradigm as a stochastic perturbation of the initial conditions, or as a stochastic forcing to the right-hand side of the dynamical equations of a deterministic climate model, then model the range of stochastic effects by creating an ensemble of simulations (see e.g., \cite{Palmer_2019}). The particular choice of the stochastic perturbation must be based on the modelling objectives in the case at hand. 
%\ADD{
A concise mathematical framework for stochastic climate modelling was developed in \cite{MajdaTimoEinden}.  This approach relies also on scale separation in fast and slow dynamics following Hasselmann's paradigm but it differs methodically from ours in modelling the nonlinear self-interaction of the fast variables by means of a linear stochastic model such as an Ornstein-Uhlenbeck process and incorporates multiplicative noise while at the same maintaining energy conservation.
 
 While following Hasselmann's view in a general sense, in this paper we incorporate the stochasticity in a novel way, applied in two stages which both deviate from the established practice.
 
Our approach deviates in several respects from common practice. First, the path of a fluid element in the Lagrangian sense is assumed to be stochastic. This assumption injects stochasticity directly into the transport velocity of the atmospheric fluid dynamics, thereby transforming the governing equations into stochastic PDEs. Second, although our stochasticity is introduced  
\emph{ab initio} and not \emph{a posteriori} via external forcing, both stages of our stochastic models are transparently related to the deterministic model by the Kelvin Circulation Theorem. This fundamental connection facilitates the physical interpretation of the two stages of the stochastic models. Our modelling approach could be seen as an implementation of Hasselmann's program, since we couple the fast and stochastic atmosphere model to the slow and deterministic ocean model and we implement this through a new coupling mechanism that passes the {\it expectation} of the atmospheric wind forcing to the ocean. However, Hasselmann's paradigm discussed in \cite{Hasselmann1976,FrankignoulHasselmann1977} now has more than three thousand citations, so the present paper could equally well be considered as a footnote to Hasselmann's program. 

The two stages of our stochastic approach are called Stochastic Advection by Lie Transport (SALT) and Lagrangian Averaged Stochastic Advection by Lie Transport (LA-SALT). The two stages of our approach represent two different viewpoints or modelling philosophies depending on the time scales of the intended application. For SALT, atmospheric `weather' produces uncertainty in advection arising from motion on unresolved time scales. In LA-SALT, atmospheric `climate' is taken as the baseline, and the atmospheric `weather' is treated as a field of fluctuations around the climate baseline, as discussed in Ed Lorenz's famous lecture \cite{Lorenz1995}. The LA-SALT approach brings us back to Hasselmann's paradigm, which decomposes a general climate model into deterministic and stochastic parts. In LA-SALT, in addition, the ideas of McKean, Vlasov and Kac \cite{McKean}, \cite{Vlasov} 
\cite{Kac}\footnote{See also the seminal discussion in Sznitman \cite{sas} of the ``propagation of chaos'' introduced in \cite{Kac}.}  are applied to the deterministic climate description. 
Namely, the LA-SALT approach results in deterministic linear fluctuation equations that govern the dynamics of the climate statistics themselves, including variance, covariance and higher statistical moments. Within our framework these higher order statistical moments are governed by linear equations. This result offers potential computational advantages and opens new perspectives for the theoretical analysis of these moments.

In summary, this paper formulates two complementary stochastic idealized climate models called SALT and LA-SALT. The SALT climate model couples a stochastic PDE for the atmospheric circulation to a deterministic PDE for the circulation of the ocean. The stochasticity is incorporated by assuming that Lagrangian particles in the atmosphere follow a stochastic path given by a Stratonovich process which appears in the motion of the material loop in Kelvin's circulation theorem. The stochastic Lagrangian path of the material loop is a semimartingale stochastic process in the SALT approach and is a McKean-Vlasov  process in the LA SALT approach. Both the SALT and LA-SALT approaches are related to an underlying deterministic model via Kelvin's circulation theorem. We substantiate our modelling choices by anchoring them within an established class of idealized climate models, as well as providing a mathematical analysis that demonstrates by proving a local well-posedness theorem that the proposed stochastic climate models rest on a firm mathematical basis. 

The numerical simulations that would demonstrate the capabilities of the SALT and LA-SALT stochastic models, however, are beyond the scope of the present paper. The key element of such a numerical experiment is the sensible specification of the stochastic process. For the purpose of mathematical analysis carried out here, though, it is sufficient to assume that the stochastic process is of Stratonovich \emph{type}. In contrast, a numerical experiment would require one to chose a specific Stratonovich process by incorporating externally obtained information either from observations or from high-resolution simulations. For an example of the latter procedure in the context of the Euler fluid equations in two dimensions, we refer to  
\cite{COTTERetal2019}.

\color{black}

\bigskip

 In the remainder of the introduction we detail our modelling approach for the deterministic model in Section \ref{subsect_DetModel} and for the stochastic 
 model in Section \ref{subsect_StochModel}.   
 
 \bigskip
 %\todo[inline,color=yellow]{DC: For Peter and Darryl to draft.}

%\color{blue}
 \begin{enumerate}
  \item
 Main content of the paper
 \begin{enumerate}
  \item
  Adaptation of the deterministic Gill-Matsuno \cite{GILL,MATSUNO} class of ocean-atmosphere climate model (OACM) to the geometric variational framework. This adaptation produces a Kelvin circulation theorem which retains the transformation properties which are the basis for the remainder of the paper. These transformation properties are inherited from the variational framework. They enable the formulation of the deterministic and stochastic models in terms of the same type of Kelvin circulation theorem. 
  \item
  Derivations of the SALT and LA-SALT stochastic versions of the OACM, whose flows all possess the same geometric transformation properties. This shared geometric structure enables the analysis to develop sequentially from deterministic to stochastic models.
  \item
  Mathematical analysis for the deterministic, SALT and LA-SALT versions of the OACM. Specifically we prove existence and uniqueness of local solutions for the deterministic  OACM, the existence of a martingale solution 
  for the SALT version of the OACM and existence and uniqueness of local solution for the LA-SALT version. 
  
  \item
  Outlook -- open problems, including further pusuit of the predictive equations derived here for the dynamics of OACM statistics. 
  \item
  Two appendices provide details of the derivations of the deterministic and stochastic models using Hamilton's variational principle. 
  \end{enumerate}
 \item
 Plan of the paper 
  \begin{enumerate}[(1)]
  \item
  The introduction in Section \ref{sect_INTRO} explains that the present work is based on Hasselmann's program of fast-slow decomposition of the climate into deterministic and stochastic components. It also introduces the deterministic climate model upon which we implement Hasselmann's program and it compares the deterministic and stochastic models we treat in terms of their individual Kelvin circulation theorems. 
  \item
  Section \ref{SectSectMathAnalysis1} proves the local existence and uniqueness properties of our variational geometric adaptation of the deterministic Gill-Matsunoclimate model.
  \item
  Section \ref{Sec3} also discusses the analytical properties of the SALT in section \ref{sec-SALT} and LA-SALT in section \ref{LA-SALT-eqns-sec} stochastic models. 
  \item
  Section \ref{Sec4} provides a summary conclusion and specification of open problems for the SALT and LA-SALT OACM.
  \end{enumerate}
 \end{enumerate}

\color{black} 
  \bigskip

 \subsection{The Deterministic Climate Model}\label{subsect_DetModel}
%%%%%%%%%%%%%%%%%%%%%%%%%%%%%%%%%%%%%%%%%%%%%%%%
The model of the atmospheric component of our idealized climate model consists of  the compressible 2D Navier-Stokes equation coupled to an advection-diffusion equation for  temperature $\theta^a$. The atmospheric velocity field $\U^a$ transports the temperature that provides the gradient term of the velocity equation. The ocean component of the coupled system consists of a 2D incompressible Navier-Stokes equations and an equation for the oceanic temperature variable $\theta^o$ that is passively advected by the ocean velocity field $\U^o$. Here the pressure acts here as a Lagrange multiplier to impose incompressibility. More specifically, the deterministic coupled PDE's for the ocean and the atmosphere are given by

%%%%%%%%%%%%%%%%%%%%%%%%%%%%%%%%%%%%%%%%%%%%%%%%
%\begin{equation}%\begin{split} \label{COUPLED_SWE}
\begin{align}
\text{Atmosphere: }\quad&\frac{\partial \U^a}{\partial t}
+(\U^a\cdot\nabla)\U^a
+\frac{1}{Ro^a}\U^{a\bot}
+\frac{1}{Ro^a}\nabla \theta^a
=
\frac{1}{Re^a}\triangle\U^a,\label{COUPLED_SWE_VELOC_A}
\\
%%%%%%%%%%%%%%%%%%%%%%%\\
& \frac{\partial \theta^a}{\partial t}
+ (\U^a\cdot\nabla)\theta^a
= 
 \gamma(\theta^a - \theta^o)
+\frac{1}{Pe^a}\triangle \theta^a.\label{COUPLED_SWE_T_A}\\
%%%%%%%%%%%%%%%%%%%%%%%
\text{Ocean: }\quad&\frac{\partial \U^o}{\partial t}
+(\U^o\cdot\nabla)\U^o
+\frac{1}{Ro^o}\U^{o\bot}
+\frac{1}{Ro^o}\nabla (p^o+q^a)
=
\sigma(\U^o-\bar{\U}^a_{sol})
+\frac{1}{Re^o}\triangle\U^o,\label{COUPLED_SWE_VELOC_O}\\
%%%%%%%%%%%
&\frac{\partial \theta^o}{\partial t}
+(\U^o\cdot\nabla)\theta^o
=\frac{1}{Pe^o}\triangle\theta^o,\label{COUPLED_SWE_T_O}\\
&\Div(\U^o)=0,\label{COUPLED_SWE_INCOMPRESS_O}\\
\text{with }& \text{ initial conditions }\nonumber\\
&\U^a(t_0)=\U^a_0,\ \theta^a(t_0)=\theta^a_0, \ \U^o(t_0)=\U^o_0,\ \theta^o(t_0)=\theta^o_0\nonumber.
\end{align}
%with $\bar{\U}^a:=\U^a-\frac{1}{|\Omega|}\int_\Omega \U^a dx$. 
%\end{split}\end{equation}
In these equations, the ocean velocity $\U^o$ is coupled to the atmospheric velocity $\U^a$ and the atmospheric temperature $\theta^a$ is coupled to the oceanic temperature $\theta^o$. The {\it coupling constants} $\gamma,\sigma<0$ regulate the strength of the interaction between the two components. 

The velocity coupling between the compressible atmosphere and the incompressible ocean model deserves some consideration. To preserve the incompressibility of the oceanic velocity field during the coupling we apply the Leray-Helmholtz Theorem to decompose the atmospheric velocity $\U^a$into a solenoidal component $\U^a_{sol}$ and a gradient term $q^a$ such that $\U^a=\U^a_{sol}+\nabla q^a$. The gradient part is combined with the oceanic pressure. In a second step we remove the space average via 
$\bar{\U}:=\U-\frac{1}{|\Omega|}\int_\Omega \U dx$ such that the oceanic velocity fields remains in the space of periodic flows with vanishing average. This property allows to determine the oceanic pressure.
%\todo[inline]{DC: Can we discuss the signs of the constants 
%$\gamma,\sigma$ ? Are they not negative ?}
%In the coupling, the mean of the atmospheric velocity is removed to guarantee that the ocean velocity remains with zero mean if it is initialized with a velocity that has mean zero. This is necessary to avoid non-uniqueness difficulties of the ocean component. 
Physically, this step removes the rapid mean velocity of the atmosphere relative to the slower ocean velocity in the frame of motion of the Earth's rotation. This means the ocean momentum responds to the shear force, which is proportional to the difference between the local ocean velocity at a given time and the local deviation of the atmospheric velocity away from its mean velocity.   

%\todo[inline, color=pink]{PK: Text below needs to be improved. The link to ENSO modelling needs to be emphasized. Also the physics of coupling need to be explained better.}
%\todo[inline, color=pink]{PK: The relation to ENSO-modelling is exciting, it provides a benchmark for our approach! Can we improve ENSO predictions ? The  Cane-Zebiak model (also idealized) improved ENSO-predictions and provided insight into this coupled phenomena (in the 90's). ENSO was also modelled as a noise driven phenomena. We may revisit/revitalize ENSO-modelling that was fashionable in the 90's and then more or less fell asleep. This does not require a GCM but just the idealized model we are formulating.}

The model above belongs to the class of  {\it intermediate coupled model}. These models are much simpler than the coupled general circulation models of the atmosphere-ocean system that are used for climate research. Intermediate coupled models allow to study fundamental aspects of the atmosphere-ocean interaction. The most prominent example is {\it El Ni\~{n}o-Southern Oscillation (ENSO)} in the tropical Pacific. As originally hypothesized by Bjerknes in 1969 \cite{BJERKNES_1969} this climate phenomenon crucially depends on the coupled interaction of both ocean and atmosphere. 
According to Bjerknes stronger trade winds increase the upwelling in the east Pacific, thereby creating a temperature gradient in the sea-surface temperature
that amplifies the trade winds. This interaction between the trade winds and sea surface temperature in the tropical Pacific generates a quasi-periodic oscillation between the three ENSO-phases: the neutral phase, El Ni\~{n}o and La Ni\~{n}a. Intermediate coupled models have been used successfully to shed light on the fundamental principle of ENSO, thereby confirming Bjerknes hypothesis. 

\begin{figure}[h]
\begin{center}
\includegraphics*[width=0.75\textwidth, height=0.3\textheight]{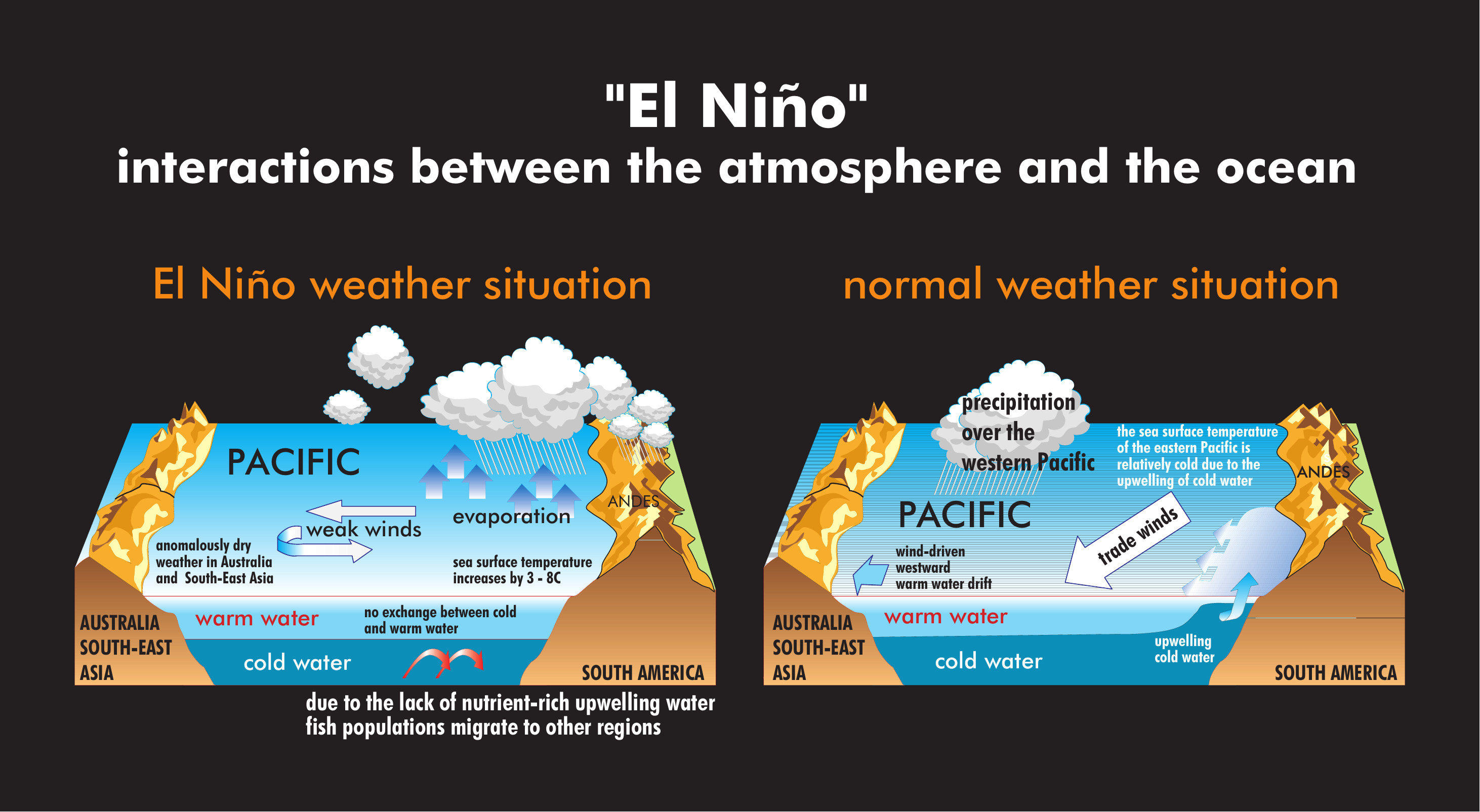}%{Idealized_coupledSystem.png}
\end{center}
\caption{\it\small Illustration of the dynamics and feedbacks of the atmosphere-ocean model that generate the El Ni\~{n}o-Southern Oscillation (ENSO). The trade winds that are part of the Walker circulation in the tropical Pacific interact with the cold/warm pools of the sea surface temperature. Local heating creates wind anomalies that in turn change the thermocline and the upwelling. During these process Rossby and Kelvin waves are emitted. The feedback of the ocean weakens the trade winds. }
  \label{fig:coupledSystem}
\end{figure}

The story of intermediate coupled models began with (uncoupled) models to study equatorial waves and their response to external forcing.
Matsuno \cite{MATSUNO} investigated an (uncoupled) divergent barotropic model (single layer of incompressible fluid of homogeneous density, with a free surface, on the beta plane)
\begin{equation}\begin{split}\label{MATSUNO}
&\frac{\partial \U}{\partial t}
+\frac{1}{Ro^a}\U^{a\bot}
+\frac{1}{Ro^a}\nabla \theta
=0,\\
%%%%%%%%%%%%
&\frac{\partial \theta}{\partial t}
+H\Div(\U)= Q.
\end{split}\end{equation}
Matsuno \cite{MATSUNO} refers to $\theta$ as {\it surface elevation} above a mean depth $H$, and in this context $Q$ appears as a source/sink of mass. Gill \cite{GILL} studied the steady response to heating anomalies of a tropical atmosphere, as described by the Matsuno model. Systems of equations in the following class are often called  {\it Gill models}
%\todo{Gill interpretation: $\theta=p$, with $p$ surface pressure} 
\begin{equation}\begin{split}\label{GILL}
&\frac{\partial \U}{\partial t}
+\frac{1}{Ro^a}\U^{a\bot}
+\frac{1}{Ro^a}\nabla \theta + a\U
=0,\\
%%%%%%%%%%%%
&\frac{\partial \theta}{\partial t}
+H\Div(\U) +b\theta= Q,
\end{split}\end{equation}
where $a,b$ are Raleigh friction and Newtonian cooling and where $Q$ is a heating term. In Gill's work $\theta$ is proportional to the {\it surface pressure}.
Since surface hydrostatic pressure is proportional to surface height, this identification is consistent with Matsuno's interpretation.

Atmospheric models of Gill-Matsuno type are often used to understand the atmospheric response during an El Ni\~{n}o to observed 
sea surface temperature anomalies. Zebiak \cite{ZEBIAK1982} parametrized the heat flux $Q$ from the ocean to the atmosphere in terms of the ocean sea surface temperature SST $Q=\alpha\, SST$. This relation can be motivated by a linearization of Clausius-Clapeyron relation (see \cite{ZEBIAK1986} and also \cite{HANEY}). 

As a next step, intermediate coupled models were constructed with atmospheric model either either of Gill-Matsuno type \cite{GILL} or as a statistical model of the atmosphere, that is for example constructed through an statistical analysis of atmospheric data (e.g. by Empirical Orthogonal Functions). Then the atmospheric model is coupled to a one or two-layer ocean model. Nonlinear terms are omitted. 
%Modelling El Nino-Southern Oscillation has been one of the main applications of this approach \cite{Latif_review}. 
The famous {\it Cane-Zebiak model} \cite{CANE_1985} applied a steady state atmosphere following Gill (\ref{GILL}) and a two-layer ocean model, with two equations for layer thickness and two equations for temperature. This model was used to issue the first ENO forecast \cite{CANE_1986}.
%These models have also been for the first successful ENSO forecast}, \cite{CANE_1986}. 
An overview can be found in chapter 7 of \cite{DIJKSTRA_BOOK}, or in \cite{ENSO_THEORY}. 

We have modified the equations in the model above by including the nonlinear terms in the velocity and temperature equations. First, we have replaced the damping terms due to Raleigh friction and Newtonian cooling in (\ref{GILL}) by Laplace operators for velocity and temperature. Furthermore, we interpret  atmospheric $\theta=h$ as surface elevation or equivalent depth and combined this with the {\it Charles Law} of thermodynamics (see e.g. \cite{DUTTON_BOOK}) according to which the volume $V$ is proportional to temperature $T$,  
$V=c\, T$, with constant $c>0$. Since the volume is also proportional to the surface elevation $V=c\, h $ one obtains $h=c\, T$. This identification allows one to interpret $\theta=T$ as the atmospheric temperature variable. This interpretation also relates the Matsuno and Gill equations (\ref{MATSUNO}) and (\ref{GILL}) to the atmospheric $\theta^a$-equation and allows one to interpret $\theta^a$ as a temperature variable.%
\footnote{For horizontal 2D models, the difference between potential and absolute temperature disappears.}
 
A distinct feature of ENSO is its pronounced irregularity, ENSO extremes occur irregularly in time and vastly differing amplitudes of sea surface temperature anomalies.
For the explanation of this irregular behaviour, two theories exist. Both theories rely on the separation of time scales shown by observations of the spectrum of variability in the tropical atmosphere-ocean system. Namely, a distinct time-scale separation is observed between the subseasonal and interannual oscillations. These observations suggest a natural decomposition of the dynamics in the tropics into short and long time scales. One theory explains the ENSO irregularity as a result of the chaotic dynamics exhibited by a nonlinear dynamical system through the interaction of slow components in which the fast components are less important. The other theory attributes the irregularity to a stochastic forcing of the slow modes by the fast modes, with applications to the Madden-Julian oscillation, and westerly/easterly wind bursts. The latter approach leads immediately to the question of which type of stochastic forcing, additive or multiplicative, is appropriate \cite{Perez_2005}. We do not aim to resolve this debate here. Rather, we suggest a new modelling approach for the stochastic theory of ENSO's irregularity. The new modelling approach is based on stochastic transport along the Lagrangian paths of advected fluid properties.

 \subsection{The Stochastic Atmospheric Climate Model}\label{subsect_StochModel}
%%%%%%%%%%%%%%%%%%%%%%%%%%%%%%%%%%%%%%%%%%%%%%%%

%%%%%%%%%%%%%%%%%%%%%%%%%%%%%%%%%%%%%%%%%%%%%%%%

%\subsubsection*{Kelvin circulation theorems for a sequence of atmosphere-ocean models}

The fundamental principle in modelling stochastic fluid advection is the Kelvin circulation theorem. As we shall see, each component of the deterministic atmosphere-ocean model in equations \eqref{COUPLED_SWE_VELOC_A} - \eqref{COUPLED_SWE_INCOMPRESS_O} above possesses its own Kelvin theorem, and the two components are coupled together by their relative velocity. The model \eqref{COUPLED_SWE_VELOC_A} - \eqref{COUPLED_SWE_INCOMPRESS_O} describes their interaction as the exchange of circulation between the atmosphere and ocean. Later we treat the atmospheric component of the model as being stochastic either in the sense of weather (SALT) or in the sense of climate (LA-SALT). In either case, the stochastic modification of the atmospheric dynamics will retain a Kelvin circulation theorems.

%\paragraph{Kelvin circulation theorems for the two components of the deterministic model.}  
\begin{theorem}[Kelvin theorem for the deterministic atmospheric model in \eqref{COUPLED_SWE_VELOC_A} - \eqref{COUPLED_SWE_INCOMPRESS_O}]$
\label{KelvinThm_atmosphere_determnistic}\,$\\
The deterministic model for atmospheric dynamics satisfies the following Kelvin theorem for circulation around a loop $c(u^a)$ moving with the flow of the atmospheric velocity $\bs{u}^a$. Namely,
\begin{align*}
\frac{d}{dt}\oint_{c(u^a)} (\bs{u}^a + \frac{1}{Ro^a}\bs{R}(\bs{x}))\cdot d\bs{x}
= \frac{1}{Re^a}\oint_{c(u^a)} \triangle\U^a \cdot d\bs{x}
\,,
\end{align*}
where ${\rm curl}\bs{R} = 2\bs{\hat{z}}\Omega(\bs{x})$ is the Coriolis parameter in nondimensional units.
\end{theorem}

\begin{proof}
By direct calculation, one shows that the deterministic atmospheric dynamics in the model above satisfies the relation in the Kelvin circulation theorem,
\begin{align*}
\frac{d}{dt}\oint_{c(u^a)} (\bs{u}^a + \frac{1}{Ro^a}\bs{R}(\bs{x}))\cdot d\bs{x}
&=
\oint_{c(u^a)} (\partial_t + \mathcal{L}_{u^a}) \big((\bs{u}^a + \frac{1}{Ro^a}\bs{R}(\bs{x}))\cdot d\bs{x}\big)
\\&= \oint_{c(u^a)} \Big( 
\partial_t \bs{u}^a + (\bs{u}^a\cdot\nabla)\bs{u}^a + u_j^a\nabla {u^a}^j
\\& \qquad- \bs{u}^a\times {\rm curl}\frac{1}{Ro^a}\bs{R}(\bs{x}) 
+ \nabla (\bs{u}^a\cdot\frac{1}{Ro^a}\bs{R})
\Big)\cdot d\bs{x} 
\\\hbox{By the model}\quad& =  \oint_{c(u^a)}\Big(- \frac{1}{Ro^a}\nabla\theta^a
+ \frac12\nabla |\bs{u}^a|^2 + \nabla (\bs{u}^a\cdot\frac{1}{Ro^a}\bs{R}) +  \triangle\U^a 
\Big)\cdot d\bs{x} 
\\& = \frac{1}{Re^a}\oint_{c(u^a)} \triangle\U^a \cdot d\bs{x}
\,.
\end{align*}
\end{proof}
\begin{remark}
In the proof above, $\mathcal{L}_u$ represents Lie derivative with respect to the vector field $u^a=\bs{u}^a \cdot\nabla$ with components $\bs{u}^a(\bs{x},t)$ and $c(u^a)$ denotes a material loop moving with the atmospheric Lagrangian transport velocity $\bs{u}^a(\bs{x},t)$.
Consequently, in the absence of viscosity, atmospheric circulation is conserved by the deterministic model because the viscous term is absent then and the loop integrals of gradients such as $u_j\nabla u^j=\frac12\nabla |\bs{u}|^2$ vanish on the right-hand side of the equation in the proof. \\
\end{remark}

Likewise, the dynamics of the ocean component of the model above satisfies the following Kelvin circulation theorem. 
\begin{theorem}[Kelvin theorem for the deterministic oceanic model in \eqref{COUPLED_SWE_VELOC_A} - \eqref{COUPLED_SWE_INCOMPRESS_O}]$\,$\\
The circulation dynamics around a loop $c(u^o)$ moving with the flow of the oceanic velocity $\bs{u}^o$ is given by
\begin{align*}
\frac{d}{dt}\oint_{c(u^o)} (\bs{u}^o + \frac{1}{Ro^o}\bs{R}(\bs{x}))\cdot d\bs{x}
= \oint_{c(u^o)} \Big( \sigma(\U^o-\bar{\U}^a)
+ \frac{1}{Re^o} \triangle\U^o 
\Big)\cdot d\bs{x} 
\,.
\end{align*}
\end{theorem}

\begin{proof}
The proof follows analogously to the proof of Theorem \ref{KelvinThm_atmosphere_determnistic}. 
%The deterministic oceanic dynamics in the model \eqref{COUPLED_SWE_VELOC_A} - \eqref{COUPLED_SWE_INCOMPRESS_O} above satisfies,
%\begin{align*}
%\frac{d}{dt}\oint_{c(u^o)}& (\bs{u}^o + \frac{1}{Ro^o}\bs{R}(\bs{x}))\cdot d\bs{x}
%=
%\oint_{c(u^o)} (\partial_t + \mathcal{L}_{u^o}) \big((\bs{u}^o + \frac{1}{Ro^o}\bs{R}(\bs{x}))\cdot d\bs{x}\big)
\\%\hbox{By the model}\quad& =  \oint_{c(u^o)}\Big(- \frac{1}{Ro^o}\nabla p^o
%+ \frac12\nabla |\bs{u}|^2 + \nabla (\bs{u}\cdot\frac{1}{Ro^o}\bs{R}) 
%+ \sigma(\U^o-\bar{\U}^a)
%+ \frac{1}{Re^o} \triangle\U^o 
%\Big)\cdot d\bs{x} 
%\\&=
%\oint_{c(u^o)} \Big( \sigma(\U^o-\bar{\U}^a)
%+ \frac{1}{Re^o} \triangle\U^o \Big)\cdot d\bs{x} \,.
%\end{align*}
\end{proof}

\subsection*{Stochastic Advection by Lie Transport (SALT) atmospheric model.} 
%\paragraph{Stochastic Advection by Lie Transport (SALT) atmospheric model.} 

Let $(\Xi ,\mathcal{F},(\mathcal{F}_{t})_{t},\mathbb{P})$ be a filtered probability space on which we have defined a sequence of independent Brownian motions  $(W^{i})_{i}$. Let $(\xi_{i})_{i}$ be a given sequence of sufficiently smooth vector 
fields that satisfies the condition in (\ref{xiassumpt}) below. In this work we assume the vector fields $(\xi_{i})_{i}$
to be given. For numerical simulations one defines these vector fields by extracting information from observational data. For an example we refer to \cite{COTTERetal2020c}.
%\footnote{To be precise, the vector fields $(\xi_{i})_{i}$ will be assumed to satisfy the condition in (\ref{xiassumpt}) below.}.
The derivation of the SALT atmospheric model introduces the stochastic Lagrangian path 
\begin{align}
\dxran := \U^a(\bs{x},t)dt-\sum_i\xi_i^a(\bs{x})\circ dW_i(t)
\,.\label{Atmos-SALT-Lag-Path}
\end{align}

Following \cite{Holm2015}, Appendix \ref{Appendix-SALT} discusses the introduction of the stochastic Lagrangian paths in \eqref{Atmos-SALT-Lag-Path} into Hamilton's variational principle for the atmospheric model equations. This step  leads to the SALT version of the idealized deterministic climate model comprising equations \eqref{COUPLED_SWE_VELOC_A} - \eqref{COUPLED_SWE_INCOMPRESS_O}. Namely, the SALT model is specified by the system of stochastic differential equations below:\\[3mm]

\noindent
Atmosphere:\footnote{%
As in the deterministic case we will write the Coriolis parameter as $\mathrm{curl}\,\mathbf{R}(%
\mathbf{x})=2\Omega \mathbf{(x}).$}
\begin{align}
& d\mathbf{u}^{a}+(d\mathbf{x}_{t}^{a}\cdot \nabla )\mathbf{u}^{a}+\frac{1}{%
Ro^{a}}d\mathbf{x}_{t}^{a\bot }+{\sum_{i}\Big(u_{j}^{a}\nabla \xi
_{i}^{j}+\frac{1}{Ro^{a}}\nabla \Big(R_{j}\mathbf{(x})\xi _{i}^{j}\Big)\Big)%
\circ dW_{t}^{i}}  \notag \\
& \hspace{3cm}+\frac{1}{Ro^{a}}\nabla \theta ^{a}=\frac{1}{Re^{a}}\triangle 
\mathbf{u}^{a},  \label{COUPLED_SWE_VELOC_A_STOCH} \\
& d\theta ^{a}+d\mathbf{x}_{t}^{a}\cdot \nabla \theta ^{a}=-\gamma (\theta
^{o}-\theta ^{a})+\frac{1}{Pe^{a}}\triangle \theta ^{a},
\label{COUPLED_SWE_T_A_STOCH} \\
& d\mathbf{x}_{t}^{a}=\mathbf{u}^{a}dt+{\sum_{i}\xi _{i}\circ
dW_{l}^{i}}  \label{stochasticpath}
\end{align}

\noindent Ocean: %%%%%%%%%%%
\begin{align}
\frac{\partial \mathbf{u}^{o}}{\partial t}+(\mathbf{u}^{o}\cdot \nabla )%
\mathbf{u}^{o}&+\frac{1}{Ro^{o}}\mathbf{u}^{o\bot }+\frac{1}{Ro^{o}}\nabla
p^{o}  \notag \\
& =\sigma (\mathbf{u}^{o}-\mathbb{E}\bar{\mathbf{u}^{a}})+\frac{1}{Re^{o}}%
\triangle \mathbf{u}^{o},  \label{COUPLED_SWE_VELOC_O_STOCH} \\
\frac{\partial \theta ^{o}}{\partial t}+(\mathbf{u}^{o}\cdot \nabla
)\theta ^{o}&=\frac{1}{Pe^{o}}\triangle \theta ^{o},
\label{COUPLED_SWE_T_O_STOCH} \\
{\ div}(\mathbf{u}^{o})&=0,  \label{COUPLED_SWE_INCOMPRESS_O_STOCH}
%\label{eq:CMSE}
\end{align}%

%\paragraph{Kelvin circulation theorem for the atmospheric component of the SALT model.}
\begin{theorem}Kelvin theorem for the SALT version of the atmospheric model in equations \eqref{COUPLED_SWE_VELOC_A_STOCH} - \eqref{stochasticpath}
\begin{align}
{\rm d}\oint_{c(\dxran)} (\bs{u}^a + \frac{1}{Ro^a}\bs{R}(\bs{x}))\cdot d\bs{x}
=  \frac{1}{Re}\oint_{c(\dxran)} \triangle\U^a \,dt \cdot d\bs{x}
\,,\label{Atmos-SALT-Kel}
\end{align}
where $c(\dxran)$ denotes any closed material loop whose line elements follow stochastic  Lagrangian paths 
as in \eqref{Atmos-SALT-Lag-Path}.
\end{theorem}
\begin{proof}
Upon suppressing the superscript $a$ in the velocity $\bs{u}^a$ for brevity of notation, we calculate
\begin{align*}
{\rm d}\oint_{c(\dxran)} (\bs{u}+ \bs{R}(\bs{x}))\cdot d\bs{x}
&=
\oint_{c(\dxran)} ({\rm d} + \mathcal{L}_{\dxran}) \big((\bs{u} + \bs{R}(\bs{x}))\cdot d\bs{x}\big)
\\&= \oint_{c(\dxran)} \Big( 
{\rm d} \bs{u} + (\rmd\bs{x}_t\cdot\nabla)\bs{u} + u_j \nabla {\rm d}x_t^j
\\& \qquad - \dxran\times {\rm curl}\bs{R}(\bs{x}) 
+ \nabla (\dxran\cdot\bs{R})
\Big)\cdot d\bs{x} 
\\\hbox{[By motion equation \eqref{COUPLED_SWE_VELOC_A_STOCH}]}\quad& =  \oint_{c(\dxran)}\Big(- \nabla\theta dt
+ \frac12\nabla |\bs{u}|^2dt - \dxran\times {\rm curl}\bs{R}(\bs{x})  + \nabla (\bs{u}\cdot\bs{R})dt
\\& \qquad  + { u_j  \nabla \sum \xi^j\circ dW(t)
+ \sum\nabla \big(\bs{\xi}\circ dW(t)\cdot\bs{R}\big)}
\Big)\cdot d\bs{x}
\\&=
\frac{1}{Re}\oint_{c(\dxran)} \triangle\U\,dt \cdot d\bs{x}
\,.
\end{align*}
%where the superscript $a$ is understood in the velocity $\bs{u}$. 
%(The additional {\color{red}red} terms above also appear in equation \eqref{COUPLED_SWE_VELOC_A_STOCH} below.)

\end{proof}

\begin{remark} The stochastic equation for the potential temperature $\theta^a$ in the atmospheric model inherits the stochasticity of the Lagrangian trajectories $\dxran$ in \eqref{Atmos-SALT-Lag-Path}, as a scalar tracer transport equation,
\begin{align}
{\rm d}\theta^a 
+ (\rmd\bs{x}_t\cdot\nabla)\theta^a
= 
 \big[ \gamma(\theta^a - \theta^o)
+\frac{1}{Pe^a}\triangle \theta^a\big]dt.
\label{Stoch_THETA_A}
\end{align}
\end{remark}

%\todo[inline, color=yellow]{PK: Isn't Kelvin's Theorem at this point of the introduction a Definition rather than a Theorem? 
%DH: I hope that by formulating the SALT equations first no confusion will now arise.}

%\todo[inline, color=yellow]{DH; I'm satisfied that the role of the SALT equations in this paper has now been made clear. Do you agree, Peter? }

%%%%%%%%%%%%%%%%%%%%%%%%%%%%%%%%%%%%%%%%%%%%%%%%%%%%%
\subsection{Lagrangian-Averaged Stochastic Advection by Lie Transport (LA-SALT) atmospheric  model.} 
%%%%%%%%%%%%%%%%%%%%%%%%%%%%%%%%%%%%%%%%%%%%%%%%%%%%%
%\paragraph{Lagrangian-Averaged Stochastic Advection by Lie Transport (LA-SALT) atmospheric  model.} 

\color{black}
We next modify the SALT approach to the two-dimensional atmospheric component of the climate system in the previous section to make it non-local in probability space, in the sense that the expected velocity will replace the drift velocity in the semimartingale for the SALT transport velocity of the stochastic fluid flow. This stochastic fluid model is derived by exploiting a novel idea introduced in \cite{DH2020} and developed further in \cite{AdLHT2020,DHL2020}, of applying Lagrangian-averaging (LA) in probability space to the fluid equations governing stochastic advection by Lie transport (SALT) which were introduced in \cite{Holm2015}.

The LA-SALT approach achieves three results of potential interest in climate modelling. These results address three different components of the climate change problem. 
\begin{itemize}
    \item First, the LA-SALT approach introduces a sense of determinism into climate science, by replacing the drift velocity of the stochastic vector field for material transport by its expected value in equation \eqref{Atmos-SALT-Lag-Path}. In this step, the expected fluid velocity becomes deterministic. 
    \item Second, the LA-SALT approach reduces the dynamical equations for the fluctuations to a \emph{linear} stochastic transport problem with a deterministic drift velocity. Such problems are well-posed. We prove here that the LA-SALT version of the SALT climate model a possesses local weak solutions. 
    \item Third, the LA-SALT approach addresses the dynamics of the variances of the fluctuations. In particular, the third result enables the variances and higher moments of the fluctuation statistics to be found deterministically, as they are driven by a certain set of correlations of the fluctuations among themselves. 
\end{itemize}
In summary, the first LA-SALT result makes the distinction between climate and weather for the case at hand. Namely, the LA-SALT fluid equations for the 2D atmosphere-ocean climate model system may be regarded as a dissipative system akin to the Navier-Stokes equations for the expected motion (climate) which is embedded into a larger conservative system which includes the statistics of the fluctuation dynamics (weather). The second result provides a set of linear stochastic transport equations for predicting the fluctuations (weather) of the physical variables, as they are driven by the deterministic expected motion. The third result produces closed deterministic evolutionary equations for the dynamics of the variances and covariances of the stochastic fluctuations.
Thus, the LA-SALT approach to investigating the 2D atmosphere-ocean climate model system treated here reveals that its statistical properties are fundamentally dynamical. 
Specifically, the LA-SALT analysis of the 2D atmosphere-ocean model presented here defines its climate, climate change, weather, and change of weather statistics, in the context of a hierarchical systems of PDEs and SPDEs with unique local weak solutions.
\color{black}

\paragraph{LA-SALT}
The expectation terms in LA-SALT induce another modification of the model which preserves the Kelvin circulation theorem, whose expectation yields a deterministic equation,
\begin{theorem}[Kelvin theorem for the LA-SALT atmospheric model]
\begin{align}
d\oint_{c(\dXran)} \big(\bs{u}+ \bs{R}(\bs{x})\big)\cdot d\bs{x}
&=
\frac{1}{Re}\oint_{c(\dXran)} \triangle\U \,dt\cdot d\bs{x}
\label{KelThm-LASALT}
\end{align}
where 
\[\dXran^a := \mathbb{E}[\U^a](\bs{x},t)dt+\sum_i\xi_i^a\circ dW_i(t).\]
\end{theorem}
\noindent 
{\bf Proof.} The proof of the Kelvin theorem for LA-SALT follows the same lines as for SALT. $\Box$

%%%%%%%%%%%%%%%%%%%%%%%%%%%%%%%%%%%%%%%%%%%%%%%%%%%%%%
\begin{comment}

\begin{proof}
\begin{align*}
d\oint_{c(\dXran)} \big(\bs{u}+ \bs{R}(\bs{x})\big)\cdot d\bs{x}
&=
\oint_{c(\dXran)} ({\rm d} + \mathcal{L}_{\dXran}) \big((\bs{u} + \bs{R}(\bs{x}))\cdot d\bs{x}\big)
\\&= \oint_{c(\dXran)} \Big( 
d \bs{u} + (\dXran \cdot\nabla)\bs{u}dt + u_j \nabla {\rm d}X_t^j
\\& \qquad - \dXran\times {\rm curl}\bs{R}(\bs{x}) 
+ \nabla (\dXran\cdot\bs{R})
\Big)\cdot d\bs{x} 
\\\hbox{By the model}\quad& =  \oint_{c(\dXran)}\Big(- \nabla\theta dt
+ (\mathbb{E}[\bs{u}]\cdot\nabla)\bs{u}dt  
+ u_j\nabla \mathbb{E}[u^j] dt + \nabla (\mathbb{E}[\bs{u}]\cdot\bs{R})dt
\\& \qquad
{\color{brown}
- \mathbb{E}[\bs{u}]\times {\rm curl}\bs{R}(\bs{x}) dt
- \bs{\xi}\circ dW(t)\times {\rm curl}\bs{R}(\bs{x})}
\\& \qquad {\color{red} + u_j \nabla \sum \xi^j\circ dW(t)
+ \sum\nabla \big(\bs{\xi}\circ dW(t)\cdot\bs{R}(\bs{x})\big)}
\Big)\cdot d\bs{x}
\\& 
=  \frac{1}{Re}\oint_{c(\dXran)} \triangle\U\,dt \cdot d\bs{x}
\,.
\end{align*}

\end{proof}

\end{comment}
%%%%%%%%%%%%%%%%%%%%%%%%%%%%%%%%%%%%%%%%%%%%%%%%%%%%%%

\paragraph{LA-SALT atmospheric equations in Stratonovich form.}

As shown in Appendix \ref{Appendix-LASALT}, the \emph{Stratonovich} LA-SALT equations are given 
in the standard notation for stochastic fluid dynamics by expanding out Kelvin's theorem in \eqref{KelThm-LASALT} to find
\begin{align}
d\mathbf{u}^{a}+({d\mathbf{X}_{t}}^{a}\cdot \nabla )\mathbf{u}%
^{a}+\frac{1}{Ro^{a}}{d\mathbf{X}_{t}}^{a\bot }& 
+{ \sum_{i}\Big(u_{j}^{a}\nabla \xi _{i}^{j}+\frac{1}{Ro^{a}}\nabla \Big(R_{j}%
\mathbf{(x})\xi _{i}^{j}\Big)\Big)\circ dW_{t}^{i}} 
 \notag\\
& \hspace{-22mm}{+\,u_{j}^{a}\nabla \mathbb{E}[{u^{a}}^{j}]dt+%
\frac{1}{Ro^{a}}\nabla (\mathbb{E}[{\mathbf{u}}^{a}]\cdot \mathbf{R})dt}+%
\frac{1}{Ro^{a}}\nabla \theta ^{a}\,dt=\frac{1}{Re^{a}}\triangle \mathbf{u}^{a}\,dt\,,
 \notag\\
& d\theta ^{a}+{d\mathbf{X}_{t}}^{a}\cdot \nabla \theta
^{a}=-\gamma (\theta ^{o}-\theta ^{a})\,dt + \frac{1}{Pe^{a}}\triangle \theta ^{a}\,dt\,.
\label{COUPLED_SWE_T_A_STOCH_LA-Strat-Intro}
\end{align}
where the \emph{Stratonovich} stochastic Lagrangian trajectory for LA-SALT is given by
\begin{equation}
{{\sf d}\mathbf{X}_{t}^{a}} 
:= {\mathbb{E}[\mathbf{u}^{a}]}(\bs{x},t)dt + \sum_{i}\xi _{i}^{a}(\bs{x})\circ dW_{i}(t).
\label{Lag-path-Strat}
\end{equation}

\paragraph{LA-SALT atmospheric equations in It\^o form.}
Likewise, the \emph{It\^o} LA-SALT equations are
given in the standard notation for stochastic fluid dynamics in Appendix \ref{Appendix-LASALT} by
\begin{align}
& d\mathbf{u}^{a}+({d\mathbf{\wh{X}}_{t}}^{a}\cdot \nabla )\mathbf{u}%
^{a}+\frac{1}{Ro^{a}}{d\mathbf{\wh{X}}_{t}}^{a\bot }
 + {
\sum_{i}\Big(u_{j}^{a}\nabla \xi _{i}^{j}+\frac{1}{Ro^{a}}\nabla \Big(R_{j}%
\mathbf{(x})\xi _{i}^{j}\Big)\Big) dW_{t}^{i}}  \nonumber \\
&+\frac12 \bigg[ \mathbf{\hat{z}}\times \xi \Big( {\rm div}\Big(\xi\,\big(\,\mathbf{\hat{z}}\cdot{\rm curl}\,(\,\mathbb{E}[{\mathbf{u}}^{a}] + \frac{1}{Ro^{a}}\mathbf{R}(\mathbf{x}) \big)\Big) \,\,\Big)  
- \nabla \bigg( \xi\cdot\nabla\Big(\xi \cdot\big(\mathbb{E}[{\mathbf{u}}^{a}] + \frac{1}{Ro^{a}}\mathbf{R}(\mathbf{x})\big) \Big)
\bigg)\bigg]dt
\nonumber \\& \hspace{22mm}{+\,u_{j}^{a}\nabla \mathbb{E}[{u^{a}}^{j}]dt+%
\frac{1}{Ro^{a}}\nabla (\mathbb{E}[{\mathbf{u}}^{a}]\cdot \mathbf{R})dt}+%
\frac{1}{Ro^{a}}\nabla \theta ^{a}\,dt = \frac{1}{Re^{a}}\triangle \mathbf{u}^{a}\,dt\,,
\label{COUPLED_SWE_VELOC_A_STOCH_LA-Ito-Intro} \\
& d\theta ^{a}+{d\mathbf{\wh{X}}_{t}}^{a}\cdot \nabla \theta
^{a} - \frac12 \Big(\xi\cdot\nabla(\xi\cdot\nabla \theta^a)  \Big)dt 
=-\gamma (\theta ^{o}-\theta ^{a})\,dt +  \frac{1}{Pe^{a}}\triangle \theta ^{a}\,dt
\,,
\label{COUPLED_SWE_T_A_STOCH_LA-Ito-Intro}
\end{align}
where the \emph{It\^o} stochastic Lagrangian trajectory for LA-SALT is given by
\begin{equation}
{{\sf d}\mathbf{\wh{X}}_{t}^{a}} 
:= {\mathbb{E}[\mathbf{u}^{a}]}(\bs{x},t)dt + \sum_{i}\xi _{i}^{a}(\bs{x})dW_{i}(t).
\label{Lag-path-Ito}
\end{equation}

\begin{remark}[Expected LA-SALT atmospheric equations.]
Taking the expectation of equations \eqref{COUPLED_SWE_VELOC_A_STOCH_LA-Ito-Intro} and  \eqref{COUPLED_SWE_T_A_STOCH_LA-Ito-Intro} yields a closed set of deterministic PDE for the expectations 
$\mathbb{E}[{\mathbf{u}}^{a}]$ and $\mathbb{E}[\theta ^{a}]$. Subtracting the expectations from equations \eqref{COUPLED_SWE_VELOC_A_STOCH_LA-Ito-Intro} and  \eqref{COUPLED_SWE_T_A_STOCH_LA-Ito-Intro} yields \emph{linear equations} for the differences,   
\begin{align}
{\mathbf{u}}^{a'}:= {\mathbf{u}}^{a} - \mathbb{E}[{\mathbf{u}}^{a}]
\quad\hbox{and}\quad
\theta ^{a'} := \theta ^{a} - \mathbb{E}[\theta ^{a}]
\,.
\label{fluctuations_u_theta}
\end{align}
Since ${\mathbf{u}}^{a'}$ and $\theta ^{a'}$ satisfy $\mathbb{E}[{\mathbf{u}}^{a'}]=0$ and $\mathbb{E}[\theta ^{a'}]=0$, one may regard these difference variables as fluctuations of ${\mathbf{u}}^{a}$ and $\theta ^{a}$ away from their expected values. From here one can calculate the dynamical equations for the statistics of the atmospheric model, e.g., its variances and its other tensor moments, as detailed in \cite{AdLHT2020,DH2020,DHL2020}. Further details of these equations can be found in Section \ref{LA-SALT-eqns-sec}
\end{remark}

\paragraph{Oceanic part of the LA-SALT model.}
The oceanic part of the LA-SALT model \emph{coincides} with the oceanic part of the SALT model (\ref{COUPLED_SWE_VELOC_O_STOCH})-(\ref{COUPLED_SWE_INCOMPRESS_O_STOCH}).

\newpage

%%%%%%%%%%%%%%%%%%%%%%%%%%%%%%%%%%%%%%%%%%%%%%%%%%%%%
\section{Local Existence and Uniqueness of the Deterministic Climate Model}\label{SectSectMathAnalysis1}
%%%%%%%%%%%%%%%%%%%%%%%%%%%%%%%%%%%%%%%%%%%%%%%%%%%%%

In the treatment below, we compactify the notation for the dynamics of the two-component system \eqref{COUPLED_SWE_VELOC_A} - \eqref{COUPLED_SWE_INCOMPRESS_O}, as follows. 
The state of the system is described by a state vector $\psi:=(\psi^a,\psi^o)$ with atmospheric component 
$\psi^a:=(\U^a,\theta^a)$ and oceanic component $\psi^o:=(\U^o,\theta^o)$. The initial state is denoted by
$\psi(t_0)=\psi_0$. where $\psi_0=(\U^a_0,\theta^a_0,\U^o_0,\theta^o_0)$.

In this notation, equations \eqref{COUPLED_SWE_VELOC_A} - \eqref{COUPLED_SWE_INCOMPRESS_O} take the operator form
\begin{equation}
d_t\psi+B\left( \psi,\psi\right) +C \psi +D(\psi^a,\psi^o)=L\psi,  \label{eq:CoupledOperator}
\end{equation}%
where one defines
\begin{itemize}
\item $B:=(B^a,B^o)$ is the usual bilinear transport operator, with
$B^a( \psi^a,\psi^a):=(\U^a\cdot\nabla\U^a,\U^a\cdot\nabla\theta^a)$ and 
$B^o( \psi^o,\psi^o):=(\U^o\cdot\nabla\U^o,\U^o\cdot\nabla\theta^o)$.
\item $C:=(C^a,C	^o)$ with $C^a\psi^a:=(\frac{1}{Ro^a}\U^{a\bot}+\nabla\theta^a,0)$ and
 $C^o\psi^o:=(\frac{1}{Ro^o}\U^{o\bot}+\nabla p^o,0)$
 \item $L:=(L^a,L^o)$ denotes the dissipation/diffusion operator for velocity and temperature with 
 $L^a\psi^a:=(\frac{1}{Re^a}\triangle\U^a,\frac{1}{Pe^a}\triangle\theta^a)$ and 
 $L^o\psi^o:=(\frac{1}{Re^o}\triangle\U^o,\frac{1}{Pe^o}\triangle\theta^o)$
 \item $D(\psi^a,\psi^o):=(0, \gamma(\theta^a - \theta^o),\sigma(\U^o-\bar{\U}^a_{sol}), 0)$ is the coupling operator.  
%\item $E_i$ are diagonal transport operators identically equal to zero for the components of $\psi$ corresponding to $\bar{\U}^o, \theta^o$.       
\end{itemize} 
%\begin{itemize}
%\item $B$ is the usual bilinear transport operator, 
%\item $C$ comprises all the linear terms (including the pressure term in the equation for the components of $\psi$ corresponding to $\bar{\U}^o$), \item $D$ is the diagonal matrix 
%$D=\mathrm{ diag}(0,0,0,1,1,0)$,\footnote{The same analysis applies to arbitrary matrices.} 
%%\item $E_i$ are diagonal transport operators identically equal to zero for the components of $\psi$ corresponding to $\bar{\U}^o, \theta^o$.       
%\end{itemize} 

\emph{Domain and Boundary Conditions:} The spatial domain is a two dimensional square $\Omega:=[0,L]\times [0,L]$ 
with $L\in \mathbb{R}^+$. We assume periodic boundary conditions.

\noindent \emph{Operators and Spaces:}  By $W^s(\Omega)$ we denote the $L^2$-Sobolev space of 
order $s\in\bb{Z}_+\cup \{0\}$ that is 
defined as the set of functions $f\in L^2(\Omega)$ 
such that its derivatives in the distributional sense $\Deriv f(x,y)=\partial^{\alpha_1}_{x}\partial^{\alpha_2}_{y} f(x,y)$ 
are in $L^2(\Omega)$ for all $|\alpha|\leq s$, with multi-index $\alpha=(\alpha_1,\alpha_2)\in\mathbb{Z}_+^2$, and
degree $\abs{\alpha}:=\alpha_1+\alpha_2$. The scalar product in 
$W^s(\Omega)$ is defined by
\begin{equation}\label{SOBOLEV_SCALATR_PROD}
 \big<f,g\big>_{W^s}:=\sum_{|\alpha|\leq s}\int_\Omega \Deriv f\cdot\Deriv g\, dx.  
\end{equation}
The vectorial counterpart of the Sobolev space $W^s(\Omega)$ is denoted by $\mathbf{W}^s(\Omega)$. 
More information about Sobolev spaces can be found for example in \cite{Evans}, \cite{Mazja}. 
We define the scalar space
\begin{equation}\begin{split}\label{V_def_scal}
{V}:=\{f:\mathbb{R}^2\to\mathbb{R}:\ &f\text{ is a trigonometric polynomial with period L}\},
%&\text{ and } \int_\Omega f\, dx=0\},
\end{split}\end{equation}
and its vector-valued equivalents for atmosphere and ocean component
\begin{equation}\begin{split}\label{V_def_vec}
&\mathbf{ V}^a:=\{u:\mathbb{R}^2\to\mathbb{R}^2:\ u\text{ is a vector-valued trigonometric polynomial with period L}\}\\
%&\text{period L\},\\
&\mathbf{ V}^o:=\{u:\mathbb{R}^2\to\mathbb{R}^2:\ u\text{ is a vector-valued trigonometric polynomial with period L}\\
&\hskip2.5cm \text{and }\int_\Omega u\, dx=0\}.
\end{split}\end{equation}
We define now the following function spaces
\begin{equation}\begin{split}\label{SOBOLEV}
&H^s(\Omega) := \text{the closure of }V\ \text{in }W^s(\Omega),\quad
\mathbf{H}^{s,a}(\Omega):= \text{the closure of }\mathbf{ V}^a\ \text{in }\mathbf{W}^s(\Omega),\\
&\mathbf{H}^{s,o}(\Omega):= \text{the closure of }\mathbf{ V}^o\ \text{in }\mathbf{W}^s(\Omega),\quad
\mathbf{H}^s_{div}(\Omega):=\{u\in {\bf H}^{s,o}(\Omega): \Div(u)=0\}.
\end{split}\end{equation}

In this notation, we define for $s\in\bb{N}\cup\{0\}$ the {\it Sobolev space of state vectors} by
\begin{equation}\label{Sobolev_StateSpace}
\mathcal{H}^s(\Omega):=
\mathbf{H}^{s,a}(\Omega)\times H^{s}(\Omega) \times \mathbf{H}^{s,o}_{div}(\Omega)\times H^{s}(\Omega),
\end{equation}
%Note that the first and third component of this state vector are vectorial Sobolev spaces, while the second and
%forth are scalar ones. 
in which the norm of $\psi=(\U^a,\theta^a,\U^o,\theta^o)\in \mathcal{H}^{{s}}$
is given by
\begin{equation}\label{Sobolev_StateSpaceNorm}
||\psi||_{\mathcal{H}^{s}}:=
(||\U^a||_{\mathbf{H}^{s}}^2 + ||\theta^a||_{H^{s}}^2+ ||\U^o||_{\mathbf{H}^{s}}^2 + ||\theta^o||_{ H^{s}}^2)^{1/2}. 
\end{equation}
We use an analogous notation for the Lebesgue spaces and denote by $L^2, {\bf L}^2, \mathcal{L}^2$ 
the sets of square-integrable scalar functions, vector fields and state vectors, respectively.

\begin{definition}\label{REGULAR_SOL}
Let $s\in\bb{N}$.
A state vector $\psi:=(\U^a,\theta^a,\U^o,\theta^o)$  is said to be a local regular solution of
(\ref{COUPLED_SWE_VELOC_A})-(\ref{COUPLED_SWE_INCOMPRESS_O}) 
on the time interval $T:=[t_0,t_1]$ if it satisfies(\ref{COUPLED_SWE_VELOC_A})-(\ref{COUPLED_SWE_INCOMPRESS_O}) 
with initial condition $\psi(t_0)=\psi_0$ and if 
\begin{equation}\begin{split}\label{REGULAR_SOL1}
\psi\in C(T,\mathcal{H}^s(\Omega))\cap L^2(T,\mathcal{H}^{s+1}(\Omega)), \qquad
{\frac{d\psi}{ dt}}\in \mathcal{H}^0 (T,\mathcal{H}^{s-1}(\Omega))
\end{split}.\end{equation}
\end{definition}
Define a cut-off function as follows
\begin{equation}\begin{split} \label{COUPLED_SWE_GALERKIN_TRUNC}
g_R(x):=
\begin{cases}
&1,\quad \text{if }0\leq x\leq R,\\
&0,\quad \text{if } x\geq R+\delta,\\
&\text{smoothly decaying}\quad \text{if }R<x<R+\delta
\end{cases}
\end{split}.\end{equation}
%\todo[inline]{DC: Why do you need $s\geq 3$.}
%\todo[color=pink,inline]{PK: I need $s\geq 3$, see eqs. (39)-(42) for the nonlinear terms and (47). See also comment near (47).}
Next, we define a finite-dimensional approximate system of equations for (\ref{eq:CoupledOperator}). Because of our assumption of periodic boundary conditions we may write the finite-dimensional approximation in terms of the Fourier basis 
$w_{\bf n}({\bf x}):=e^{\frac{2\pi i {\bf n}\cdot {\bf x}}{L}}$. We remark that the specific form of the basis does not play a role in our proofs. We must take into account the incompressibility of the ocean flow. This is imposed through the Leray projection, which projects the ocean equation onto the space of divergence-free vector fields. For periodic boundary conditions, the Leray projection commutes with the Laplace operator, so the Stokes operator coincides with the Laplacian. 

The Galerkin approximations for the atmospheric component $\psi^a=(\U^a,\theta^a)$ of the state vector are given by 
\begin{equation}\begin{split}\label{GALERKIN_U}
& %\U_m^a(x,t):=
{\bf P}_m^a\U^a(x,t):=\sum_{{\bf n}\in\{\mathbb{Z}^2, |{\bf n}|\leq m\}} 
\widehat{\U}^a_{\bf n}(t) w_{\bf n}(x)\\ \text{ and }&\  
 %\theta_m^a(x,t):=
 P_m^a\theta^a(x,t):=\sum_{{\bf n}\in\{\mathbb{Z}^2, |{\bf n}|\leq m\}} \hat{\theta}_{\bf n}^a(t) w_{\bf n}(x),
\end{split}\end{equation}
with $\widehat{\U}^a_{\bf n}:=\int_\Omega \U^a(x,t)w_{\bf n}(x)\, dx$, 
$\hat{\theta}_{\bf n}^a:=\int_\Omega \theta^a(x,t)w_{\bf n}(x)\, dx$. 
Incompressibility must be into account for the oceanic component $\psi^o=(\U^o,\theta^o)$, so the basis combines the Galerkin approximation with the Leray projection onto the space of divergence-free vector fields
\begin{equation}\begin{split}\label{GALERKIN_O}
& {\bf P}_m^o\U^o(x,t):=\sum_{{\bf n}\in\{\mathbb{Z}^2\setminus \{ 0\}, |{\bf n}|\leq m\}} 
\big(\widehat{\U}^o_{\bf n}(t)-\frac{\widehat{\U}^o_{\bf n}(t)\cdot{\bf n}}{|{\bf n}|^2}{\bf n} \big) w_{\bf n}(x), \\
%\qquad\text{with }\quad \widehat{\U}^o_{\bf n}\cdot {\bf n}=0,\ \widehat{\U}^o_{-{\bf n}}=\bar{\widehat{\U}}^o_{{\bf n}},\\ 
\text{ and }&\  
 P_m^o\theta^o(x,t):=\sum_{{\bf n}\in\{\mathbb{Z}^2, |n|\leq m\}} \hat{\theta}^o_{\bf n}(t) w_{\bf n}(x),
\end{split}\end{equation}
where $\widehat{\U}^o_{\bf n}:=\int_\Omega \U^o(x,t)w_{\bf n}(x)\, dx$, 
$\hat{\theta}_{\bf n}^o:=\int_\Omega \theta^o(x,t)w_{\bf n}(x)\, dx$. 

The Galerkin approximation of the state vector $\psi=(\psi^a,\psi^o)$ is defined as
\begin{equation}\begin{split}\label{GALERKIN_STATE}
P_m\psi:=({\bf P}_m^a\psi^a, {\bf P}^o_m\psi^o):=({\bf P}_m^a\U^a,P_m^a\theta^a, {\bf P}^o_m\U^o, P_m^o\theta^o).
\end{split}\end{equation}
We also use the notation 
\begin{equation}\begin{split}\label{GALERKIN_NOTATATION}
\psi_m:=P_m\psi, \text{ with } \U^a_m:={\bf P}_m^a\U^a, \theta_m^a:=P_m^a\theta, \text{ and }\U^o_m:={\bf P}^o_m\U^o, \theta_m^o:=P_m^o\theta^o.
\end{split}\end{equation}
In preparation for establishing the local existence of the stochastic version of the coupled model, we prove the following theorem on the global 
existence in time of the truncated approximation to the coupled model. 
%%%%%%%%%%%%%%%%%%%%%%%%%%%%%%%
\begin{theorem}[Global well-posedness of the truncated coupled model]
\label{THM_DETERMINISTIC_Galerkin}
For the time interval $[0,T]$, let $s\geq 2$ and suppose the initial conditions of (\ref{eq:CoupledOperatorGalerkin}) 
satisfy $\psi_{0}=(\U^a_{},\theta^a_{0},\U^o_{0},\theta^o_{0})\in \mathcal{H}^{s}(\Omega)$. 
The truncation of the coupled model (\ref{eq:CoupledOperator}) is given by
\begin{equation}\label{eq:CoupledOperatorTruncated}
d_t\psi+g_R(||\psi||_{H^s})B\left( \psi,\psi\right) +C \psi +D(\psi^a,\psi^o)=L\psi.
\end{equation}%
Then there exists a unique solution to  (\ref{eq:CoupledOperatorTruncated})  in the sense of Definition \ref{REGULAR_SOL}. This solution depends continuously with respect to the $\mathcal{L}^2$-norm on the initial conditions. 
\end{theorem}
%%%%%%%%%%%%%%%%%%%%%%%%%%%%%%%
\begin{proof} 
We first show the local existence in time of solutions for the truncated Galerkin system. Next, we prove the global existence via $\mathcal{H}^s$-estimates. Then, we pass to the limit and prove the corresponding assertions for the truncated system 
(\ref{eq:CoupledOperatorTruncated}). Finally, we show uniqueness and continuous dependency on the initial condition.\\
The truncated Galerkin system is given by
\begin{equation}\label{eq:CoupledOperatorGalerkin}
d_t\psi_m+g_R(||\psi_m||_{H^s})\mathbb{P}_m B\left( \psi_m,\psi_m\right) +C \psi_m +D(\psi^a_m,\psi^o_m)=L\psi_m,
\end{equation}%
where ${P}_m$ was defined in (\ref{GALERKIN_STATE}). 

%%%%%%%%%%%%%%%%%%%%%%%%%%%%%%%
\noindent{\it Step 1: Local existence of the truncated Galerkin approximation.}\\ 
%%%%%%%%%%%%%%%%%%%%%%%%%%%%%%%
The truncated Galerkin system can be written as 
\begin{equation}\begin{split}\label{eq:CoupledOperatorGalerkin_Proof1}
&d_t\psi_m=K(\psi_m)\\
\text{with }&K(\psi_m):=-g_R(||\psi_m||_{H^s})P_mB\left( \psi_m,\psi_m\right) -C \psi_m -D(\psi^a_m,\psi^o_m)+L\psi_m.
\end{split}\end{equation}
The right-hand side $K$ of  (\ref{eq:CoupledOperatorGalerkin_Proof1}) is a Lipschitz continuous mapping from 
$\mathcal{H}^{{s}}$ into itself. It follows from the Picard Theorem that a unique solution $\psi_m\in C^1([t_0,t_1^m], \mathcal{H}^{{s}})$ of (\ref{eq:CoupledOperatorGalerkin}) exists on time intervals $[t_0^m,t_1^m]$ that depend on $m$.
\noindent{\it Step 2: Global existence of the truncated Galerkin approximation.}\\
%%%%%%%%%%%%%%%%%%%%%%%%%%%%%%%
We show now that a solution exists globally in time. 

We apply the derivative $\Deriv $ to (\ref{eq:CoupledOperatorGalerkin}) %, where $|\alpha|=0,1$, 
and then take $L^2$-scalar product of equation (\ref{eq:CoupledOperatorGalerkin})  with $\Deriv\psi_m$. This yields
\begin{equation}\begin{split}\label{eq:CoupledOperatorGalerkin_Proof_s1_0}
&\frac{1}{2}d_t||\Deriv\psi_m||^2_{L^2}
+ \big<g_R(||\psi_m||_{H^1})\Deriv P_mB\left( \psi_m,\psi_m\right),\Deriv\psi_m\big>_{L^2} \\
&+\big<C \Deriv\psi_m,\Deriv\psi_m\big>_{L^2} +\big<D(\Deriv\psi^a_m,\Deriv\psi^o_m),\Deriv\psi_M\big>_{L^2}
-\big<L\Deriv\psi_m,\Deriv\psi_m\big>_{L^2}=0.
\end{split}\end{equation}

%%%%%%%%%%%%%%%%%%%%%%%%%%
\paragraph*{The nonlinear transport operator $B$.}
We rewrite the operator $B$ %we obtain with the inequalities of Agmon and Young (exponent 4, 4/3)
\begin{equation}\begin{split}\label{eq:CoupledOperatorGalerkin_Proof_s1_4}
&| \big<g_R(||\psi_m||_{H^s})\Deriv P_m\left( \psi_m\cdot\nabla\right )\psi_m,\Deriv\psi_m\big>_{L^2}\\
&=g_R(||\psi_m||_{H^s})| \big<\left( \Deriv\psi_m\cdot\nabla\right )\psi_m+ \left( \psi_m\cdot\nabla\right )\Deriv\psi_m,\Deriv\psi_m\big>_{L^2}
%\leq g_R(||\psi_m||_{H^2}) ||\Delta\psi_m||_{L^2}||\left( \psi_m\cdot\nabla\right )\psi_m||_{L^2}\\
%&\leq cg_R(||\psi_m||_{H^2}) ||\psi_m||_{H^2}||\psi||_{L^\infty}||\nabla\psi||_{L^2}
%\leq cg_R(||\psi_m||_{H^2}) ||\psi_m||_{H^2}^{3/2}||\nabla\psi||_{L^2}^{3/2}\\
%&\leq 
%cg_R(||\psi_m||_{H^2})\big( \frac{\epsilon_0}{2}||\psi_m||_{H^2}^2+\frac{1}{2\epsilon_0}||\nabla\psi||_{L^2}^{6}\big).
\end{split}\end{equation}
For the oceanic component the second term on the right-hand side vanishes as a consequence of the incompressibility. 
For the (compressible) atmospheric model the second term can be estimated with the inequalities of H\"older and Young
\begin{equation}\begin{split}\label{eq:CoupledOperatorGalerkin_Proof_s1_5}
&g_R(||\psi_m||_{H^s})| \big< \left( \psi_m\cdot\nabla\right )\Deriv\psi_m,\Deriv\psi_m\big>_{L^2}|
\leq
g_R(||\psi_m||_{H^s})|| \psi_m||_{L^6} ||\nabla\Deriv\psi_m||_{L^2} ||\Deriv\psi_m||_{L^3}\\
&\leq
g_R(||\psi_m||_{H^s})|| \psi_m||_{H^1} ||\mathcal{D}^{\alpha +1}\psi_m||_{L^2} ||\Deriv\psi_m||_{L^2}^{1/2}||\mathcal{D}^{\alpha +1}\psi_m||_{L^2}^{1/2}\\
&\leq
g_R(||\psi_m||_{H^s})||\Deriv\psi_m||_{L^2}^{3/2}||\mathcal{D}^{\alpha +1}\psi_m||_{L^2}^{3/2}\\
&\leq
\frac{c}{2\epsilon_1}g_R^4(||\psi_m||_{H^s})||\Deriv\psi_m||_{L^2}^{6}+ \frac{\epsilon_1}{2}||\mathcal{D}^{\alpha +1}\psi_m||_{L^2}^{2}.
\end{split}\end{equation}
For the first term in (\ref{eq:CoupledOperatorGalerkin_Proof_s1_4}) the inequalities of H\"older and Young imply that both the atmospheric and oceanic components satisfy

\begin{equation}\begin{split}\label{eq:CoupledOperatorGalerkin_Proof_s1_6}
&g_R(||\psi_m||_{H^s})| \big<\left( \Deriv\psi_m\cdot\nabla\right )\psi_m,\Deriv\psi_m\big>_{L^2}
\leq 
g_R(||\psi_m||_{H^s})||\Deriv\psi_m||_{L ^2}||\nabla\psi_m||_{L^6}||\Deriv\psi_m||_{L^3}\\
&\leq 
g_R(||\psi_m||_{H^s})||\Deriv\psi_m||_{L ^2}||\nabla\psi_m||_{H^1}||\Deriv\psi_m||_{L^2}^{1/2}||\mathcal{D}^{\alpha+1}\psi_m||_{L^2}^{1/2}\\
&\leq 
\frac{\epsilon_2}{2}||\Deriv\psi_m||_{L ^2}^2
+\frac{1}{\epsilon_2}g_R^2(||\psi_m||_{H^s})||\psi_m||_{H^2}^2||\Deriv\psi_m||_{L^2}||\mathcal{D}^{\alpha+1}\psi_m||_{L^2}\\
&\leq 
\frac{\epsilon_2}{2}||\Deriv\psi_m||_{L ^2}^2
+\frac{1}{\epsilon_2\epsilon_3}g_R^4(||\psi_m||_{H^s})||\psi_m||_{H^2}^4||\Deriv\psi_m||_{L^2}^2
+\frac{\epsilon_3}{2}||\mathcal{D}^{\alpha+1}\psi_m||_{L^2}^2.
\end{split}\end{equation}
\paragraph*{The linear operator $C$.}
%%%%%%%%%%%%%%%%%%%%%%%%%
We obtain using the inequalities of Cauchy-Schwarz and Young that
\begin{equation}\begin{split}\label{eq:CoupledOperatorGalerkin_Proof5}
&\int_\Omega \big(C\Deriv\psi_m\big)\cdot\Deriv\psi_m\, dx
=
\int_\Omega \big(C^a\Deriv\psi^a_m\big)\cdot\Deriv\psi^a_m\, dx
+
\int_\Omega \big(C^o\Deriv\psi^o_m\big)\cdot\Deriv\psi^o_m\, dx\\
&=\int_\Omega \big(\frac{1}{Ro^a}\Deriv\U^{a\bot}_m+\nabla\Deriv\theta^a_m\big)\cdot\Deriv\U^a_m\, dx
+\int_\Omega \big(\frac{1}{Ro^o}\Deriv\U^{o\bot}_m+\nabla\Deriv p^o_m\big)\cdot\Deriv\U^o_m\, dx\\
&\leq
|\int_\Omega \nabla\Deriv\theta^a_m\cdot\Deriv\U^a_m\, dx|
\leq \frac{1}{Ro^a}||\Deriv \theta^a_m||_{L^2}||\Deriv\nabla\U^a_m||_{L^2}
\leq 
\frac{1}{2\epsilon _4Ro^a}||\Deriv \theta^a_m||_{L^2}^2+\frac{\epsilon_4}{2}||\mathcal{D}^{\alpha+1}\U^a_m||_{L^2}^2,
\end{split}\end{equation}
in which the ocean component of the pressure term has vanished, upon using incompressibility of the ocean flow.

%%%%%%%%%%%%%%%%%%%%%%%%%
\paragraph*{The coupling operator $D$.}
%%%%%%%%%%%%%%%%%%%%%%%%%
For the coupling term in the atmospheric temperature equation we find by using the inequalities of Cauchy-Schwarz and Young that
\begin{equation}\begin{split}\label{CoupledModel_Estimate04}
&
|\gamma\int_\Omega\Deriv\big(\theta^o_m-\theta^a_m\big)\cdot\Deriv\theta^a_m\,dx|
\leq
%|\gamma|||\Deriv\big(\theta^o_m-\theta^a_m\big)||_{L^2}||\Deriv\theta^a_m||_{L^2}\\
%&\leq
|\gamma|\,(||\Deriv\theta^o_m||_{L^2}+||\Deriv\theta^a_m||_{L^2})||\Deriv\theta^a_m||_{L^2}\\
%&=\gamma(||\Deriv\theta^o_m||_{L^2}||\Deriv\theta^a_m||_{L^2}+||\Deriv\theta^a_m||_{L^2}^2)\\
%&=\gamma(\frac{1}{2}||\Deriv\theta^o_m||_{L^2}^2+\frac{1}{2}||\Deriv\theta^a_m||_{L^2}^2+||\Deriv\theta^a_m||_{L^2}^2)\\
&\leq
\frac{3|\gamma|}{2}(||\Deriv\theta^o_m||_{L^2}^2+||\Deriv\theta^a_m||_{L^2}^2).
\end{split}\end{equation} 
%\begin{equation}\begin{split}\label{CoupledModel_Estimate04}
%&
%%\int_\Omega \big(D(\Deriv\psi^a_m,\Deriv\psi^o_m)\big)\cdot\Deriv\psi\, dx
%|\gamma\int_\Omega\Deriv\big((\theta^o_m-\theta^a_m)\big)\cdot\Deriv\theta^a_m\,dx|
%\leq
%\frac{c\gamma}{2\epsilon_3}%||\Deriv\gamma||_{L^2}^2
%\big(||\Deriv\theta^o_m||_{L^2}^2+||\Deriv\theta^a_m||_{{\bf L}^2}^2\big)
%+\frac{\epsilon_3}{2}||\Deriv\nabla\theta^a_m||_{{\bf L}^2}^2.
%\end{split}\end{equation} 
The coupling term in the oceanic velocity equation can be estimated as follows
\begin{equation}\begin{split}\label{DERIV_MOMENTUM_O_12}
&|\sigma\int_\Omega\Deriv(\U^o_m-\bar{\U}^a_{sol,m})\cdot\Deriv\U^o_m\,dx|
=
|\sigma\int_\Omega\Deriv(\U^o_m-\U^a_{sol,m}+\frac{1}{|\Omega|}\int_{\Omega}\U^a_{sol,m}\,dx   )\cdot\Deriv\U^o_m\,dx|\\
%&\leq
%\sigma||\Deriv\U^o_m||_{{\bf L}^2}^2
%+\sigma\int_\Omega|\frac{1}{|\Omega|}\int_{\Omega}\Deriv\U^a_m\,dx \Deriv\U^o_m|\,dx\\
&\leq
|\sigma|\big(\,||\Deriv\U^o_m||_{{\bf L}^2}^2
+ ||\Deriv\U^o_{sol,m}||_{\bf L^2}||\U^a_{sol,m}||_{\bf L^2}
+\int_\Omega\frac{1}{|\Omega|}\int_{\Omega}|\Deriv\U^a_{sol,m}(x)|\,dx\big)\, |\Deriv\U^o_m(y)|\,dy\\
&\leq
|\sigma|\big(\,||\Deriv\U^o_m||_{{\bf L}^2}^2
+||\Deriv\U^o_m||_{\bf L^2}^2+||\U^a_{sol,m}||_{\bf L^2}^2
+\,\int_\Omega\frac{1}{|\Omega|}||\Deriv\U^a_{sol,m}||_{L^2}|\Omega| |\Deriv\U^o_m|\,dx\big)\\
&\leq
C|\sigma|\,(||\Deriv\U^o_m||_{L^2}^2+||\Deriv\U^a_{m}||_{L^2}^2).
\end{split}\end{equation} 
This estimate implies for the coupling operator
\begin{equation}\begin{split}\label{eq:CoupledOperatorGalerkin_Proof6}
&\int_\Omega \big(D(\Deriv\psi^a_m,\Deriv\psi^o_m)\big)\cdot\Deriv\psi\, dx
\leq 
%\frac{3\gamma}{2}(||\Deriv\theta^o_m||_{L^2}^2+||\Deriv\theta^a_m||_{L^2}^2)
%+C\sigma(||\Deriv\U^o_m||_{L^2}^2+||\Deriv\U^a_m||_{L^2}^2)\\
%&\leq
C(|\gamma|+|\sigma|)(||\Deriv\psi^o_m||_{L^2}^2+||\Deriv\psi^a_m||_{L^2}^2).
\end{split}\end{equation}
After summing over $|\alpha|$ up to $s$ this yields for (\ref{eq:CoupledOperatorGalerkin_Proof_s1_0})
\begin{equation}\begin{split}\label{eq:CoupledOperatorGalerkin_Proof12}
\frac{1}{2}d_t||\psi_m||_{\mathcal{H}^s}^2
+\frac{1}{P}||\nabla\psi_m||^2_{\mathcal{H}^s}
&\leq
Cg_R^4(||\psi_m||_{\mathcal{H}^s})||\psi_m||^6_{\mathcal{H}^s}
+C(|\gamma|+|\sigma|)(||\psi^o_m||_{\mathcal{H}^s}^2+||\psi^a_m||_{\mathcal{H}^s}^2).
%&\leq
%C\big(g_R(||\psi_m||_{\mathcal{H}^s})||\psi_m||_{\mathcal{H}^3}||\psi_m||^2_{\mathcal{H}^s}+||\psi_m||^2_{\mathcal{H}^s}\big),
%&\leq
%C\bigg(g_R(||\psi_m||_{\mathcal{H}^s})||\psi_m||_{\mathcal{H}^3}+1\bigg)||\psi_m||^2_{\mathcal{H}^s},
\end{split}\end{equation}
where $\frac{1}{P}:=\min\{\frac{1}{Re^a},\frac{1}{Re^o}, \frac{1}{Pe^a}, \frac{1}{Pe^o}\}$. 
Upon using the truncation $g_R$ it follows 
%\todo[color=pink,inline]{PK: for the deterministic case it seems sufficient to truncate just the $H^3$-norm. This allow already the application of Gronwall. Not sure if this is true for the stochastic model (?) }
\begin{equation}\begin{split}\label{eq:CoupledOperatorGalerkin_Proof13}
\frac{1}{2}d_t||\psi_m||_{\mathcal{H}^s}^2
+\frac{1}{P}||\nabla\psi_m||^2_{\mathcal{H}^s}
&\leq
CR||\psi_m||^2_{\mathcal{H}^s}.
%&\frac{C}{2Ro^a}||\theta^a_m||_{\mathcal{H}^s}^2\\
%+C(|\gamma|+|\sigma|)(||\psi^o_m||_{\mathcal{H}^s}^2+||\psi^a_m||_{\mathcal{H}^s}^2),
\end{split}\end{equation}
With Gronwall's inequality it follows that
\begin{equation}\begin{split}\label{eq:CoupledOperatorGalerkin_Proof14}
||\psi_m(t)||_{\mathcal{H}^{{s}}}^2
&\leq
||\psi_m(t_0)||_{\mathcal{H}^{{s}}}^2e^{ CR (t-t_0)}
\leq ||\psi(t_0)||_{\mathcal{H}^{{s}}}^2e^{ CR (t-t_0)},
\end{split}\end{equation}
where $\psi(t_0)$ denotes the initial condition of the coupled equations (\ref{COUPLED_SWE_VELOC_A})-(\ref{COUPLED_SWE_INCOMPRESS_O}). This estimate implies in particular that $||\psi_m(t)||_{\mathcal{H}^{{s}}}^2$ is bounded uniformly in $m$. 
Integrating (\ref{eq:CoupledOperatorGalerkin_Proof13}) over the time interval $[t_0,t]$ yields with 
(\ref{eq:CoupledOperatorGalerkin_Proof14}) 
\begin{equation}\begin{split}\label{eq:CoupledOperatorGalerkin_Proof15}
\frac{1}{R}\int_{t_0}^t||\nabla\psi_m(s)||_{\mathcal{H}^{{s}}}^2ds
&\leq
CR\int_{t_0}^t||\psi_m(s)||_{{\bf H}^{s}}^2ds+||\psi_m(t_0)||_{\mathcal{H}^{{s}}}^2\\
&\leq
CR||\psi(t_0)||_{\mathcal{H}^{{s}}}^2e^{ 2CR (t-t_0)}+||\psi(t_0)||_{\mathcal{H}^{{s}}}^2.
\end{split}\end{equation}
From (\ref{eq:CoupledOperatorGalerkin_Proof14}) and (\ref{eq:CoupledOperatorGalerkin_Proof15}) it follows that $(\psi_m)_m$ 
is uniformly bounded in $L^\infty(T,\mathcal{H}^{{s}})\cap L^2(T,\mathcal{H}^{{s}+1})$ with 
time derivative $(\frac{d\psi_m}{dt})_m$ which is according to (\ref{eq:CoupledOperatorGalerkin_Proof13}) uniformly bounded in 
$L^2(T,\mathcal{H}^{{s}})$.

The oceanic pressure $p$ can be recovered analogously to the Navier-Stokes Equations by solving the elliptic equation
\begin{equation}\label{pressure_eq}
\Delta p_m=div\big((\U^o_m\cdot\nabla)\U^o_m+\nabla q^a_m\big),
\end{equation}
where $q^a_m$ is the gradient part of the Leray-Helmholtz decomposition of the atmospheric velocity $\U^a_m$.
%Through this equation the pressure is adjusted accordingly to the nonlinearity and to the atmospheric velocity such that the resulting oceanic velocity is divergence-free. Equation (\ref{pressure_eq}) has a unique solution \textcolor{blue}{(to be continued with regularity of rhs. and implied regularity for pressure)} 

%%%%%%%%%%%%%%%%%%%%%%%%%%%%%%%%%%%%%%%%%%%%%%%%%%%%%%%%
\noindent{\it Step 3:} {\it Passage to the limit.}\\
%%%%%%%%%%%%%%%%%%%%%%%%%%%%%%%%%%%%%%%%%%%%%%%%%%%%%%%
The uniform boundedness of $(\psi_m)_m$ in $L^2(T,\mathcal{H}^{{s}+1})$ implies with the compact embedding of $L^2(T,\mathcal{H}^{{s}+1})$
into $L^2(T,\mathcal{H}^{{s}})$
that a subsequence $(\psi_{k})_k$ exists that converges strongly to $\psi\in L^2(T,\mathcal{H}^s)$. This subsequence converges also weakly in
$L^\infty(T,\mathcal{H}^{{s}})$.
%The uniform boundedness of $(\psi_m)_m$ in $L^\infty(T,\mathcal{H}^{{s}})\cap L^2(T,\mathcal{H}^{{s}+1})$
% implies together with 
%(\ref{eq:CoupledOperatorGalerkin_Proof13}) that $(\frac{d\psi_m}{dt})_m$ is uniformly bounded in 
%$L^2(T;\mathcal{H}^{{s}})$.
%The sequence $(\frac{d\psi_m}{dt})_m$ is in particular uniformly bounded in $L^2(T,\mathcal{L}^2)$ and 
%in $L^2(T;\mathcal{H}^{-s})$. The last fact follows from the continuous injection 
%$\mathcal{H}^{{s}}\subseteq \mathcal{L}^2\subset \mathcal{H}^{-(s+1)}$.
%Since the sequence $(\psi_m)_m$ is bounded in $L^2(T,\mathcal{H}^{{s}+1})$,  
%$(\frac{d\psi_m}{dt})_m$ is bounded in $L^2(T,\mathcal{H}^{-({s}+1)})$ and 
%since $\mathcal{H}^{{s}}$ is compactly embedded in $\mathcal{H}^{{s}+1}$
%Together with the continuous injection $\mathcal{H}^{{s}}\subseteq \mathcal{L}^2\subset \mathcal{H}^{-s}$
%follows with the  Aubin compactness theorem (cf. \cite{FoiasConstantin}, Lemma 8.2)
%that a subsequence $(\psi_{k})_k$ of $(\psi_m)_m$ exists that converges strongly to $\psi\in L^2(T,\mathcal{H}^{{s}^{'}})$
%for all ${s}^{'}\leq {s}-1$. 
We show now that the limit $\psi$ satisfies the truncated equations (\ref{eq:CoupledOperatorTruncated}).
%\todo[inline]{Dc: Peter see formulas (88) and (89).}
For the coupling term it holds for all $\phi\in [H^2(\Omega)]^6$
\begin{equation}\begin{split}\label{Coupled_CONV0}
&\lim_{k\to\infty}\int_T\big<D(\psi^a_k,\psi^o_k)-D(\psi^a,\psi^o) ,\phi\big>_{L^2}dt\\
%&=
%\lim_{k\to\infty}\int_T\big<\gamma(\theta^a_k-\theta^o_k)-\gamma(\theta^a-\theta^o)+\sigma(\U^o_k-\bar{\U}^a_k))-\sigma(\U^o-\bar{\U}^a)) ,\phi\big>_{L^2}dt\\
&=
\lim_{k\to\infty}\int_T\big<\gamma\big((\theta^a_k- \theta^a)+(\theta^o- \theta^o_k)\big)+\sigma\big((\U^o_k-\U^o)+(\bar{\U}^a_{sol}-\bar{\U}^a_{sol,k})\big) ,\phi\big>_{L^2}dt\\
\end{split}\end{equation}
For the velocity coupling involving the solenoidal part of the atmospheric velocity field in the ocean component it follows by using the Cauchy-Schwarz inequality that
\begin{equation}\begin{split}\label{Coupled_CONV01}
&|\int_T\big<\sigma\big(\bar{\U}^a_{sol}-\bar{\U}^a_{sol,k}\big) ,\phi\big>_{L^2}dt|
=|\sigma\int_T\int_\Omega
(\bar{\U}^a_{sol}(x,t)-\bar{\U}^a_{sol,k}(x,t)) \cdot\phi(x)\, dxdt|\\
=&|\sigma\int_T\int_\Omega\big(\U^a_{sol}(x,t)-
\U^a_{sol,k}(x,t)
+
\frac{1}{|\Omega|}\int_\Omega  \U^a_{sol,k}(z,t)- \U^a_{sol}(z,t)\, dz\big) \cdot\phi(x)\, dxdt|\\
\leq&
|\sigma|\,\int_T\int_\Omega| \U^a_{sol}(x,t)-\U^a_{sol,k}(x,t)|\, |\phi(x)|\,dxdt\\
\left.\right.&+|\sigma|\,\frac{1}{|\Omega|}\int_\Omega(\int_T\int_\Omega|  \U^a_{sol,k}(z,t)- \U^a_{sol}(z,t)|\, dzdt)   |\phi(x)|\,dx\\
\leq&
|\sigma|\,\int_T ||\U^a_{sol}(t)-\U^a_{sol,k}(t)||_{L^2}||\phi||_{L^2}dt
+|\sigma|\,\int_\Omega(\int_T||  \U^a_{sol,k}(t)- \U^a_{sol}(t)||_{L^2}dt) |\phi(x)|\,dx.
\end{split}\end{equation}
From (\ref{Coupled_CONV01}) and the convergence of $(\psi_{k})_k$ in $\psi\in L^2(T,\mathcal{H}^s)$,
there follows the convergence of the integral in  (\ref{Coupled_CONV0}). The convergence of the remaining linear terms in the equations is obvious. Next, we focus on the nonlinear terms for which we have to show that
\begin{equation}\begin{split}\label{Coupled_CONV1}
\lim_{k\to\infty}\int_T\big< g_R(||\psi_k||_{H^s})P_m B(\psi_k,\psi_k)-g_R(||\psi||_{H^s})P_mB(\psi,\psi),\phi\big>_{L^2}dt=0,\quad \text{for all }
\phi\in [C^\infty(\Omega)]^6.
\end{split}\end{equation}
The integral above can be written as
\begin{equation}\begin{split}\label{Coupled_CONV2}
&\int_T\big< g_R(||\psi_k||_{\mathcal{H}^s})P_m B(\psi_k,\psi_k)-g_R(||\psi||_{\mathcal{H}^s})P_m B(\psi,\psi),\phi\big>_{L^2}dt\\
&=
\int_T \big(g_R(||\psi_k||_{H^s})-g_R(||\psi||_{\mathcal{H}^s})\big)\big<P_m B(\psi_k,\psi_k),\phi\big>_{L^2}dt\\
&+\int_T g_R(||\psi||_{H^s})\big<P_m B(\psi_k,\psi_k) -P_mB(\psi,\psi),\phi\big>_{L^2}dt
\end{split}\end{equation}
For the first integral on the right-hand side it follows with the H\"older inequality
\begin{equation}\begin{split}\label{Coupled_CONV3}
&\int_T \big(g_R(||\psi_k||_{\mathcal{H}^s})-g_R(||\psi||_{\mathcal{H}^s})\big)\big<P_mB(\psi_k,\psi_k),\phi\big>_{L^2}dt\\
&\leq 
\int_T \big(g_R(||\psi_k||_{\mathcal{H}^s})-g_R(||\psi||_{\mathcal{H}^s})\big)||\psi_k||_{L^3} ||\nabla\psi_k||_{L^2} \,||\phi||_{L^6} dt\\
&\leq
c\int_T \big(g_R(||\psi_k||_{\mathcal{H}^s})-g_R(||\psi||_{\mathcal{H}^s})\big)||\psi_k||_{H^1}^2 \,||\phi||_{H^1} dt\\
&\leq
c\sup_{t\in T} \big(g_R(||\psi_k(t)||_{\mathcal{H}^s})-g_R(||\psi(t)||_{\mathcal{H}^s})\big)\int_T||\psi_k||_{H^1}^2 \,||\phi||_{H^1} dt.
\end{split}\end{equation}
%\begin{equation}\begin{split}\label{Coupled_CONV3}
%&\int_T \big(g_R(||\psi_k||_{H^s})-g_R(||\psi||_{H^1})\big)\big<B(\psi_k,\psi_k),\phi\big>_{L^2}dt\\
%&\leq 
%\int_T \big(g_R(||\psi_k||_{H^s})-g_R(||\psi||_{H^s})\big)||\psi_k||_{L^3} ||\nabla\psi_k||_{L^2} \,||\phi||_{L^6} dt\\
%&\leq
%c\int_T \big(g_R(||\psi_k||_{H^s})-g_R(||\psi||_{H^s})\big)||\psi_k||_{H^1} ||\nabla\psi_k||_{L^2} \,||\phi||_{H^1} dt\\
%&\leq
%c\int_T \big(g_R(||\psi_k||_{H^s})-g_R(||\psi||_{H^s})dt\big)^{1/2}\big(\int_T||\psi_k||_{H^1} ||\nabla\psi_k||_{L^2} \,||\phi||_{H^1} dt\big)^{1/2}\\
%&\leq
%K \big(\int_T||\psi_k||_{H^s}-||\psi||_{H^s}dt\big)^{1/2}\big(\int_T||\psi_k||_{H^1} ||\nabla\psi_k||_{L^2} \,||\phi||_{H^1} dt\big)^{1/2}
%\end{split}\end{equation}
The sequence $(\psi_k)_k$ converges weakly to $\psi$ in $L^\infty(T,\mathcal{H}^{{s}})$, i.e.
$||\psi_k||_{H^s}$ converges to $||\psi||_{\mathcal{H}^s}$ and with the continuity of the truncation function $g_R$ follows that the first term on the right-hand side converges to zero. Since 
$(\psi_k)_k$ is bounded in $L^\infty(T,\mathcal{H}^{{s}})\cap L^2(T,\mathcal{H}^{{s}+1})$ the right-hand side of (\ref{Coupled_CONV3}) converges for $k\to\infty$ to zero. 

For the second integral in  (\ref{Coupled_CONV2}) it follows with H\"older's inequality that
\begin{equation}\begin{split}\label{Coupled_CONV4}
&\int_Tg_R(||\psi||_{\mathcal{H}^s})\big<P_mB(\psi_k,\psi_k) -P_mB(\psi,\psi),\phi\big>_{L^2}dt\\
&=
\int_Tg_R(||\psi||_{\mathcal{H}^s})\big<P_mB(\psi_k-\psi,\psi_k) +P_mB(\psi,\psi_k-\psi),\phi\big>_{L^2}dt\\
&\leq
\int_Tg_R(||\psi||_{\mathcal{H}^s})|| \psi_k-\psi||_{L^4}||\nabla\psi_k||_{L^2}\, ||\phi||_{L^4} \,dt
+\int_Tg_R(||\psi||_{\mathcal{H}^s})||\psi||_{L^4}||\nabla(\psi_k-\psi)||_{L^2}||\phi||_{L^4} \,dt\\
&\leq
\int_Tg_R(||\psi||_{\mathcal{H}^s})|| \psi_k-\psi||_{H^1}||\nabla\psi_k||_{L^2}\, ||\phi||_{H^1} \,dt
+
\int_Tg_R(||\psi||_{\mathcal{H}^s})||\psi||_{H^1}||\psi_k-\psi||_{H^1}||\phi||_{H^1} \,dt,
\end{split}\end{equation}
where the right-hand side tends to zero for $k\to\infty$ as a consequence of the  
boundedness of the sequence $(\psi_k)_k$ in $L^\infty(T,\mathcal{H}^s)\cap L^2(T,\mathcal{H}^{s+1})$ and its converges in $L^2(T,\mathcal{H}^s)$.
\noindent{\it Step 4: Uniqueness of solutions of the truncated system}\\
%%%%%%%%%%%%%%%%%%%%%%%%%%%%%%%%%%%%%%%%%%%%%
Let $\psi_1,\psi_2$ be two solutions of (\ref{eq:CoupledOperatorTruncated}) with respective initial conditions $\psi_1(t_0=0)$ and 
$\psi_2(t_0=0)$. We assume that $||\psi_1(t)||_{\mathcal{H}^s},||\psi_2(t)||_{\mathcal{H}^s}\leq R$ for $t\in T$.
This implies for the  difference by $\widehat{\psi}:=\psi_1-\psi_2$ that $||\psi(t)||_{\mathcal{H}^s}\leq R$.

The difference $\widehat{\psi}$ satisfies the following equation
\begin{equation}\begin{split}\label{eq:CoupledOperatorGalerkin_uniq0}
&d_t\widehat{\psi}
+g_R(||\psi_{1}||_{\mathcal{H}^s})  \big( B\left( \psi_{1},\psi_{1}\right) -  B\left( \psi_{2},\psi_{2}\right)\big)\\
&+B\left( \psi_{2},\psi_{2}\right)\big(g_R(||\psi_{1}||_{\mathcal{H}^s})-g_R(||\psi_{2}||_{\mathcal{H}^s})\big)
 +C \widehat{\psi} +D(\widehat{\psi}^a,\widehat{\psi}^o)
=L\widehat{\psi}.
\end{split}\end{equation}
Taking the $L^2$-inner product with $\widehat{\psi}$ yields
\begin{equation}\begin{split}\label{eq:CoupledOperatorGalerkin_uniq1}
&d_t||\widehat{\psi}||_{\mathcal{L}^2}^2
+\int_\Omega g_R(||\psi_{1}||_{\mathcal{H}^s})  \big( B\left( \psi_{1},\psi_{1}\right) -  B\left( \psi_{2},\psi_{2}\right)\big)\cdot \widehat{\psi}dx\\
&+\int_\Omega B\left( \psi_{2},\psi_{2}\right)\cdot\widehat{\psi}\big(g_R(||\psi_{1}||_{\mathcal{H}^s})-g_R(||\psi_{2}||_{\mathcal{H}^s})\big)dx\\
& +\int_\Omega \big(C \widehat{\psi} +D(\widehat{\psi}^a,\widehat{\psi}^o)\big)\cdot\widehat{\psi}dx
=\int_\Omega (L\widehat{\psi})\cdot\widehat{\psi} dx.
\end{split}\end{equation}

For the difference of the two nonlinear terms we have
%\begin{equation}\begin{split}\label{eq:CoupledOperatorGalerkin_uniq2}
%&B\left( \psi_{1},\psi_{1}\right)- B\left( \psi_{2},\psi_{2}\right)
%=
%B( \widehat{\psi},\psi_{1})+B(\psi_{2},\widehat{\psi}).
%\end{split}\end{equation}
%Taking the $L^2$-inner product with $\widehat{\psi}_m$ yields
\begin{equation}\begin{split}\label{eq:CoupledOperatorGalerkin_uniq3}
&\int_\Omega g_R(||\psi_{1}||_{\mathcal{H}^s}) \big(B\left( \psi_{1},\psi_{1}\right)
-B\left( \psi_{2},\psi_{2}\right)\big)\cdot\widehat{\psi}\, dx
=
\int_\Omega g_R(||\psi_{1}||_{H^s}) \big( B( \widehat{\psi},\psi_{1})+B(\psi_{2},\widehat{\psi})\big)\cdot\widehat{\psi}\, dx
%\leq K||\widehat{\psi}||_{\mathcal{H}^s}
\end{split}\end{equation}
In order to estimate the right-hand side consider the atmospheric velocity component of (\ref{eq:CoupledOperatorGalerkin_uniq3}), 
with the inequalities of H\"older, Agmon and Young follows
\begin{equation}\begin{split}\label{eq:CoupledOperatorGalerkin_uniq4}
&\int_\Omega g_R(||\U_{1}^a||_{{\bf H}^s}) 
\big( (\widehat{\U}^a\cdot\nabla)\U_{1}^a)
+( \U^a_{2}\cdot\nabla)\widehat{\U}^a)\big)\cdot\widehat{\U}^a\, dx\\
&\leq
cg_R(||\U_{1}^a||_{{\bf H}^s}) ||\nabla\U_{1}^a ||_{{\bf L}^6}||\widehat{\U}^a||_{{\bf L}^3}||\widehat{\U}^a||_{{\bf L}^2}
+||\U_{2}^a ||_{{\bf L}^\infty}||\nabla\widehat{\U}^a||_{{\bf L}^2}||\widehat{\U}^a||_{{\bf L}^2}\\
&\leq
cg_R(||\U_{1}^a||_{{\bf H}^s}) ||\U_{1}^a ||_{{\bf H}^2}||\widehat{\U}^a||_{{\bf H}^1}||\widehat{\U}^a||_{{\bf L}^2}
+c||\U_{2}^a ||_{{\bf H}^2}||\widehat{\U}^a||_{{\bf H}^1}||\widehat{\U}^a||_{{\bf L}^2}\\
&\leq
\frac{cg_R(||\U_{1}^a||_{{\bf H}^s})}{2\epsilon_1} ||\U_{1}^a ||_{{\bf H}^2}^2||\widehat{\U}^a||_{{\bf L}^2}^2
+\frac{c}{\epsilon_2}||\U_{2}^a ||_{{\bf H}^2}^2||\widehat{\U}^a||_{{\bf L}^2}^2
+(\frac{\epsilon_1}{2}+\frac{\epsilon_2}{2})||\widehat{\U}^a||_{{\bf H}^1}^2\\
&\leq
M_0 ||\widehat{\U}^a||_{{\bf L}^2}^2
+(\frac{\epsilon_1}{2}+\frac{\epsilon_2}{2})||\widehat{\U}^a||_{{\bf H}^1}^2,
\end{split}\end{equation}
where $M_0=M_0(||\U_{1}^a ||_{{\bf H}^2}^2,||\U_{2}^a ||_{{\bf H}^2}^2)$.
Similarly we derive for the oceanic velocity component
\begin{equation}\begin{split}\label{eq:CoupledOperatorGalerkin_uniq5}
&\int_\Omega \big( (\widehat{\U}^o\cdot\nabla)\U_{1}^o)
+( \U^o_{2}\cdot\nabla)\widehat{\U}^o)\big)\cdot\widehat{\U}^o\, dx
\leq
M_1||\widehat{\U}^o||_{L^2}^2,
\end{split}\end{equation}
where $M_1=M_1(||\U_{1}^o ||_{H^2}^2)$, because the dependency on $||\U^o_2||_{H^s}$ vanishes due to  the incompressibility. Analogous estimates hold for the atmospheric and oceanic temperature transport terms, such that the nonlinear operator difference in (\ref{eq:CoupledOperatorGalerkin_uniq3}) can be bounded by
\begin{equation}\begin{split}\label{eq:CoupledOperatorGalerkin_uniq6}
&\int_\Omega g_R(||\psi_{1}||_{\mathcal{H}^s})\big(B\left( \psi_{1},\psi_{1}\right)
-B\left( \psi_{2},\psi_{2}\right)\big)\cdot\widehat{\psi}\, dx\leq K_0||\widehat{\psi}||_{\mathcal{L}^2}^2,
\end{split}\end{equation}
where $K_0=K_0(||\psi_1||_{\mathcal{H}^s}^2,(||\psi_2||_{\mathcal{H}^s}^2)$. 
The second term in (\ref{eq:CoupledOperatorGalerkin_uniq1}) can be estimated analogously as above
\begin{equation}\begin{split}\label{eq:CoupledOperatorGalerkin_uniq6a}
&\int_\Omega B\left( \psi_{2},\psi_{2}\right)\cdot\widehat{\psi}\big(g_R(||\psi_{1}||_{H^s})-g_R(||\psi_{2}||_{H^s})\big)dx\\
&\leq
c||\psi_2||_{\mathcal{H}^2}||\psi_2||_{\mathcal{H}^1}||\widehat{\psi}||_{\mathcal{L}^2}\big|||\psi_{1}||_{\mathcal{H}^s}-||\psi_{2}||_{\mathcal{H}^s}\big|
\leq
K_1||\widehat{\psi}||_{\mathcal{L}^2}^2,
\end{split}\end{equation}
where $K_1=K_1(||\psi_2||_{\mathcal{H}^2})$ and where we have used the reverse triangle inequality in the last step.

For the (linear) coupling operator it holds that
\begin{equation}\label{eq:CoupledOperatorGalerkin_uniq7}
|\int_\Omega D(\widehat{\psi}^a,\widehat{\psi}^o)\cdot \widehat{\psi}\, dx|
=
|\int_\Omega \gamma(\hat{\theta}^a - \hat{\theta}^o)^2+\sigma(\widehat{\U}^o-\overline{\widehat{\U}^a})^2\, dx|
\leq K_3||\widehat{\psi}||_{\mathcal{L}^2}^2
\end{equation}
where we have used that $\bar{\U}^a_1-\bar{\U}^a_2=\overline{\U^a_1-\U^a_2}=\overline{\widehat{\U}^a}$ and where $K_3$ depends on the coupling constants. This implies the following estimate for the difference equation (\ref{eq:CoupledOperatorGalerkin_uniq1})
\begin{equation}\label{eq:CoupledOperatorGalerkin_uniq8}
\frac{1}{2}d_t||\widehat{\psi}||^2_{\mathcal{L}^2}
+\frac{1}{R}||\nabla\widehat{\psi}||^2_{\mathcal{L}^2}
\leq K||\widehat{\psi}||_{\mathcal{L}^2}^2,
\end{equation}
where $K=K(||\psi_1||_{\mathcal{H}^2},||\psi_2||_{\mathcal{H}^2}), \sigma,\gamma)$.
From Gronwall's inequality we obtain
\begin{equation}\label{eq:CoupledOperatorGalerkin_uniq9}
||\widehat{\psi}(t)||^2_{\mathcal{L}^2}
\leq 
||\widehat{\psi}(t_0)||^2_{\mathcal{L}^2}e^{\int_{t_0}^t K(s)ds}.
\end{equation}
Since $\psi_1\psi_2\in L^2(T,\mathcal{H}^2(\Omega))$ the function $K$ is integrable and the right-hand side is bounded.
This proves the continuous dependency on the initial condition. If the 
two solutions have the same initial conditions, then the solutions coincide on $T$ and uniqueness follows. 

%\noindent{\it Step 5: Recovery of pressure in the ocean equations.} Analogously to the Navier-Stokes Equations the oceanic pressure $p$ can be recovered by solving the elliptic equation
%\begin{equation}\label{pressure_eq}
%\Delta p=div\big((\U^o\cdot\nabla)\U^o+\nabla q^a\big),
%\end{equation}
%where $q^a$ is the gradient part of the Leray-Helmholtz decomposition of the atmospheric velocity.
%%Through this equation the pressure is adjusted accordingly to the nonlinearity and to the atmospheric velocity such that the resulting oceanic velocity is divergence-free. Equation (\ref{pressure_eq}) has a unique solution \textcolor{blue}{(to be continued with regularity of rhs. and implied regularity for pressure)} 

\end{proof}

The following theorem is the main result for the deterministic version of the coupled model.

%%%%%%%%%%%%%%%%%%%%%%%%%%%%%%% 
\begin{theorem}[Local well-posedness of the coupled model]\label{THM_DETERMINISTIC_LOCAL} 
Let $s\geq 2$ and suppose the initial condition of the coupled equations (\ref{COUPLED_SWE_VELOC_A})-(\ref{COUPLED_SWE_INCOMPRESS_O}) 
satisfy $\psi_0=(\U^a_0,\theta^a_0,\U^o_0,\theta^o_0)\in \mathcal{H}^{s}(\Omega)$. 
Then there exists a unique time $t_1^{*}\in (t_{0}, \infty]$ such that a local regular solution $\psi$ of 
(\ref{COUPLED_SWE_VELOC_A})-(\ref{COUPLED_SWE_INCOMPRESS_O})
in the sense of Definition \ref{REGULAR_SOL} exists and is unique on any interval $T:=[t_0,t_1],$ where $t_0<t_1<t_1^*$ and that, if  $t_1^{*}< \infty$, then
\begin{equation}\label{explosionede} 
\lim_{t\nearrow t_1^{*}} ||\psi(t)||_{\mathcal{H}^{s}}=\infty.
\end{equation}
\end{theorem}
\begin{proof}
We define $t_R:=\inf\{t\geq t_0: ||\psi(t)||_{H^s}>R \}$, for $R>0$ and $\tau:=\lim_{R\to\infty} t_R]$. By $\psi_R$ we denote the solution of (\ref{eq:CoupledOperatorTruncated}) with initial condition $\psi_0$. We define $\psi(t):=\psi_R(t)$ for $t\in [t_0, t_R]$. 
On any time interval $[t_0,t_1]$ with $t_0<t_1< t_R$ the solutions
$\psi$ and $\psi_R$ coincide, as a consequence of the uniqueness of solutions of the truncated equations (\ref{eq:CoupledOperatorTruncated}).  
If $\tau=\infty$ then $\psi$ is a global solution of (\ref{eq:CoupledOperator}). If $\tau<\infty$ then $||\psi(\tau)||_{H^s}=R$ and $[t_0,\tau]$ is the maximal 
interval of existence of the solution $\psi$.
\end{proof}
\section{The Stochastic Idealized Atmospheric Climate Model}\label{Sec3}
%%%%%%%%%%%%%%%%%%%%%%%%%%%%%%%%%%%%%%%%%%%%%%%%%%%%%

\subsection{SALT atmospheric climate model}\label{sec-SALT}

Recall that the state of the system is described by a state vector $\psi
:=(\psi ^{a},\psi ^{o})$ with atmospheric component $\psi ^{a}:=(\mathbf{u}%
^{a},\theta ^{a})$ and oceanic component $\psi ^{o}:=(\mathbf{u}^{o},\theta
^{o})$. The initial state is denoted by $\psi (t_{0})=\psi _{0}$. where $%
\psi _{0}=(\mathbf{u}_{0}^{a},\theta _{0}^{a},\mathbf{u}_{0}^{o},\theta
_{0}^{o})$, with six entries. 

We  summarize equations \eqref{COUPLED_SWE_VELOC_A_STOCH}-%
\eqref{COUPLED_SWE_T_O_STOCH} for the SALT version of the idealized climate model as

\begin{equation}
d\psi _{t}+(B\left( \psi _{t},\psi _{t}\right) +C\left( \psi _{t}\right) +D%
\mathbb{E}[\bar{\psi}])dt+\sum_{i=1}^{\infty }E_{i}(\psi _{t})\circ
dW_{t}^{i}=\nu \Delta \psi _{t}dt,  \label{eq:CMSE}
\end{equation}%
where

\begin{itemize}
\item The process $\psi $ gathers all variables in %
\eqref{COUPLED_SWE_VELOC_A_STOCH}-\eqref{COUPLED_SWE_T_O_STOCH}, i.e., 
\[
\psi
_{t}=\{\psi
_{t}^i\}_{i=1}^6:=(\mathbf{u}_{t}^{a},\theta _{t}^{a},\mathbf{u}_{t}^{o},\theta
_{t}^{o})=({u}_{t}^{a,1},{u}_{t}^{a,2},\theta _{t}^{a},{u}_{t}^{o,1},{u}_{t}^{o,2},\theta
_{t}^{o}),
\]
as in the deterministic case.

\item $B$ is the usual bilinear transport operator,

\item $C$ comprises all the linear terms (including the pressure term in the
equation for the components of $\psi $ corresponding to $\bar{\mathbf{u}}%
^{o} $ as well as the term $-\gamma(\theta^a-\theta^o$) from 
\eqref{COUPLED_SWE_T_A_STOCH}).

\item $D=\{D_{ij}\}_{i,j=1}^6$ is a $6\times 6$-matrix that captures the influence of $\mathbb{E}[\bar{\psi}]$ on the various components of   $\psi
_{t}$. More precisely,  $D_{ij}$ is the coefficient appearing in front of 
$\mathbb{E}[\bar{\psi^j}]$ in the equation satisfied by $\psi^i
_{t}$. For the  SALT equations \eqref{COUPLED_SWE_VELOC_A_STOCH}-\eqref{COUPLED_SWE_T_O_STOCH},    
we have $D_{41}=D_{52}=-\sigma $,  with all the other
entries equal to 0. Therefore the pair $({u}_{t}^{o,1},{u}_{t}^{o,2})$ 
is affected by $(\mathbb E [\bar{u}_{t}^{a,1}], \mathbb E [\bar{u}_{t}^{a,2}]) $.

\item $E_{i}$ are diagonal operators given by 
\begin{equation*}
E_{i}(\mathbf{u}^{a},\theta ^{a},\mathbf{u}^{o},\theta ^{o})=\mathrm{diag}%
(\xi _{i}\cdot \nabla \mathbf{u}^{a}+\frac{1}{Ro^{a}}\xi _{i}+{u}%
_{j}^{a}\nabla \xi _{i}^{j}+\frac{1}{Ro^{a}}\nabla (R_{j}(\mathbf{x})\xi
_{i}^{j}),\ \xi _{i}\cdot \nabla \theta ^{a},\ 0,\ 0)
\end{equation*}

\item $\mathrm{curl}R({\mathbf{x}})=2\Omega ({\mathbf{x}})$.

\item $\bar{\psi}:=\psi -\frac{1}{|\Omega |}\int_{\Omega }\psi dx$. \footnote{%
Recall that subtraction of the mean $\frac{1}{|\Omega |}\int_{\Omega }\psi dx$
places the oceanic and atmospheric variables all into the same frame of motion relative to the
Earth's rotation. Note that this subtraction is only applied to $\mathbf{u}%
^{o}\,$. This is ensured by the multiplication by the matrix $D$.}
\end{itemize}

\noindent We start by giving a rigorous definition of the solution of \eqref{eq:CMSE}.
Let $(\Xi ,\mathcal{F},(\mathcal{F}_{t})_{t},\mathbb{P},(W^{i})_{i})$ be a
fixed stochastic basis. In addition to the Sobolev spaces defined above we
introduce $C^{m}(\Omega ;\mathbb{R}^{p})$ to be the (vector) space of all $%
\mathbb{R}^{p}$-valued functions $f$ which are continuous on $\Omega $ with
continuous partial derivatives $D^{\alpha }f$ of orders $|\alpha |\leq m$,
for fixed $m\geq 0$. Notice that on the torus all continuous functions are
bounded. The space $C^{m}(\Omega ;\mathbb{R}^{p})$ is a Banach space when
endowed with the usual supremum norm 
\begin{equation*}
\Vert f\Vert _{m,\infty }=\sum_{|\alpha |\leq m}\Vert D^{\alpha }f\Vert
_{\infty }.
\end{equation*}%
The space $C^{\infty }(\Omega ;\mathbb{R}^{p})$ is regarded as the
intersection of all spaces $C^{m}(\Omega ;\mathbb{R}^{p})$. Let $(\xi
_{i})_{i}$ a sequence of vector fields which satisfy the following
condition: 
\begin{equation}
\sum_{i=1}^{\infty }\Vert \xi _{i}\Vert _{s+3,\infty }^{2}<\infty .
\label{xiassumpt}
\end{equation}%
We will work directly with the It\^{o} version of (\ref{eq:CMSE}). In
this version, the Stratonovich integrals in (\ref{eq:CMSE}) are recast as It%
\^{o}~integrals with the required It\^{o} correction added in the drift term of the equation. The equivalence between
the two versions is straightforward, see e.g. \cite{cfh} for details.

\begin{definition}
\label{def:localglobalsolution}$\left. \right. $

\begin{enumerate}
\item[a.] A pathwise \textsl{local}{\ solution} of the system (\ref{eq:CMSE}%
) is given by a pair $(\psi ,\tau ),$ where $\tau :\Xi \rightarrow \lbrack
t_{0},\infty ]$ is a strictly positive stopping time and $%
\psi:\Omega \times \lbrack t_{0},\infty ]\rightarrow 
\mathcal{H}^{s}(\Omega )$, is an $\mathcal{F}_{t}$-adapted process with
initial condition $\psi _{t_{0}}\in \mathcal{H}^{s}(\Omega )$ such that 
\begin{equation*}
\psi \in L^{2}\left( \Xi ;C\left( [{t_{0}},T ];%
\mathcal{H}^{s}(\Omega )\right) \right) 
\ \ \psi 1_{[t_0,\tau]}
\in  L^{2}\left( \Xi ;L^{2}\left( [{t_{0}},T ];\mathcal{%
H}^{s+1}(\Omega )\right) \right)
\end{equation*}%
%for any $T>{t_{0}}$ 
for any $T\ge t_0$ and the system (\ref{eq:CMSE}) is satisfied locally
i.e., the following identity%
\begin{equation}
\psi _{t}=\psi _{{t_{0}}}-\int_{{t_{0}}}^{t\wedge \tau%
}\left( F(\psi _{s})+D\mathbb{E}\left[ \bar{\psi}_{s}\right]
\right) ds-\sum_{i=1}^{\infty }\int_{{t_{0}}}^{t\wedge \tau}E_{i}(\psi
_{s})dW_{s}^{i},  \label{eq:CMSElocal}
\end{equation}%
holds $\mathbb{P}$-almost surely, as an identity in $L^{2}\left( \Xi ;%
\mathcal{H}^{0}(\Omega )\right) $ for any $t\in [t_0,\infty)$. In (\ref{eq:CMSElocal}), the mapping $F(\psi _{s})$
is defined as 
\begin{equation}
F(\psi _{s})=B\left( \psi _{s},\psi _{s}\right) +C\left( \psi _{s}\right) -%
\frac{1}{2}\sum_{i=1}^{\infty }E_{i}^{2}(\psi _{s})-\nu \Delta \psi _{s}.
\label{eq:Fofpsi}
\end{equation}

\item[b.] A \emph{martingale local} solution of equation \eqref{eq:CMSE} is a triple $(\check\Omega, \mathcal{\check F}, \check{\mathbb{P}}), (\mathcal{\check F}_t)_t, (\check\psi, \check\tau, (\check W^i)_i)$ such that\\ $(\check\Omega, \mathcal{\check F}, \check{\mathbb{P}})$ is a probability space, $(\mathcal{\check F}_t)_t$ is a filtration defined on this space,   
where $\check\tau :\Xi \rightarrow \lbrack
t_{0},\infty ]$ is a strictly positive $\mathcal{\check F}_{t}$-stopping time and $%
\psi:\Omega \times \lbrack t_{0},\infty ]\rightarrow 
\mathcal{H}^{s}(\Omega )$, is an $\mathcal{\check F}_{t}$-adapted process with
initial condition $\psi _{t_{0}}\in \mathcal{H}^{s}(\Omega )$ such that 
\begin{equation*}
\check\psi \in L^{2}\left( \Xi ;C\left( [{t_{0}},T ];%
\mathcal{H}^{s}(\Omega )\right) \right) 
\ \ \check\psi 1_{[t_0,\check\tau]}
\in  L^{2}\left( \Xi ;L^{2}\left( [{t_{0}},T ];\mathcal{%
H}^{s+1}(\Omega )\right) \right)
\end{equation*}%
%for any $T>{t_{0}}$ 
for any $T\ge t_0$ and which satisfies equation \eqref{eq:CMSElocal}+\eqref{eq:Fofpsi} with 
$\psi$ replaced by  $\check\psi$ 
\footnote{We use the "check" notation $(\check{\,\phantom\,})$ in the description of the various components of a martingale solution, to emphasize that the existence of a martingale solution does not guarantee that, for a \emph{given} set of Brownian motions $(W^i)_i$ defined on the (possibly different) probability space $(\Omega, \mathcal{F}, {\mathbb{P}})$ a solution of \eqref{eq:CMSE} will exist. Clearly the existence of a strong solution implies the existence of a martingale solution.
}.
\item[c.] If $\tau=\infty $, then we say that the system (\ref%
{eq:CMSE}) has a \emph{global}{\ solution}.
In this case can remove the usage of the stopping time from  equation %
(\ref{eq:CMSElocal}). In other words, we have that%
\begin{equation}
\psi _{t}=\psi _{{t_{0}}}-\int_{{t_{0}}}^{t}\left( F(\psi _{s})+D\mathbb{E}[%
\bar{\psi}_{s}]\right) ds-\sum_{i=1}^{\infty }\int_{{t_{0}}}^{t}E_{i}(\psi
_{s})dW_{s}^{i},  \label{eq:CMSEglobal}
\end{equation}%
holds as an identity in $L^{2}\left( \Xi ;\mathcal{H}^{0}(\Omega )\right) $.\end{enumerate}
\end{definition}

\begin{remark}\label{impl}
Observe that $\psi_t=\psi_\tau$ for any $t\ge \tau$, it other words the solution remains constant once it hits the defining stopping time. We require this to be able to make sense of the quantity 
$\mathbb{E}[\bar{\psi}_{s}]$ even for temporal values $s$ larger than $\tau$. In fact, equation (\ref{eq:CMSElocal}) can be re-written as   
\begin{equation}
\psi _{t}=\psi _{{t_{0}}}-\int_{{t_{0}}}^{t\wedge \tau}%
\left( F(\psi _{s})+D\mathbb{E}\left[ \bar{\psi}_{s\wedge \tau}\right] \right)
ds-\sum_{i=1}^{\infty }\int_{{t_{0}}}^{t\wedge \tau}E_{i}(\psi
_{s})dW_{s}^{i},  \label{eq:CMSERlocalwrong}
\end{equation}%
\end{remark}
  
{\bf \noindent Roadmap of the section}

In the following we show that equation (\ref{eq:CMSE}) has a martingale solution provided the additional condition \eqref{adp} is satisfied. To do so, we follow the same route as in the deterministic case. In Theorem \ref{truncatedgalerkinscm} we show that the Galerkin approximations of a truncated version of equation 
(\ref{eq:CMSElocal}) are well defined globally.
Moreover we show that we can control the Sobolev norms of these approximations uniformly in the 
level of approximation, see
\eqref{eq:goodboundG} and \eqref{lpc}. 
The truncation is 
done by multiplying each the coefficients 
of (\ref{eq:CMSElocal}) with the function
$g_{R,\delta}(\left\vert \left\vert
\cdot\right\vert \right\vert _{\mathcal{H}^{s}(\Omega )})$, where, as in the deterministic case, we use the cut-off function $g_{R,\delta}:\mathbb{R}_{+}\rightarrow
\lbrack 0,1]$
\begin{equation*}
g_{R,\delta}(x):=%
\begin{cases}
& 1,\quad \text{if }0\leq x\leq R \\ 
& 0,\quad \text{if }x\geq R+\delta \\ 
& \text{smoothly decaying}\quad \text{if }R<x<R+\delta%
\end{cases}.%
\end{equation*}%
for arbitrary $R\ge 0$ and $\delta \in [0,1]$. We then show that the laws of the Galerkin approximations are relatively compact. Using this we deduce that equation (\ref{eq:truncatedsystem}) which is the truncated version of  equation (\ref{eq:CMSElocal}) has a global solution, see Theorem  \ref{truncatedscm}. We also show that 
(\ref{eq:truncatedsystem}) has a unique solution. In the deterministic case, the existence of a solution of the truncated equation would immediately imply the existence of a local solution of the original equation on the interval $[t_0,\tau_R]$, where $\tau_R$ is the first time when  $\left\vert \left\vert
\psi _{s}^{R,\delta}\right\vert \right\vert _{\mathcal{H}^{s}(\Omega )}$ reaches the value $R$. We cannot do this here as the solution of the truncated equation (\ref{eq:truncatedsystem}) as well as that of the original equation (\ref{eq:CMSElocal}) depend on temporal values $s$ larger than $\tau_R$ through the quantity 
$\mathbb{E}[\bar{\psi}_{s}^{R,\delta}]$, respectively, $\mathbb{E}[\bar{\psi}_{s}]$.
A final convergence argument is required: We consider a sequence  $\psi^{R,\delta_n}$ with $\delta_n$ tending to $0$. Then the laws of the elements of this sequence are relatively compact and we deduce from here that any limit point of the sequence that satisfies the additional property 
\eqref{adp} is a martingale solution of the original equation (\ref{eq:CMSE}). 
%We can prove that, if $(\psi_1 ,\tau ),$ and $(\psi_2 ,\tau ),$ are two local %solutions that share the same defining stopping time $\tau$, then %$\psi_1=\psi_2$.\footnote{In particular, global solutions of %(\ref{eq:CMSElocal}), if they exists, are unique.}. 

With regards to the uniqueness of the solutions of (\ref{eq:CMSE}): If $(\psi_1 ,\tau_1 ),$ and $(\psi_2 ,\tau_2 ),$ are local solutions that are defined with  \emph{different} stopping times, then we cannot deduce that $\psi_1=\psi_2$ on the common interval of existence $[t_0,\tau_1\wedge\tau_2]$. The reason for this is that, in contrast with the deterministic case, the choice of the stopping time influences $\psi_s$  even for values $s\le \tau$ as a result of the expectation term in (\ref{eq:CMSE}). This lack of consistency between local solutions deters us to construct a pathwise local solution of equation  (\ref{eq:CMSE}) and also a corresponding maximal solution for (\ref{eq:CMSE}).

\begin{remark}
The assumption $\psi\in L^{2}\left( \Xi ;C(\left( [{%
t_{0}},T];\mathcal{H}^{s}(\Omega )\right) \right) $ for any $T\ge 0$ insures that the
term 
\begin{equation*}
\int_{{t_{0}}}^{t\wedge \tau}\mathbb{E}\left[ \bar{\psi}_{s}\right] ds
\end{equation*}%
is well defined as an element of $\mathcal{H}^{0}(\Omega )$. Observe that,
since $\bar{\psi}:=\psi -\frac{1}{|\Omega |}\int_{\Omega }\psi dx$, we have that, for $t\in [t_0,T],$ 
\begin{eqnarray*}
\left\vert \int_{{t_{0}}}^{t\wedge \tau}\mathbb{E}\left[ \bar{\psi}_{s}\right] ds
\right\vert &\leq& \left( t%
-t_{0}\right) \mathbb{E}\left[ \sup_{s\in \left[ t_{0},t\wedge\tau \right]
}\left\vert \left\vert \bar{\psi}_{s}\right\vert \right\vert _{\mathcal{H}%
^{0 }(\Omega )}\right] \\
&\leq& 2\left( T-t_{0}\right) 
\mathbb{E}\left[ \sup_{s\in \left[ t_{0},T \right] }\left\vert \left\vert
\psi _{s}\right\vert \right\vert _{\mathcal{H}^{s}(\Omega )}\right] <\infty .
\end{eqnarray*}%
Moreover, for $0<\left\vert \alpha \right\vert \leq s$, we have that 
\begin{eqnarray*}
\mathcal{D}^{\alpha }\mathbb{E}\left[ \bar{\psi}_{s}\right] 
 &=&\mathbb{E}\left[\mathcal{D}^{\alpha } \psi_{s}\right]   \\
\left\vert \left\vert \mathbb{E}\left[ \mathcal{D}_{s%
}^{\alpha }\psi \right] \right\vert \right\vert _{\mathcal{H}^{0}(\Omega
)}^{2} &\leq &\mathbb{E}\left[ \left\vert \left\vert \mathcal{D}_{s}^{\alpha }\psi \right\vert \right\vert _{\mathcal{H}^{0}(\Omega
)}^{2}\right] \leq \mathbb{E}\left[ \left\vert \left\vert \psi _{s}\right\vert \right\vert _{\mathcal{H}^{s}(\Omega )}^{2}\right]
\end{eqnarray*}%
as the centralizing term vanishes when differentiated.
\end{remark}

The proof of the existence and uniquence of a local solution for the
system (\ref{eq:CMSE}) shares many of the steps with that the existence and
uniquence of the corresponding deterministic case. It uses the same truncated
procedure as in the deterministic case and the same Galerkin approximations.
However several technical difficulties need to be overcome. In the arguments
below we will emphasize these difficulties and the methodology used to
resolve them and omit the arguments that coincide with the deterministic
case. 

We begin by introducing the Galerkin approximation to a truncated
version of the system (\ref{eq:CMSE}). More precisely, let $\psi
^{m,R,\delta}=\left\{ \psi _{t}^{m,R,\delta},t\geq 0\right\} $ be the solution of the
following stochastic differential system%
\begin{equation}
\psi _{t}^{m,R,\delta}=P_{m}\left[ \psi _{{t_{0}}}\right] -\int_{{%
t_{0}}}^{t} F^{m,R,\delta}(\psi _{s}^{m,R,\delta}) ds-\sum_{i=1}^{\infty
}\int_{{t_{0}}}^{t}g_{R,\delta}(||\psi _{s}^{m,R,\delta}||_{\mathcal{H}^{s}(\Omega
)})P_mE_{i}(\psi _{s}^{m,R,\delta})dW_{s}^{i}.
\label{eq:CMSEGalerkinlocal}
\end{equation}%
Just as in (\ref%
{eq:CMSElocal}), the identity in (\ref{eq:CMSEGalerkinlocal}) is assumed to hold, $%
\mathbb{P}$-almost surely, in $L^{2}\left( \Xi ;\mathcal{H}^{0}(\Omega
)\right) $. In (\ref{eq:CMSEGalerkinlocal}), the mapping $F^{m,R,\delta}(\psi
_{s}^{m,R,\delta})$ is defined as 
\begin{eqnarray}
F^{m,R,\delta}(\psi _{s}^{m,R,\delta})&=&g_{R,\delta}(||\psi _{s}^{m,R,\delta}||_{\mathcal{H}^{s}(\Omega
)})\bigg( P_{m}\left[ B\left( \psi _{s}^{m,R,\delta},\psi _{s}^{m,R,\delta}\right) \right]
+C\left( \psi _{s}^{m,R,\delta}\right) \nonumber\\
&&\left.+\frac{1}{2}\sum_{i=1}^{\infty
}P_mE_{i}^{2}(\psi _{s}^{m,R,\delta})-\nu \Delta \psi _{s}^{m,R,\delta}
+ D\mathbb{E}\left[
\bar{\psi}_{s}^{m,R,\delta}\right] \right).
\label{eq:FofGalerkinpsi}
\end{eqnarray}%
The projection operator $P_m$
is defined as in \eqref{GALERKIN_U}+\eqref{GALERKIN_O}. Recall that, when defining the projection corrresponding to the ocean velocity component,  we have taken the incompressibility into account and projected onto the space of divergence-free vector fields.

We then have the following:

\begin{theorem}
\label{truncatedgalerkinscm} Assume that $\psi _{t_{0}}\in \mathcal{H}%
^{s}(\Omega )$. Then the stochastic differential system (\ref%
{eq:CMSEGalerkinlocal}) admits a unique global solution with values in the space 
\begin{equation*}
L^{2}\left( \Xi ;C(\left( [{t_{0}},T];\mathcal{H}^{s}(\Omega )\right)
\right) \cap L^{2}\left( \Xi ;L^{2}(\left( [{t_{0}},T];\mathcal{H}%
^{s+1}(\Omega )\right) \right) .
\end{equation*}%
for any $T>0$. Moreover, there exists a constant $C=C\left( R,T\right) $
independent of $m$ and $\delta$ such that%
\begin{equation}
\mathbb{E}\left[ \sup_{t\in \left[ t_{0},T\right] }||\psi _{s}^{m,R,\delta}||_{%
\mathcal{H}^{s}(\Omega )}^{2}\right] +\mathbb{E}\left[ \int_{{t_{0}}%
}^{T}||\psi _{s}^{m,R,\delta}||_{\mathcal{H}^{s+1}(\Omega )}^{2}\right] \leq C
\label{eq:goodboundG}
\end{equation}%
for any $T>0$.
\end{theorem}

\begin{proof}
Similar to the deterministic case, the system (\ref{eq:CMSEGalerkinlocal}) is
equivalent to a finite dimensional system of stochastic differential
equations of McKean-Vlasov type with Lipschitz continuous coefficients. The same holds true for the system satisfied by $\left( 
\mathcal{D}^{\alpha }\psi ^{m,R,\delta}\right) _{\left\vert \alpha \right\vert \leq
s}$ which involves $\psi ^{m,R,\delta}$ as well as all of its partial derivatives
up to order $s$. The existence and uniqueness of a solution of the system system (\ref{eq:CMSEGalerkinlocal}) follows for example from \cite{sas}. The bound (\ref{eq:goodboundG}) is obtained, as in the
deterministic case, via a Gronwall type argument. 
\end{proof}

\begin{remark}
In addition, one can prove that there exists a constant $C=C\left(
p,R,T\right) $ independent of $m$ such that%
\begin{equation}
\label{lpc}
\mathbb{E}\left[ \sup_{t\in \left[ t_{0},T\right] }||\psi _{s}^{m,R,\delta}||_{%
\mathcal{H}^{s}(\Omega )}^{p}\right] \leq C
\end{equation}%
for any $T>0$ and $p\geq 2$. This bound is useful to show the continuity of
the limit of the Galerkin approximation.
\end{remark}

%%%%%%%

Introduce next $\psi ^{R,\delta}=\left\{ \psi _{t}^{R,\delta},t\geq 0\right\} $ to be the
solution of the following stochastic differential system%
\begin{equation}
\psi _{t}^{R,\delta}= \psi _{{t_{0}}}-\int_{{%
t_{0}}}^{t} F^{R,\delta}(\psi _{s}^{R,\delta}) ds-\sum_{i=1}^{\infty
}\int_{{t_{0}}}^{t}g_{R,\delta}(||\psi _{s}^{R,\delta}||_{\mathcal{H}^{s}(\Omega
)})E_{i}(\psi _{s}^{R,\delta})dW_{s}^{i},
 \label{eq:truncatedsystem}
\end{equation}%
Just as in (\ref%
{eq:CMSElocal}), the identity  (\ref{eq:truncatedsystem}) is assumed to hold, $%
\mathbb{P}$-almost surely, in $L^{2}\left( \Xi ;\mathcal{H}^{0}(\Omega
)\right) $. In (\ref{eq:truncatedsystem}), the mapping $F^{R,\delta}(\psi
_{s}^{R,\delta})$ is defined as 
\begin{eqnarray}
F^{R,\delta}(\psi _{s}^{R,\delta})&=&g_{R,\delta}(||\psi _{s}^{R,\delta}||_{\mathcal{H}^{s}(\Omega
)})\bigg( B\left( \psi _{s}^{R,\delta},\psi _{s}^{R,\delta}\right) 
+C\left( \psi _{s}^{R,\delta}\right) \nonumber\\
&&\hspace{3cm}\left.+\frac{1}{2}\sum_{i=1}^{\infty
}E_{i}^{2}(\psi _{s}^{R,\delta})-\nu \Delta \psi _{s}^{R,\delta}
+ D\mathbb{E}\left[
\bar{\psi}_{s}^{R,\delta}\right] \right).
\label{eq:FofTpsi}
\end{eqnarray}
In the following, we will also need the weak version of the systems (\ref%
{eq:CMSEGalerkinlocal}) and (\ref{eq:truncatedsystem}). These are standard:
For example, the weak version of (\ref{eq:truncatedsystem}) reads as%
\begin{eqnarray}
\left\langle \psi _{t}^{R,\delta},\varphi \right\rangle _{\mathcal{%
H}^{0}(\Omega )^{6}} &=&\left\langle \psi _{{t_{0}}},\varphi \right\rangle _{%
\mathcal{H}^{0}(\Omega )^{6}}-\int_{{t_{0}}}^{t}\left(
\left\langle  F^{R,\varphi }(\psi _{s}^{R,\delta}) ,\varphi \right\rangle _{\mathcal{H}%
^{0}(\Omega )^{6}}\right) ds  \notag \\
&-&\sum_{i=1}^{\infty }\int_{{t_{0}}}^{t}g_{R,\delta}(||\psi _{s}^{R,\delta}||_{\mathcal{H}^{s}(\Omega
)})\left\langle \psi
_{s}^{R,\delta},E_{i}^{\ast }\varphi \right\rangle _{\mathcal{H}^{0}(\Omega
)^{6}}dW_{s}^{i},~~~\varphi \in \left( \mathcal{H}^{2}(\Omega )\right) ^{6},
\label{eq:truncatedsystemweak}
\end{eqnarray}%
where, for $\psi _{t}=(\mathbf{u}_{t}^{a},\theta _{t}^{a},\mathbf{u}%
_{t}^{o},\theta _{t}^{o})$ and $\varphi =\left( \varphi _{i}\right)
_{i=1}^{6}$, $\varphi _{i}\in \mathcal{H}^{2}(\Omega )$ we have the
composite inner product 
\begin{eqnarray*}
\left\langle \psi _{{t}},\varphi \right\rangle _{\mathcal{H}^{0}(\Omega
)^{6}} &=&\left\langle \mathbf{u}_{t}^{a,1},\varphi _{1}\right\rangle _{%
\mathcal{H}^{0}(\Omega )}+\left\langle \mathbf{u}_{t}^{a,2},\varphi
_{2}\right\rangle _{\mathcal{H}^{0}(\Omega )}+\left\langle \theta
_{t}^{a},\varphi _{3}\right\rangle _{\mathcal{H}^{0}(\Omega )} \\
&&+\left\langle \mathbf{u}_{t}^{o,1},\varphi _{4}\right\rangle _{\mathcal{H}%
^{0}(\Omega )}+\left\langle \mathbf{u}_{t}^{o,2},\varphi _{5}\right\rangle _{%
\mathcal{H}^{0}(\Omega )}+\left\langle \theta _{t}^{o},\varphi
_{6}\right\rangle _{\mathcal{H}^{0}(\Omega )}.
\end{eqnarray*}%
As the ocean velocity component takes values in the the space of divergence-free vector fields, we will take the corresponding pair of test function $(\varphi_4,\varphi_5)$ to take value in the same space. The operators $E_{i}^{\ast }$ are adjoint operators corresponding to the
operators $E_{i}$, so that 
\begin{equation*}
\left\langle E_{i}\psi _{s}^{R,\delta},\varphi \right\rangle _{\mathcal{H}%
^{0}(\Omega )^{6}}=\left\langle \psi _{s}^{R,\delta},E_{i}^{\ast }\varphi
\right\rangle _{\mathcal{H}^{0}(\Omega )^{6}}.
\end{equation*}%
In other words, $E_{i}^{\ast }$ are diagonal operators given by 
\begin{eqnarray*}
E_{i}^{\ast }\varphi  &:&=\mathrm{diag}(-\bar{E}^{1}\left( \varphi
_{1},\varphi _{2}\right) ,\ -\bar{E}^{2}\left( \varphi _{1},\varphi
_{2}\right) ,-\mathrm{div}\left( \xi _{i}\varphi _{3}\right) ,\ 0,\ 0,0) \\
\bar{E}^{k}\left( \varphi _{1},\varphi _{2}\right)  &:&=\mathrm{div}\left( \xi
_{i}\varphi _{k}\right) +\frac{1}{Ro^{a}}\xi _{i}^{k}+{\varphi }%
_{j}^{a}\partial _{k}\xi _{i}^{j}+\frac{1}{Ro^{a}}\partial ^{k}(R_{j}(%
\mathbf{x})\xi _{i}^{j}).
\end{eqnarray*}%
The identity in (\ref{eq:truncatedsystemweak}) holds, $\mathbb{P}$-almost
surely, in $L^{2}\left( \Xi ;\mathbb{R}\right) $.

For the following theorem we need to introduce the additional space $%
\mathcal{W}^{\alpha ,p}\left( [{t_{0}},T];\mathcal{H}^{0}(\Omega )\right) $, where $\beta\in (0,1)$ and $p>2$ with 
$\beta p>1$,
defined as 
\[
\mathcal{W}^{\beta ,p}\left( [{t_{0}},T];\mathcal{H}^{0}(\Omega )\right)
:=\left\{ a\in \mathcal{L}_{p}(\left( [{t_{0}},T];\mathcal{H}^{0}(\Omega
)\right) |~~~\Vert a\Vert _{\mathcal{W}^{\beta ,p}(0,T;\mathcal{H}%
^{0}(\Omega ))}<\infty \right\} 
\]%
where the norm $\Vert \cdot \Vert _{\mathcal{W}^{\beta ,p}(0,T;\mathcal{H}%
^{0}(\Omega ))}$ is defined as 
\[
\Vert a\Vert _{\mathcal{W}^{\beta ,p}(0,T;\mathcal{H}^{0}(\Omega
))}^{p}:=\int_{0}^{T}\left\Vert a_{t}\right\Vert _{\mathcal{H}^{0}(\Omega
))}^{p}dt+\int_{0}^{T}\int_{0}^{T}\frac{\left\Vert a_{t}-a_{s}\right\Vert _{%
\mathcal{H}^{0}(\Omega ))}^{p}}{|t-s|^{1+\beta p}}dtds.
\]

\begin{theorem}
\label{truncatedscm}Assume that $\psi _{t_{0}}\in \mathcal{H}^{s}(\Omega )$.
Then the stochastic differential system
(\ref{eq:truncatedsystem}) admits a
unique global solution (in the sense of Definition \ref%
{def:localglobalsolution}) with values in the space 
\begin{equation*}
C\left( \Xi ;C(\left( [{t_{0}},T];\mathcal{H}^{s}(\Omega )\right) \right)
\cap L^{2}\left( \Xi ;L^{2}(\left( [{t_{0}},T];\mathcal{H}^{s+1}(\Omega
)\right) \right) 
\end{equation*}%
for any $T>0$ and any $p\geq 2$. Moreover, there exists a constant $%
C=C\left( R,T\right) $ independent of $\delta$ such that%
\begin{equation}
\mathbb{E}\left[ \sup_{t\in \left[ t_{0},T\right] }||\psi _{s}^{R,\delta}||_{%
\mathcal{H}^{s}(\Omega )}^{p}\right] +\mathbb{E}\left[ \int_{{t_{0}}%
}^{T}||\psi _{s}^{R,\delta}||_{\mathcal{H}^{s+1}(\Omega )}^{2}\right] \leq C
\label{q:goodboundp}
\end{equation}%
for any $T>0$.
\end{theorem}

\begin{proof}$\left. \right. $\newline
\textbf{Existence. }We follow the same steps as in \cite{cfh}, underlying the main
differences below. Let $\left\{ Q^{m}\right\} $ be the family of the
probability laws of the processes $\left\{ \psi ^{m,R,\delta}\right\} .$ These laws
supported on the space%
\begin{eqnarray*}
\mathcal{E}_{0} &:&=\mathcal{E}_{1}\cap \mathcal{E}_{2} \\
\mathcal{E}_{1} &:&=
L^{p}\left( \Xi ; \mathcal{W}^{\beta ,p}\left( [{t_{0}},T];\mathcal{H}^{0}(\Omega )\right)
\right)  \\
\mathcal{E}_{2} &:&=L^{p}\left( \Xi ;C(\left( [{t_{0}},T];\mathcal{H}%
^{s}(\Omega )\right) \right) \cap L^{2}\left( \Xi ;L^{2}(\left( [{t_{0}},T];%
\mathcal{H}^{s+1}(\Omega )\right) \right) ,
\end{eqnarray*}%
where $\beta $ is an arbitrary positive constant such that $\beta <1/2-1/p$ and $p>2$%
. Since $\mathcal{E}_{0}~$is compactly embedded in $L^{p}\left( \Xi
;C(\left( [{t_{0}},T];\mathcal{H}^{0}(\Omega )\right) \right) $, we deduce
that these laws are\ relatively compact in the space of probability measures
$L^{p}\left( \Xi ;C(\left( [{t_{0}},T];\mathcal{H}^{0}(\Omega )\right)\right) .$
We add to the processes $\left\{ \psi ^{m,R,\delta}\right\} $ the driving Brownian
motions $\mathcal{W}=\left\{ W^{i}\right\} _{i=1}^{\infty }$. Then the pairs 
$\{\psi ^{m,R,\delta},\mathcal{W}\}$ have probability laws $\left\{ \widetilde{Q}%
^{m}\right\} $ that are relatively compact in the space of probability
measures over the state space 
\begin{equation*}
L^{p}\left( \Xi ;C(\left( [{t_{0}},T];\mathcal{H}^{0}(\Omega )\right)
\right) \times L^{2}\left( \Xi ;C(\left( [{t_{0}},T];\mathbb{R}\right)
^{\infty }\right) .
\end{equation*}%
Let $\widetilde{Q}$ be a limit point and $\left\{ \widetilde{Q}^{m_{n}}\right\} $ be
a subsequence of measures converging to $\widetilde{Q}$. Also let $\left\{
\varphi _{k}\right\} _{k}$ be a countable dense set of $\left( \mathcal{H}%
^{2}(\Omega )\right) ^{6}.$ By Theorem 2.2 in \cite{kurtzprotter}, it follows that the processes 
\begin{equation*}
\{\psi ^{m_n,R,\delta},\int_{{t_{0}}}^{\cdot }\left\langle \psi _{s}^{{m_n,R,\delta}%
},E_{i}^{\ast }P_{m_{n}}\varphi _{k}\right\rangle _{\mathcal{H}^{0}(\Omega
)^{6}}dW_{s}^{i},i,k=1,...\infty ,\mathcal{W}\}
\end{equation*}%
converge in distribution. By using the Skorohod representation theorem,
there exists a probability space $\left( \widetilde{\Xi},\mathcal{\widetilde{F}},%
\widetilde{P}\right) $ on which we can find

\begin{itemize}
\item a sequence 
\begin{equation*}
\widetilde{A}^{m_{n}}=\left\{ \widetilde{\psi}^{m_n,R,\delta},\int_{{t_{0}}}^{\cdot
}\left\langle \widetilde{\psi}_{s}^{m_n,R,\delta},E_{i}^{\ast }P_{m_{n}}\varphi
_{k}\right\rangle _{\mathcal{H}^{0}(\Omega )^{6}}d\widetilde{W}%
_{s}^{m_{n},i},~~i,k=1,...\infty ,~~\mathcal{\widetilde{W}}^{m_{n}}\right\} 
\end{equation*}%
such that $\left\{ \widetilde{\psi}^{m_n,R,\delta},\mathcal{\widetilde{W}}%
^{m_{n}}\right\} $ has law $\widetilde{Q}^{m_{n}}$

\item a process 
\begin{equation*}
\widetilde{A}=\left\{ \widetilde{\psi}^{R,\delta},\int_{{t_{0}}}^{\cdot }\left\langle 
\widetilde{\psi}_{s}^{R,\delta},E_{i}^{\ast }\varphi _{k}\right\rangle _{\mathcal{H}%
^{0}(\Omega )^{6}}d\widetilde{W}_{s}^{i},~~~i,k=1,...\infty ,~~\mathcal{\widetilde{W}%
}\right\}
\end{equation*}%
with values in the product space 
\begin{equation*}
\widetilde{E}:=L^{2}\left( \widetilde{\Xi};C(\left( [{t_{0}},T];\mathcal{H}%
^{0}(\Omega )\right) \right) \times L^{2}\left( \widetilde{\Xi};C(\left( [{t_{0}}%
,T];\mathbb{R}\right) ^{\mathbb{N}\times \mathbb{N}}\right) \times
L^{2}\left( \widetilde{\Xi};C(\left( [{t_{0}},T];\mathbb{R}\right) ^{\mathbb{N}%
}\right) .
\end{equation*}%
such that the component $\left\{ \widetilde{\psi}^{R},\mathcal{\widetilde{W}}\right\} $
from $\widetilde{A}$ has law $\widetilde{Q}.$

\item The sequence $\left( \widetilde{A}^{m_{n}}\right) _{n}$ converges to $%
\widetilde{A}$ as elements in the space $\widetilde{E}$.

\item Since $Q^{m_{m}}$, the law of $\widetilde{\psi}^{m_n,R,\delta}$ is supported on
the space 
\begin{equation*}
L^{p}\left( \widetilde{\Xi};L^{\infty }(\left( [{t_{0}},T];\mathcal{H}%
^{s}(\Omega )\right) \right) \cap L^{2}\left( \widetilde{\Xi};L^{2}(\left( [{%
t_{0}},T];\mathcal{H}^{s+1}(\Omega )\right) \right) ,
\end{equation*}%
the law of $\widetilde{\psi}^{R,\delta}$ has the same property.
\end{itemize}

Next we take limits of all the terms in the weak version of the equation (%
\ref{eq:CMSEGalerkinlocal}) and show that $\widetilde{\psi}^{R}$ satisfies (\ref%
{eq:truncatedsystemweak})\footnote{%
In (\ref{eq:truncatedsystemweak}), the set of original Brownian motions $%
\mathcal{W}=\left\{ W^{i}\right\} _{i=1}^{\infty }$ is replaced by $\mathcal{%
\widetilde{W}}=\left\{ \widetilde{W}^{i}\right\} _{i=1}^{\infty }$ $.$} for any $%
\varphi _{k}$ in a countable dense set of $\left( \mathcal{H}^{2}(\Omega
)\right) ^{6}$, and therefore, by a density argument, for an arbitrary $%
\varphi $ $\in \left( \mathcal{H}^{2}(\Omega )\right) ^{6}$. The convergence
of all the linear terms is straightforward. We only discuss the convergence
of the nonlinear term, in other words, the limit%
\begin{eqnarray*}
&&\lim_{n\rightarrow \infty }\int_{{t_{0}}}^{t}\left\langle
g_{R}(||\widetilde{\psi}_{s}^{m_{n},R,\delta}||_{\mathcal{H}^{s}(\Omega
)})P_{m_{n}}B\left( \widetilde{\psi}_{s}^{m_n,R,\delta},\widetilde{\psi}%
_{s}^{m_n,R,\delta}\right) ,\varphi _{k}\right\rangle _{\mathcal{H}^{0}(\Omega
)^{6}}ds \\
&&\hspace{3cm}=\int_{{t_{0}}}^{t\wedge \bar{\tau}}\left\langle g_{R}(||\widetilde{\psi}%
_{s}^{R,\delta}||_{\mathcal{H}^{s}(\Omega )})B\left( \widetilde{\psi}_{s}^{R,\delta},\widetilde{%
\psi}_{s}^{R,\delta}\right) ,\varphi _{k}\right\rangle _{\mathcal{H}^{0}(\Omega
)^{6}}ds.
\end{eqnarray*}%
Recall that we took the corresponding pair of test function $(\varphi_4,\varphi_5)$ to be divergence free so we don't need to apply the Leray projection on the corresponding ocean velocity component of the bilinear form $B\left( \widetilde{\psi}_{s}^{R,\delta},\widetilde{%
\psi}_{s}^{R,\delta}\right) $.

The convergence follows via a Sobolev space interpolation argument from the following 
\begin{eqnarray*}
\lim_{n\rightarrow \infty }\mathbb{E}\left[ \sup_{s\in \left[ t_{0},T\right] }\left\vert \left\vert \widetilde{\psi}_{s}^{m_n,R,\delta}-\widetilde{%
\psi}_{s}^{R,\delta}\right\vert \right\vert _{\mathcal{H}^{0}(\Omega )^{6}}^{2}%
\right]  &=&0 \\
\mathbb{E}\left[ \sup_{t\in \left[ t_{0},T\right] }||\widetilde{\psi}%
_{s}^{m_n,R,\delta}||_{\mathcal{H}^{s}(\Omega )}^{p}\right] +\mathbb{E}\left[
\int_{{t_{0}}}^{T}||\widetilde{\psi}_{s}^{m_n,R,\delta}||_{\mathcal{H}^{s+1}(\Omega
)}^{2}\right]  &\leq &C,~~T\geq t_{0} \\
\mathbb{E}\left[ \sup_{t\in \left[ t_{0},T\right] }||\widetilde{\psi}_{s}^{R,\delta}||_{%
\mathcal{H}^{s}(\Omega )}^{p}\right] +\mathbb{E}\left[ \int_{{t_{0}}}^{T}||%
\widetilde{\psi}_{s}^{R,\delta}||_{\mathcal{H}^{s+1}(\Omega )}^{2}\right]  &\leq
&C,~~T\geq t_{0},
\end{eqnarray*}%
where $C=C\left( p,R,T\right) $ is independent of $m_{n}.$

The above argument justifies the existence of a solution of (\ref%
{eq:truncatedsystem}) which is weak in probability sense. From this and the
pathwise uniqueness (see the argument below) of (\ref{eq:truncatedsystem}),
by Yamada-Watanabe theorem, see e.g. \cite{Rockner}
we deduce the existence of a solution of in the
original space driven by the original set of Brownian motions $\mathcal{W}%
=\left\{ W^{i}\right\} _{i=1}^{\infty }$. Note that both the solution of the
equation (\ref{eq:truncatedsystem}) on the original space and 
$\widetilde{\psi}^{R,\delta}$ have the same distribution. Since the law of $%
\widetilde{\psi}^{R,\delta}$ has support on the space 
\begin{equation*}
L^{p}\left( \widetilde{\Xi};C(\left( [{t_{0}},T];\mathcal{H}^{s}(\Omega )\right)
\right) \cap L^{2}\left( \widetilde{\Xi};L^{2}(\left( [{t_{0}},T];\mathcal{H}%
^{s+1}(\Omega )\right) \right) ,
\end{equation*}%
it follows that the solution of the equation (\ref{eq:truncatedsystem})
satisfies (\ref{eq:goodboundG}).\newline
\textbf{Uniqueness.} Let $\psi ^{R,\delta,1}$ and $\psi ^{R,\delta,2}$ be two solutions of
(\ref{eq:truncatedsystem}), in other words, 
\[
\psi _{t}^{R,\delta,i}= \psi _{{t_{0}}}-\int_{{%
t_{0}}}^{t} F^{R,\delta,i}(\psi _{s}^{R,\delta,i}) ds-\sum_{i=1}^{\infty
}\int_{{t_{0}}}^{t}g_{R,\delta,i}(||\psi _{s}^{R,\delta,i}||_{\mathcal{H}^{s}(\Omega
)})E_{i}(\psi _{s}^{R,\delta,i})dW_{s}^{i},\ \ i=1,2
\]
where $F^{R,\delta,i}(\psi
_{s}^{R,\delta,i})$ are defined as 
\begin{eqnarray}
F^{R,\delta,i}(\psi _{s}^{R,\delta,i})&=&g_{R,\delta}(||\psi _{s}^{R,\delta,i}||_{\mathcal{H}^{s}(\Omega
)})\bigg( B\left( \psi _{s}^{R,\delta,i},\psi _{s}^{R,\delta,i}\right) 
+C\left( \psi _{s}^{R,\delta,i}\right) \nonumber\\
&&\hspace{3cm}\left.+\frac{1}{2}\sum_{i=1}^{\infty
}E_{i}^{2}(\psi _{s}^{R,\delta,i})-\nu \Delta \psi _{s}^{R,\delta,i}
+ D\mathbb{E}\left[
\bar{\psi}_{s}^{R,\delta,i}\right] \right).
\label{eq:FofTpsi2}
\end{eqnarray}
The uniqueness argument is now standard: We
use a Gronwall argument. We introduce the following notation 
\begin{eqnarray*}
\psi ^{R,\delta,1,2} &=&\psi ^{R,\delta,1}-\psi ^{R,\delta,2},~~~~\bar{F}^{R,\delta,1,2}=\bar{F}%
^{R}(\psi _{s}^{R,\delta,1})-\bar{F}^{R}(\psi _{s}^{R,\delta,2}) \\
\bar{\psi}_{s}^{R,\delta,1,2} &=&\bar{\psi}_{s}^{R,\delta,1}-\bar{\psi}_{s}^{R,\delta,2}.
\end{eqnarray*}%
Then 
\begin{equation*}
\psi _{t}^{R,\delta,1,2}=\psi _{{t_{0}}}^{1,2}-\int_{{t_{0}}}^{t}\left( \bar{F}%
^{R,\delta,1,2}+\frac{1}{2}\sum_{i=1}^{\infty }E_{i}^{2}(\psi _{s}^{R,\delta,1,2})+D%
\mathbb{E}\left[ \bar{\psi}_{s}^{R,\delta,1,2}\right] \right) ds-\sum_{i=1}^{\infty
}\int_{{t_{0}}}^{t}E_{i}(\psi _{s}^{R,\delta,1,2})dW_{s}^{i},
\end{equation*}%
from which we deduce that 
\begin{eqnarray}
\mathbb{E}\left[ ||\psi _{t}^{R,\delta,1,2}||_{\mathcal{H}^{0}(\Omega )}^{2}\right]
&=&\mathbb{E}\left[ ||\psi _{t_{0}}^{R,\delta,1,2}||_{\mathcal{H}^{0}(\Omega )}^{2}%
\right] -\int_{{t_{0}}}^{t}\mathbb{E}\left[ \left\langle \psi _{s}^{R,\delta,1,2},%
\bar{F}^{R,\delta,1,2}+D\mathbb{E}\left[ \bar{\psi}_{s}^{R,\delta,1,2}\right]
\right\rangle _{\mathcal{H}^{0}(\Omega )}\right] ds  \notag \\
&&+\sum_{i=1}^{\infty }\frac{1}{2}\int_{{t_{0}}}^{t}\mathbb{E}\left[
\left\langle \psi _{s}^{R,\delta,1,2},E_{i}^{2}(\psi _{s}^{R,\delta,1,2})\right\rangle _{%
\mathcal{H}^{0}(\Omega )}+||E_{i}(\psi _{s}^{R,\delta,1,2})||_{\mathcal{H}%
^{0}(\Omega )}^{2}\right] ds.  \label{eq:e1}
\end{eqnarray}%
Observe that 
\begin{eqnarray}
\sum_{i=1}^{\infty }\left( \left\langle \psi _{s}^{R,\delta,1,2},E_{i}^{2}(\psi
_{s}^{R,\delta,1,2})\right\rangle _{\mathcal{H}^{0}(\Omega )}+||E_{i}(\psi
_{s}^{R,\delta,1,2})||_{\mathcal{H}^{0}(\Omega )}^{2}\right)  &\leq &C||\psi
_{s}^{R,\delta,1,2}||_{\mathcal{H}^{0}(\Omega )}^{2}  \label{eq:e2} \\
\mathbb{E}\left[ \left( \left\langle \psi _{s}^{R,\delta,1,2},D\mathbb{E}\left[ 
\bar{\psi}_{s}^{R,\delta,1,2}\right] \right\rangle _{\mathcal{H}^{0}(\Omega
)}\right) \right]  &\leq &C\mathbb{E}\left[ ||\psi _{s}^{R,\delta,1,2}||_{\mathcal{H%
}^{0}(\Omega )}||\mathbb{E}\left[ \bar{\psi}_{s}^{R,\delta,1,2}\right] ||_{\mathcal{%
H}^{0}(\Omega )}\right]   \notag \\
&\leq &C\mathbb{E}\left[ ||\psi _{s}^{R,\delta,1,2}||_{\mathcal{H}^{0}(\Omega )}%
\right] ^{2}  \notag \\
&\leq &C\mathbb{E}\left[ ||\psi _{s}^{R,\delta,1,2}||_{\mathcal{H}^{0}(\Omega )}^{2}%
\right]   \label{eq:e3}
\end{eqnarray}%
Finally, similar to the deterministic case, we deduce that 
\begin{equation}
\left\langle \psi _{s}^{R,\delta,1,2},\bar{F}^{R,\delta,1,2}\right\rangle _{\mathcal{H}%
^{0}(\Omega )}\leq C\left( R\right) ||\psi _{s}^{R,\delta,1,2}||_{\mathcal{H}%
^{0}(\Omega )}^{2}.  \label{eq:e4}
\end{equation}%
From (\ref{eq:e1}), (\ref{eq:e2}), (\ref{eq:e3}) and (\ref{eq:e4}), we
deduce that there exists a constant $C\left( R,T\right) $ such that 
\begin{equation*}
\mathbb{E}\left[ ||\psi _{t}^{R,\delta,1,2}||_{\mathcal{H}^{0}(\Omega )}^{2}\right]
\leq \mathbb{E}\left[ ||\psi _{t_{0}}^{R,\delta,1,2}||_{\mathcal{H}^{0}(\Omega
)}^{2}\right] +C\int_{{t_{0}}}^{t}\mathbb{E}\left[ ||\psi _{s}^{R,\delta,1,2}||_{%
\mathcal{H}^{0}(\Omega )}^{2}\right] ds,~~~t\in \left[ t_{0},T\right] ,
\end{equation*}%
and, by Gronwall's ineguality, we deduce that 
\begin{equation*}
\mathbb{E}\left[ ||\psi _{t}^{R,\delta,1,2}||_{\mathcal{H}^{0}(\Omega )}^{2}\right]
\leq e^{Ct}\mathbb{E}\left[ ||\psi _{t_{0}}^{R,\delta,1,2}||_{\mathcal{H}%
^{0}(\Omega )}^{2}\right] ~~~t\in \left[ t_{0},T\right] .
\end{equation*}%
The continuous dependence of the initial condition implies the uniqueness of
the solution of (\ref{eq:truncatedsystem}).

\end{proof}

We choose next a sequence $\psi ^{R,\delta_n}=\left\{ \psi _{t}^{R,\delta_n},t\geq 0\right\} $ of solutions of the truncated equation \eqref{eq:truncatedsystem} such that $\lim \delta_n=0$. Using arguments similar to those 
applied to the sequence of Galerkin approximations $\psi ^{m, R,\delta}$ one shows that 
these laws of the elements of the  sequence are  relatively compact in the space of probability measures
$L^{p}\left( \Xi ;C(\left( [{t_{0}},T];\mathcal{H}^{0}(\Omega )\right)\right) .$ Via a Skorohod representation theorem, there exists a probability space $\left( \widetilde{\Xi},\mathcal{\widetilde{F}},%
\widetilde{P}\right) $ on which we can find 
a sequence 
$\{ \widetilde{\psi}^{R,\delta_n}\}$ with the same law as the original sequence which converges in $L^{p}\left( \Xi ;C(\left( [{t_{0}},T];\mathcal{H}^{0}(\Omega )\right)\right)$ to a process $\{ \widetilde{\psi}^{R}\}$. 
that satisfies 
\[
\mathbb{E}\left[ \sup_{t\in \left[ t_{0},T\right] }||\widetilde{\psi}_{s}^{R,\delta}||_{%
\mathcal{H}^{s}(\Omega )}^{p}\right] +\mathbb{E}\left[ \int_{{t_{0}}}^{T}||%
\widetilde{\psi}_{s}^{R,\delta}||_{\mathcal{H}^{s+1}(\Omega )}^{2}\right]  \leq
C,~~T\geq t_{0},
\]%
where $C=C\left( p,R,T\right) $.
Via a Sobolev interpolation argument we can also deduce that 
\begin{eqnarray*}
\lim_{n\rightarrow \infty }\mathbb{E}\left[ \sup_{s\in \left[ t_{0},T\right] }\left\vert \left\vert \widetilde{\psi}_{s}^{R,\delta_n}-\widetilde{%
\psi}_{s}^{R}\right\vert \right\vert _{\mathcal{H}^{s-1}(\Omega )^{6}}^{p}%
\right]  &=&0 \\
\lim_{n\rightarrow \infty }\mathbb{E}\left[ \int_{t_{0}}^T\left\vert \left\vert \widetilde{\psi}_{s}^{R,\delta_n}-\widetilde{%
\psi}_{s}^{R}\right\vert \right\vert _{\mathcal{H}^{s}(\Omega )^{6}}^{p}%
ds\right]  &=&0
\end{eqnarray*}
Next let $g_{R}:\mathbb{R}_{+}\rightarrow
\lbrack 0,1]$ be the cut-off function as follows 
\begin{equation*}
g_{R}(x):=%
\begin{cases}
& 1,\quad \text{if }0\leq x\leq R \\ 
& 0,\quad \text{if }x\geq R\\ 
\end{cases}.%
\end{equation*}%
and assume that 
\begin{equation}\label{adp}
\lim_{n\rightarrow \infty }\mathbb{E}\left[ \int_{t_0}^T
\left(g_{R,\delta_n}(\| \widetilde{\psi}_{s}^{R,\delta_n}\| _{\mathcal{H}^{s}(\Omega )^{6}})
-  
g_{R}(\| \widetilde{%
\psi}_{s}^{R}\| _{\mathcal{H}^{s}(\Omega )^{6}})\right)^p dt
\right]  =0.
\end{equation}
If condition \eqref{adp} is satisfied then, by taking 
the limit of each term in equation 
\eqref{eq:truncatedsystemweak}, we can deduce that  %\ref{truncatedscm}
%Let us define the stopping times 
%\begin{eqnarray}\label{st}
%\tau_R &=& \inf\{t\ge t_0 | \ \ %\|\widetilde{\psi}_{s}^{R}\|_{\mathcal{H}^{s}(\Omega )^{6}}^{p} \ge R \} \\
%\tau_{R+\delta_n} &=& \inf\{t\ge t_0 | \ \ %\|\widetilde{\psi}_{s}^{R}\|_{\mathcal{H}^{s}(\Omega )^{6}}^{p} \ge R+\delta_n %\}\label{std}
%\end{eqnarray}
%and, by convention, the stopping times are infinite if the sets on the right %hand side 
%of \eqref{st}, respectively, \eqref{std} are empty. Obviously %$\tau_{R+\delta_n}\ge \tau_{R}$. Let us assume that %$\lim_{n\rightarrow\infty}\tau_{R+\delta_n}=\tau_{R}$ in this case it is %immediate that
the limiting process 
$ \widetilde{\psi}^{R}$ solves the following stochastic differential system%
\begin{equation}
\widetilde\psi _{t}^{R}= \widetilde\psi _{{t_{0}}}-\int_{{%
t_{0}}}^{t} F^{R}(\widetilde\psi _{s}^{R}) ds-\sum_{i=1}^{\infty
}\int_{{t_{0}}}^{t}g_{R}(||\widetilde\psi _{s}^{R}||_{\mathcal{H}^{s}(\Omega
)})E_{i}(\widetilde\psi _{s}^{R})dW_{s}^{i}.
 \label{eq:truncatedsystem0}
\end{equation}%
In (\ref{eq:truncatedsystem0}), the mapping $F^{R}(\widetilde\psi
_{s}^{R})$ is defined as 
\begin{eqnarray}
F^{R}(\widetilde\psi _{s}^{R})&=&g_{R}(||\widetilde\psi _{s}^{R}||_{\mathcal{H}^{s}(\Omega
)})\bigg( B\left( \widetilde\psi _{s}^{R},\widetilde\psi _{s}^{R}\right) 
+C\left( \widetilde\psi _{s}^{R}\right) \nonumber\\
&&\hspace{3cm}\left.+\frac{1}{2}\sum_{i=1}^{\infty
}E_{i}^{2}(\widetilde\psi _{s}^{R})-\nu \Delta \widetilde\psi _{s}^{R}
+ D\mathbb{E}\left[
\bar{\widetilde\psi}_{s}^{R}\right] \right).
\label{eq:FofTpsi3}
\end{eqnarray}
It is then immediate that equation \eqref{eq:truncatedsystem0} is equivalent to \eqref{eq:CMSElocal} where we choose the stopping time
\[
\tau_R := \inf\{t\ge t_0 | \ \ \|\widetilde{\psi}_{s}^{R}\|_{\mathcal{H}^{s}(\Omega )^{6}} \ge R \} 
\]

\begin{remark}By using standard Sobolev interpolation results,
one can prove that 
\begin{equation}\label{adp'}
\lim_{n\rightarrow \infty }\mathbb{E}\left[ \int_{t_0}^T
\left |\| \widetilde{\psi}_{s}^{R,\delta_n}\| _{\mathcal{H}^{s}(\Omega )^{6}}
-  
\| \widetilde{%
\psi}_{s}^{R}\| _{\mathcal{H}^{s}(\Omega )^{6}}\right | dt
\right]  =0.
\end{equation}
This implies that the sequence $\widetilde{\psi}^{R,\delta_n}$ has a subsequence that converges to $\widetilde{\psi}^{R}$ on a set 
$$
\{(\xi\times t)\in \Xi\times [t_0,T]\}
$$
of full $\mathbb P \otimes \ell_{[t_0,T]}$-measure\footnote{That is the complement of the set has null measure}, where $\mathbb \ell_{[t_0,T]}$ is the Lebesgue measure on the interval $[t_0,T]$ . From the definition of $\widetilde{%
\psi}_{s}^{R,\delta_n}$ we can deduce that 
$$
\{(\xi, t)\in \Xi\times [t_0,T]| | 
\| \widetilde{%
\psi}_{t}^{R,\delta_n}(\xi)\| _{\mathcal{H}^{s}(\Omega )^{6}}\le R+\delta_n
\}
$$
and therefore that  
\begin{eqnarray*}
\Xi\times [t_0,T]&\supseteq&\{(\xi,t) | 
\| \widetilde{\psi}_{t}^{R}(\xi)\| _{\mathcal{H}^{s}(\Omega )^{6}}\le R
\}\\
&&=
\{(\xi,t) | 
\| \widetilde{\psi}_{t}^{R}(\xi)\| _{\mathcal{H}^{s}(\Omega )^{6}}< R
\}
\cup
\{(\xi,t) | 
\| \widetilde{\psi}_{t}^{R}(\xi)\| _{\mathcal{H}^{s}(\Omega )^{6}}= R
\}\\
&&=:\mathcal A_{<R}\cup\mathcal A_{=R}
\end{eqnarray*}
has full $\mathbb P \otimes \ell_{[t_0,T]}$-measure. We can deduce from here, via Egorov's theorem
that \eqref{adp} holds true provided $\mathcal A_{=R}$ is a set of null $\mathbb P \otimes \ell_{[t_0,T]}$-measure.\footnote{We thank Tom Kurtz for pointing out this Remark.}   
\end{remark}
\begin{remark}
The argument sofar only shows the existence of a martingale solution (equivalently, a probabilistically weak solution). To prove the existence of
a (probabilistically) strong solution one would need to show the uniqueness of equation 
\eqref{eq:truncatedsystem0}. This cannot be done  as the cut-off function is no longer Lipschitz over the positive half-line. One can try to control the difference between two solutions $\psi^{R,1}$ and $\psi^{R,2}$ up to the minimum of their corresponding hitting times $\tau^1_R \wedge \tau^1_R$. On the interval $[t_0,\tau^1_R \wedge \tau^1_R]$, $g_{R}(||\widetilde\psi _{s}^{R,1}||_{\mathcal{H}^{s}(\Omega)})=g_{R}(||\widetilde\psi _{s}^{R,2}||_{\mathcal{H}^{s}(\Omega)})$  so we can avoid the difficulty raised by $g_{R}$ being non-Lipschitz. However, due to the  
the expectation terms in \eqref{eq:truncatedsystem0}. the solution depends on temporal values beyond $\tau^1_R \wedge \tau^1_R$ which we cannot control. Overcoming this difficulty is beyond the scope of the current paper. 
\end{remark}

\subsection{Lagrangian-Averaged Stochastic Advection by Lie Transport Climate \\Model }\label{LA-SALT-eqns-sec}

Using the same notation as in Section \ref{sec-SALT}, we describe the state of the system is described by a state vector $\psi
:=(\psi ^{a},\psi ^{o})$ with atmospheric component $\psi ^{a}:=(\mathbf{u}%
^{a},\theta ^{a})$ and oceanic component $\psi ^{o}:=(\mathbf{u}^{o},\theta
^{o})$ with initial state is denoted by $\psi (t_{0})=\psi _{0}$. where $%
\psi _{0}=(\mathbf{u}_{0}^{a},\theta _{0}^{a},\mathbf{u}_{0}^{o},\theta
_{0}^{o})$. The LA-SALT equations for the oceanic component $\psi ^{o}$ are the same as the corresponding  SALT equations, in other words $(\mathbf{u}^{o},\theta^{o})$ satisfy equations \eqref{COUPLED_SWE_VELOC_O_STOCH} + \eqref{COUPLED_SWE_T_O_STOCH}. The LA-SALT equations for the atmospheric component $\psi ^{a}$ differ from the the corresponding  SALT equations. More precisely, $  \psi ^{a}:=(\mathbf{u}%
^{a},\theta ^{a})$ satisfy equations 
\eqref{COUPLED_SWE_VELOC_A_STOCH_LA-Ito-Intro} + \eqref{COUPLED_SWE_T_A_STOCH_LA-Ito-Intro}.

In the following we work with the same stochastic basis $(\Xi ,\mathcal{F},(%
\mathcal{F}_{t})_{t},\mathbb{P},(W^{i})_{i})$ and use the same
sequence of vector fields $(\xi _{i})_{i}$  as in Section \ref{sec-SALT}. Similar to (\ref%
{eq:CMSE}), we summarize the LA-CMSE equations \eqref{COUPLED_SWE_VELOC_O_STOCH},\eqref{COUPLED_SWE_T_O_STOCH}, 
\eqref{COUPLED_SWE_VELOC_A_STOCH_LA-Ito-Intro} and \eqref{COUPLED_SWE_T_A_STOCH_LA-Ito-Intro}  as 
\begin{equation}
d\psi _{t}+(B^{L}\left( \psi _{t},\psi _{t}\right) +C^{L}\left( \psi
_{t}\right) )dt+\sum_{i=1}^{\infty }E_{i}^{L}(\psi _{t})\circ dW_{t}^{i}=\nu
\Delta \psi _{t}dt,  \label{eq:LACMSE}
\end{equation}%
where the process $\psi $ gathers all variables in the LASALT model, i.e., $
\psi _{t}:=(\mathbf{u}_{t}^{a},\theta _{t}^{a},\mathbf{u}_{t}^{o},\theta
_{t}^{o})$ (as in the deterministic and SALT cases), $\mathrm{curl}R({\mathbf{x}})=2\Omega ({\mathbf{x}})$, $\bar{\psi}:=\psi -\frac{1}{|\Omega |}\int_{\Omega }\psi dx$ and:\footnote{ As in Section \ref{sec-SALT}, the notation $\left( \cdot \right) ^{T}$ indicates that the
operators $B^{L},C^{L},E_{i}^{L}$ are column vectors.}
\begin{itemize}

\item $B^{L}$ is the bilinear transport operator%
\begin{equation*}
B^{L}\left( \psi _{t},\psi _{t}\right) =(\mathbb{E}[\mathbf{u}^{a}]\cdot
\nabla \mathbf{u}^{a}+{u}_{j}^{a}\nabla \mathbb{E}[\mathbf{u}^{a}]^{j},\ 
\mathbb{E}[\mathbf{u}^{a}]\cdot \nabla \theta ^{a},\ \mathbf{u}^{o}\cdot
\nabla \mathbf{u}^{o},\ \mathbf{u}^{o}\cdot \nabla \theta ^{o})^{T}.
\end{equation*}

\item $C^{L}$ comprises all the linear terms (including the pressure term in
the equation for the components of $\psi $ corresponding to $\bar{\mathbf{u}}%
^{o}$)%
\begin{eqnarray*}
C^{L}\left( \psi _{t}\right)  &=&\left( \frac{1}{Ro^{a}}\mathbb{E}[\mathbf{u}%
^{a}]+\frac{1}{Ro^{a}}\nabla (\mathbb{E}[\mathbf{u}^{a}]\cdot R(\mathbf{x}))+%
\frac{1}{Ro^{a}}\nabla \theta ^{a},\ -\gamma \left( \theta ^{o}-\theta
^{a}\right) ,\right. \  \\
&&\left. \frac{1}{Ro^{o}}\left( \mathbf{u}^{o}\right) ^{\perp }+\frac{1}{%
Ro^{o}}\nabla p^{o}+\sigma \left( \mathbf{u}^{o}-\mathbb{E}[\mathbf{\bar{u}}%
^{a}]\right) ,\ 0\right) ^{T}
\end{eqnarray*}
%\todo[inline]{DC: Why is there a free term in the definition of $E_i^L$ ? }
\item $E_{i}^{L}$ are operators given by 
\begin{equation*}
E_{i}^{L}\left( \psi _{t}\right) =(\xi _{i}\cdot \nabla \mathbf{u}^{a}+\frac{%
1}{Ro^{a}}\xi _{i}+{u}_{j}^{a}\nabla \xi _{i}^{j}+\frac{1}{Ro^{a}}\nabla
(R_{j}(\mathbf{x})\xi _{i}^{j}),\ \xi _{i}\cdot \nabla \theta ^{a},\ 0,\
0)^{T}
\end{equation*}
\end{itemize}

The treatment of equation (\ref{eq:LACMSE}%
) differs slightly from that of (\ref{eq:CMSE}). The reason is that the expected value of state vector $\psi$ satisfies a closed form equation. 
More precisely, as the oceanic component $\psi ^{o}$ is not random, we only need to take expectation in the equations satisfied by the atmospheric component and deduce that 
$\widehat\psi^{a}:=\mathbb{E}[\psi ^{a}]=(\mathbb{E}[\mathbf{u}^{a}],
\mathbb{E}[\theta^{a}])^{T}=:(\widehat{\mathbf{u}}^{a}, \widehat{\theta}^{a})^{T}$ satisfies 
\begin{equation}
d_t\widehat\psi _{t}^a+B^{L,a}( \widehat\psi _{t}^a,\widehat\psi _{t}^a) +C^{L}( \widehat\psi _{t}^a) =\frac{1}{2}\sum_{i=1}^{\infty }E_{i}^{L,a,2}(\widehat\psi _{t}^a)+\nu
\Delta \widehat\psi _{t}^a,  \label{eq:LACMSEE}
\end{equation}%
where 
\begin{eqnarray*}
B^{L,a}\left( \widehat\psi _{t}^a,\widehat\psi _{t}^a\right) &=&(\widehat{\mathbf{u}}^{a}\cdot
\nabla \widehat{\mathbf{u}}^{a}+\widehat{u}_{j}^{a}\nabla (\widehat{\mathbf{u}}^{a})^{j},\ 
\widehat{\mathbf{u}}^{a}\cdot \nabla \widehat{\theta}^{a})^{T}.\\
C^{L}\left( \widehat{\psi}_{t}\right)  &=& \left(\frac{1}{Ro^{a}}
\widehat{\mathbf{u}}^{a}
+\frac{1}{Ro^{a}}\nabla (\widehat{\mathbf{u}}^{a}\cdot R(\mathbf{x}))+%
\frac{1}{Ro^{a}}\nabla \widehat{\theta} ^{a}],\ -\gamma \left( \theta ^{o}-\widehat{\theta}
^{a}\right) \right)^{T}\\
E_{i}^{L,a,2}(\widehat\psi _{t}^a)&=&
 \bigg( 
\frac12  \mathbf{\widehat{z}}\times \xi \Big( {\rm div}\Big(\xi\,\big(\,\mathbf{\widehat{z}}\cdot{\rm curl}\,(\,\widehat{\mathbf{u}}^{a} + \frac{1}{Ro^{a}}\mathbf{R}(\mathbf{x}) \big)\Big) \,\,\Big)  \\
&&\ \ \ \hspace{2mm}- \nabla \bigg( \xi\cdot\nabla\Big(\xi \cdot\big(\widehat{\mathbf{u}}^{a} + \frac{1}{Ro^{a}}\mathbf{R}(\mathbf{x})\big) \Big)
\bigg), \\
&&\ \ \ \hspace{2mm}- \xi\cdot\nabla(\xi_i\cdot\nabla \widehat{\theta}^a) ) \bigg)^{T}.
\end{eqnarray*}

It is immediate that the pair $(\widehat{\psi} ^{a}, \psi ^{o})$ satisfies a system of equations of the form (\ref{COUPLED_SWE_VELOC_A})-(\ref{COUPLED_SWE_INCOMPRESS_O}). More precisely, the only difference between the system of equations satisfied by 
the pair $(\widehat\psi ^{a}, \psi ^{o})$ and the system 
of equations  (\ref{COUPLED_SWE_VELOC_A})-(\ref{COUPLED_SWE_INCOMPRESS_O}) is the linear 
term $\frac{1}{2}\sum_{i=1}^{\infty }E_{i}^{L,a,2}(\widehat\psi _{t}^a)$. This term does not hinder (or help) the analysis of the system  
\eqref{eq:LACMSEE}+\eqref{COUPLED_SWE_VELOC_O_STOCH}+\eqref{COUPLED_SWE_T_O_STOCH}, where, in 
\eqref{COUPLED_SWE_VELOC_O_STOCH}+\eqref{COUPLED_SWE_T_O_STOCH} we replace 
$\mathbb{E}[\bar\psi ^{a}]$ by $\overline{\widehat{\psi}^{a}}:=\widehat{\psi} -\frac{1}{|\Omega |}\int_{\Omega }\widehat{\psi} dx$. 

As for the deterministic case in Theorem \ref{THM_DETERMINISTIC_LOCAL} we have the following 
\begin{theorem}\label{THM_EXPECTATION} 
Let $s\geq 2$ and suppose the initial condition of the 
system 
\eqref{eq:LACMSEE}+\eqref{COUPLED_SWE_VELOC_O_STOCH}+\eqref{COUPLED_SWE_T_O_STOCH} satisfies $\psi_0=(\widehat{\psi}_0^a,\psi_0^o)=(\psi_0^a,\psi_0^o)\in \mathcal{H}^{s}(\Omega)$. 
Then there exists a unique time $t_{e,1}^{*}\in (t_{0}, \infty]$ such that a local regular solution $(\widehat\psi^a,\widehat\psi^o)$ of 
\eqref{eq:LACMSEE}+\eqref{COUPLED_SWE_VELOC_O_STOCH}+\eqref{COUPLED_SWE_T_O_STOCH} in the sense of Definition \ref{REGULAR_SOL} exists and is unique on any interval $T:=[t_0,t_1],$ where $t_0<t_1<t_{e,1}^*$ and that, if  $t_{e,1}^{*}< \infty$, then
\begin{equation}\label{explosione} 
\lim_{t\nearrow t_{e,1}^{*}} ||(\widehat\psi^a,\widehat\psi^o)||_{\mathcal{H}^{s}}=\infty.
\end{equation}
\end{theorem}

The fact that the equation satisfied by the coupled system $(\mathbb E [\psi ^{a}], \psi ^{o})$ has a closed form and, following Theorem \ref{THM_EXPECTATION}, has a local regular solution enables us to show that there exists a solution of the system (\ref{eq:LACMSE}) up to a time  $t_{e,2}^{*}\in (t_{0}, \infty]$. We have the following 

\begin{theorem}\label{th:LAglobal} 
Let $s\geq 2$ and suppose the initial condition of the 
system (\ref{eq:LACMSE}) satisfies $\psi_0\in \mathcal{H}^{s}(\Omega)$. 
Then there exists a unique time $t_{e,2}^{*}\in (t_{0}, \infty]$, such that  on any interval $T:=[t_0,t_1],$ where $t_0<t_1<t_{e,2}^*$, the system 
 (\ref{eq:LACMSE}) has a unique solution with the property that  
\begin{equation}\label{pair0}
\psi^a \in L^{2}\left( \Xi ;C(T;%
\mathbf{H}^{s-2,a}(\Omega)\times H^{s-2}(\Omega) \right) \cup L^{2}\left( \Xi ;L^{2}(T;\mathbf{H}^{s-1,a}(\Omega)\times H^{s-1}(\Omega) \right),
\end{equation}%
 \begin{equation}\label{pair}
(\mathbb E[\psi^a],\psi^o) \in C(T;%
\mathcal{H}^{s}(\Omega )  \cup L^{2}(T;\mathcal{H}^{s+1}(\Omega )  .
\end{equation}%
and that, if  $t_{e,2}^{*}< \infty$, then
\begin{equation}\label{explosionew} 
\lim_{t\nearrow t_{e,2}^{*}} ||(\mathbb E[\psi^a],\psi^o)||_{\mathcal{H}^{s}}=\infty.
\end{equation}
\end{theorem}

\begin{remark}The loss of regularity in the atmospheric component $\psi^a$ as compared to the regularity of the pair  $(\mathbb E[\psi^a],\psi^o)$ is an artifact of our proof. We use Theorems 1 and 2 in \cite{Rozovskii}, Chapter 4, to justify the existence of $\psi^a$ which require additional regularity on the coefficients. This can only be ensured by defining the solution in a lower Sobolev space.
\end{remark}

%The proof of Theorem \ref{th:LAglobal} follows the similar steps with those %in
%the proof of Theorem 1 in \cite{DHL2020}. They are as follows: 
\begin{proof}$\left.\right.$\\
{\bf Existence.} From Theorem 
\ref{THM_EXPECTATION}, we have that there exists a solution of the system of equations \eqref{eq:LACMSEE}+\eqref{COUPLED_SWE_VELOC_O_STOCH}+\eqref{COUPLED_SWE_T_O_STOCH} up to $t_{e,1}^{*}\in (t_{0}, \infty]$ and that, if  $t_{e,1}^{*}< \infty$, then \eqref{explosione} holds.
%\begin{equation}\label{explosionee} 
%\lim_{t\nearrow t_{e,1}^{*}} %||(\widehat\psi^a,\psi^o)||_{\mathcal{H}^{s}}=\infty.
%\end{equation}

We consider next the system 
of equations (\ref{eq:LACMSE}) where we replace $\mathbb{E}[\psi ^{a}]$ 
by the ``atmospheric'' component  $\widehat\psi^a$ which is part of the solution of the system \eqref{eq:LACMSEE}+\eqref{COUPLED_SWE_VELOC_O_STOCH}+\eqref{COUPLED_SWE_T_O_STOCH} . The resulting system $(\check\psi^a,\check\psi^o)$ 
is linear in the stochastic component $\check\psi^a$ and 
nonlinear in the deterministic component $\check\psi^o$. Moreover, the oceanic component $\check\psi^o$ is decoupled from the atmospheric component
as we have replaced the dependence on the atmospheric component  
$\mathbb{E}[\bar\psi ^{a}]$ by $\overline{\widehat{\psi}^{a}}:=\widehat{\psi} -\frac{1}{|\Omega |}\int_{\Omega }\widehat{\psi} dx$. Similar to the proof 
of Theorem \ref{THM_DETERMINISTIC_LOCAL} we deduce that there exists a unique time $t_{e,2}^{*}\in (t_{0}, \infty]$, $t_{e,2}^{*}\le t_{e,1}^{*}$ such that   on any interval $T:=[t_0,t_1],$ where $t_0<t_1<t_{e,2}^*$, the equation satisfied by the 
oceanic component $\check\psi^o$ has a unique solution with the property that 
\begin{equation*}
\check\psi^o \in C(T;%
\mathbf{H}^{s,o}_{div}(\Omega)\times H^{s}(\Omega))  \cup L^{2}(T;\mathbf{H}^{s+1,o}_{div}(\Omega)\times H^{s+1}(\Omega))  .
\end{equation*}%
and, if $t_{e,2}^{*}< t_{e,1}^{*}$, then 
\begin{equation}\label{explosiono} 
\lim_{t\nearrow t_{e,2}^{*}} ||\check\psi^o||_{\mathbf{H}^{s,o}_{div}(\Omega)\times H^{s}(\Omega)}=\infty.
\end{equation}
Crucially, we have that $t_{e,2}^{*}\le t_{e,1}^{*}$ (beyond $t_{e,1}^{*}$ the coefficient of the system satisfied by $(\check\psi^a,\check\psi^o)$ may not be defined, because of blowup exhibited in \eqref{explosione}. 

The linear equation satisfied by the \emph{stochastic component} $\check\psi^a$ 
is a particular case of the equation $(1.1)-(1.2)$ in Chapter 4, Section 4.1, pp.129 in \cite{Rozovskii}. It is easy to check that all assumptions required by Theorem 1 and Theorem 2 in \cite{Rozovskii}, Chapter 4, are fulfilled. 
Therefore the equation satisfied by the stochastic component $\check\psi^a$ has a unique solution defined on the same interval $[t_0,t_{e,2}^{*})$ such that on any interval $T:=[t_0,t_1],$ where $t_0<t_1<t_{e,2}^*$, we have
$\check\psi^a$ belongs to the space stated in \eqref{pair0}.

Because of the linearity of the equation, no blow-up of $\check\psi^a$ is possible before $t_{e,2}^{*}$. 
%The pair $(\check\psi^a,\check\psi^o)$  therefore satisfies \eqref{pair}. 
Moreover $\mathbb E[\check\psi^a]$ satisfies a deterministic linear equation 
which will have a unique solution on the interval $[t_0,t_{e,2}^{*})$ 
in the same space as $\check\psi^a$. However this equations is also satisfied
by $\widehat\psi^a$. It follows that $\mathbb E[\check\psi^a]\equiv\widehat\psi^a$. Moreover the pair $(E[\check\psi^a],\check\psi^o)$ is a solution of the system 
(\ref{eq:LACMSE}) and the pair $(\check\psi^a,\check\psi^o)$  also satisfies \eqref{pair} (because $(\widehat\psi^a,\widehat\psi^o)$ does).

If $t_{e,2}^{*}< t_{e,1}^{*}$, the blow-up at $t_{e,2}^{*}$
holds because of \eqref{explosiono}. 
If $t_{e,2}^{*}= t_{e,1}^{*}<\infty$, then the blow-up at $t_{e,2}^{*}$ holds because of \eqref{explosione}. Hence 
\eqref{explosionew} holds if $t_{e,2}^{*}<\infty$.

{\bf \noindent Uniqueness}. Assume that we have another time $\tilde t_{e,2}^{*}\in (t_{0}, \infty]$, such that  on any interval $T:=[t_0,t_1],$ where $t_0<t_1<t_{e,2}^*$, the system 
 (\ref{eq:LACMSE}) has a unique solution $\tilde\psi$ with the property that  
 
 \begin{equation}\label{pair0dd}
\tilde\psi^a \in L^{2}\left( \Xi ;C(T;%
\mathbf{H}^{s-2,a}(\Omega)\times H^{s-2}(\Omega) \right) \cup L^{2}\left( \Xi ;L^{2}(T;\mathbf{H}^{s-1,a}(\Omega)\times H^{s-1}(\Omega) \right),
\end{equation}%
 \begin{equation}\label{pairdd}
(\mathbb E[\tilde\psi^a],\tilde\psi^o) \in C(T;%
\mathcal{H}^{s}(\Omega )  \cup L^{2}(T;\mathcal{H}^{s+1}(\Omega )  .
\end{equation}%
and that, if  $\tilde t_{e,2}^{*}< \infty$, then
\begin{equation}\label{explosionewt} 
\lim_{ t\nearrow \tilde t_{e,2}^{*}} ||(\mathbb E[\tilde\psi^a],\tilde\psi^o)]||_{\mathcal{H}^{s}}=\infty.
\end{equation}
If  $\tilde t_{e,2}^{*}< t_{e,2}^{*}$, then the system 
 (\ref{eq:LACMSE}) has a unique solution in the interval $[t_0,\tilde t_{e,2}^{*}]$,  
hence $\tilde\psi=\psi$ on $[t_0,\tilde t_{e,2}^{*}]$ and it must be that 
\begin{equation}\label{explosionewtt} 
\lim_{ t\nearrow \tilde t_{e,2}^{*}} ||(\mathbb E[\tilde \psi^a],\tilde\psi^o)]||_{\mathcal{H}^{s}}=\lim_{ t\nearrow \tilde t_{e,2}^{*}} ||(\mathbb E[\psi^a],\psi^o)]||_{\mathcal{H}^{s}}<\infty.
\end{equation}
which contradicts \eqref{explosionewt}. Similarly we cannot have 
$\tilde t_{e,2}^{*}> t_{e,2}^{*}$, so we must have $\tilde t_{e,2}^{*}= t_{e,2}^{*}$ and, by the local uniqueness of the system (\ref{eq:LACMSE}), we also get that $\tilde\psi=\psi$ on the maximal interval of existence. 
\end{proof}
%\todo[inline]{DC: I can add details here.}

From (\ref{eq:LACMSE}) and \eqref{eq:LACMSEE}, we can deduce the equation satisfied by the
fluctuations of the system, i.e., 
\begin{equation*}
\widetilde{\psi}_t:=\psi_t -\mathbb{E}\left[ \psi_t\right] =\psi_t-\widehat \psi_t, \ \ \ t_0\le t< t_{e,1}^{*}.
\end{equation*}%
Then 
\begin{equation}
\widetilde{\psi}_{t}=\widetilde{\psi}_{{t_{0}}}-\int_{{t_{0}}}^{t}\widetilde{F}^{L}(\psi
_{s})ds-\sum_{i=1}^{\infty }\int_{{t_{0}}}^{t}E_{i}^{L}(\psi _{s})dW_{s}^{i},
\label{globalfluctuations}
\end{equation}%
where $\widetilde{F}^{L}(\psi _{s})=F^{L}(\psi _{s})-\mathbb{E}[\left( F^{L}(\psi _{s})\right) ].$

Since only the atmospheric component is random, we can restrict (\ref%
{globalfluctuations}) to $\widetilde{\psi}_{t}^{a}:=(\mathbf{\widetilde{u}}%
_{t}^{a},\ \widetilde{\theta}_{t}^{a})^{T}$ to deduce that 
\begin{equation*}
d\widetilde{\psi}_{t}^{a}=-\widetilde{F}^{L,a}\widetilde{\psi}_{t}^{a}dt-\sum_{i=1}^{%
\infty }E_{i}^{L,a}\psi _{t}^{a}dW_{s}^{i}
\end{equation*}

where 
\begin{eqnarray*}
\widetilde{F}^{L,a}\widetilde{\psi}_{t}^{a} &=&B^{L,a}(\mathbf{\widetilde{u}}^{a},\ 
\widetilde{\theta}^{a})+C^{L,a}(\mathbf{\widetilde{u}}^{a},\ \widetilde{\theta}^{a})-%
\frac{1}{2}\sum_{i=1}^{\infty }\left( E_{i}^{L,a}\right) ^{2}(\mathbf{\widetilde{%
u}}^{a},\ \widetilde{\theta}^{a})-\nu \Delta (\mathbf{\widetilde{u}}^{a},\ \widetilde{%
\theta}^{a})^{T} \\
B^{L,a}\widetilde{\psi}_{t}^{a} &=&\Big(\mathbb{E}[\mathbf{u}^{a}]\cdot \nabla 
\mathbf{\widetilde{u}}^{a}+{\widetilde{u}}_{j}^{a}\nabla \mathbb{E}[\mathbf{u}%
^{a}]^{j},\ \mathbb{E}[\mathbf{u}^{a}]\cdot \nabla \widetilde{\theta}^{a}\Big)^{T}%
\\
&\simeq& \big(\mathcal{L}_{\mathbb{E}[\mathbf{u}^{a}]} (\mathbf{\widetilde{u}}^{a}\cdot d\mathbf{x})%
\,,\,\mathcal{L}_{\mathbb{E}[\mathbf{u}^{a}]}\widetilde{\theta}^{a}\big)^{T}%
\\
C^{L,a}\widetilde{\psi}_{t}^{a} &=&\left( \frac{1}{Ro^{a}}\nabla \widetilde{\theta}%
^{a},\ \gamma \widetilde{\theta}^{a}\right) ^{T} 
\simeq \left(
\frac{1}{Ro^{a}}d\widetilde{\theta}^{a}
,\ 
\gamma \widetilde{\theta}^{a}\right) ^{T} 
\\
E_{i}^{L,a}\psi _{t}^{a} &=&(\xi _{i}\cdot \nabla \mathbf{u}^{a}+{u}%
_{j}^{a}\nabla \xi _{i}^{j},\ \xi _{i}\cdot \nabla \theta ^{a})^{T} %
\simeq \big(\mathcal{L}_{\xi _{i}} (\mathbf{u}^{a}\cdot d\mathbf{x})%
\,,\,\mathcal{L}_{\xi _{i}}\theta^{a}\big)^{T}%
\\
\left( E_{i}^{L,a}\right) ^{2}\widetilde{\psi}_{t}^{a} &=&\big(\xi _{i}\cdot \nabla
\left( \xi _{i}\cdot \nabla \mathbf{\widetilde{u}}^{a}+{\widetilde{u}}_{j}^{a}\nabla
\xi _{i}^{j}\right) 
+\left( \xi _{i}\cdot \nabla \widetilde{u}_j^{a}+{%
\widetilde{u}}_{k}^{a}\partial_j \xi _{i}^{k}\right) \nabla \xi _{i}^{j}\,\,\,,
\ \xi_{i}\cdot \nabla ( \xi _{i}\cdot \nabla \widetilde{\theta}^{a}) \big)^{T}
\\
&\simeq& \big( \mathcal{L}_{\xi _{i}}( \mathcal{L}_{\xi _{i}} (\mathbf{\widetilde{u}}^{a}\cdot d\mathbf{x}))%
\,,\,\mathcal{L}_{\xi _{i}}( \mathcal{L}_{\xi _{i}}\widetilde{\theta}^{a} )\big)^{T}%
\end{eqnarray*}
Here the symbol $\simeq$ recalls the geometric meanings of the coefficients appearing in the fluid equations above. Namely, the left component in the pairs $(\cdot\,\cdot)$ above is understood as the Lie derivative of a 1-form, while the right component is understood as a scalar function. For more discussion of the geometric meanings of these equations, see equation \eqref{EPSD-geom-LASAM-Ito} in Appendix \ref{Appendix-LASALT}.

Define next the variance of the atmospheric component $\Theta ^{a}=\left\{z\Theta _{t}^{a},t\geq 0\right\} $ as  
\begin{equation*}
\Theta _{t}=\mathbb{E}\left[ \left\vert \left\vert \widetilde{\psi}%
_{t}^{a}\right\vert \right\vert _{\mathcal{H}^{0}}^{2}\right] 
\,.\end{equation*}%
Then the previous set of equations implies the following dynamics of the variance %
\begin{equation}
\frac{d\Theta _{t}}{dt}=2\mathbb{E}\left[ \left\langle \widetilde{\psi}_{t}^{a},%
\widetilde{F}^{L,a}\widetilde{\psi}_{t}^{a}\right\rangle _{\mathcal{H}^{0}}\right]
+\sum_{i=1}^{\infty }\mathbb{E}\left[ \left\vert \left\vert E_{i}^{L,a}\psi
_{t}^{a}\right\vert \right\vert _{\mathcal{H}^{0}}^{2}\right] 
\label{variancel}
\end{equation}%
Formula (\ref{variancel}), along with the definitions from above, is important as it can be used to simulate the
dynamics of the variance of the models. In other words, equation %
\eqref{variancel} enables one to compute the statistical dynamics of the
deviations of the fluctuations of the weather that are consistent with the
climatological expectation dynamics. Note that all of the quantities in the original equation combine to influence the fluctuations of the atmospheric component.

In spite of the brevity of formula (\ref{variancel}), we observe that the variance of the fluctuations of the atmospheric component is influenced by all the components of the coupled system $\psi$ as can be seen from the explicit description of the operators $\widetilde{F}^{L,a}$ and  $E_{i}^{L,a}$.

This section has distinguished between the stochastic atmospheric weather model described by the SALT equations in \eqref{COUPLED_SWE_VELOC_A_STOCH} - \eqref{COUPLED_SWE_INCOMPRESS_O_STOCH} as summarised formally in equation \eqref{eq:CMSE} and the atmospheric expectation-fluctuation climate model described by the LA-SALT dynamical equations which are formulated in (\ref{eq:LACMSE}). The latter system of equations has led to the evolution equation \eqref{variancel} which predicts how the evolution of the variances of the atmospheric fluctuations are affected by statistical correlations in their fluctuating dynamics. 

Our analysis for both the SALT and LA-SALT models has determined the local well-posedness properties for the dynamics of the corresponding physical variables. In this well-posed mathematical setting, we have shown that the LA-SALT expectation dynamics can combine with an intricate array of correlations in the fluctuation dynamics to determine the evolution of the mean statistics of the LA-SALT atmospheric climate model.

\section{Summary conclusion and outlook}\label{Sec4}
%Here lies our conclusion.

\begin{enumerate}
    \item We have shown that the Ocean-Atmosphere Climate Model (OACM) in equation set \eqref{COUPLED_SWE_VELOC_A} - \eqref{COUPLED_SWE_INCOMPRESS_O} analysed here for the well-known Gill-Matsuno class of models is simple enough to successfully admit the mathematical analysis required to prove the local well-posedness of these models. The physics underlying these models can be improved, of course. For example, one could naturally include heating by the Greenhouse Effect, and this heating would drive the statistical properties of the climate model. 

    \item In addition to proving well-posedness for both deterministic and stochastic OACM, we have developed a new tool for climate science for predicting the \emph{evolution of climate statistics} such as the variance. Indeed, the application of LA-SALT to the OACM here has established a method for also predicting the evolution of climate statistics such as the expectation of tensor moments of the fluctuations. We believe the new tools introduced here in the context of LA-SALT show considerable promise for future applications.
    
    \item We also expect that the shared geometric structure of these OACM will facilitate the parallel development of their numerical simulations; for example, in undertaking stochastic ENSO simulations. These stochastic simulations would be excellent candidates for the Data Analysis, Uncertainty Quantification and Particle Filtering methods for Data Assimilation which are already under development for SALT. See, e.g.,  \cite{COTTERetal2019,COTTERetal2020a,COTTERetal2020b,COTTERetal2020c}.
\end{enumerate}

%\newpage 

%%%%%%%%%%%%%%%%%%%%%%%%%%%%%%%%%%%%%%%%%%%%%%%%%%%%%%%%%

\appendix
\section{Variational derivations of the deterministic and \\ stochastic atmospheric models} 
\label{sec: Appendix}

\paragraph{Summary.}
In this appendix we explain the Euler-Poincar\'e variational principle and use it to rederive the equations 
of the standard ideal 2D compressible atmospheric and incompressible oceanic flows. 
These ideal models separately conserve their corresponding energy and potential vorticity.
Next, we couple these models using the reduced Lagrange-d'Alembert method. Finally, we introduce stochasticity into the atmospheric flow and derive the full stochastic SALT and LA-SALT climate models which we analyse in the text in both their Stratonovich and It\^o forms.\\

\subsection{Mathematical setting}\label{Appendix-GeoMech}

\begin{definition}[Fluid trajectory]
A fluid trajectory starting from $X\in M$ in the flow domain manifold $M$ at time $t=0$ is given by $x(t)=g_t(X)=g(X,t)$, with $g:M\times\mathbb{R}^+\to M$ being a smooth one-parameter submanifold (i.e., a curve parameterised by time $t$) in the manifold of diffeomorphisms acting on $M$, denoted ${\rm Diff}(M)$. In the deterministic case, computing the time derivative, i.e., the tangent to the curve with initial data $g(X,0)=X$ along $g_t(X)$, leads to the following \emph{reconstruction equation}, given by
\begin{equation}
{\partial_t}g_t(X) = u(g_t(X),t),
\label{eq:reconstructiondeterministic}
\end{equation}
where $u_t(\,\cdot\,)=u(\,\cdot\,,t)$ is a time-dependent vector field whose flow, $g_t(\,\cdot\,)=g(\,\cdot\,,t)$, is defined by the characteristic curves of the vector field $u_t(\,\cdot\,)\in \mathfrak{X}(M)$. The vector fields in $\mathfrak{X}(M)$ comprise the Lie algebra associated to the class of time-dependent maps $g_t\in {\rm Diff}(M)$.  
\end{definition}

\begin{definition}[Advected quantities and Lie derivatives]
A fluid variable $a\in V^*$ defined in a vector space $V^*$ is said to be \emph{advected}, if it keeps its value $a_t=a_0$ along the fluid trajectories. Advected quantities are sometimes called \emph{tracers}, because the histories of scalar advected quantities with different initial values (labels) trace out the Lagrangian trajectories of each label, or initial value, via the \emph{push-forward} by the flow group, i.e., $a_t=g_{t\,*}a_0= a_0g_t^{-1}$, where $g_t\in {\rm Diff}(M)$ is the time-dependent curve on the manifold of diffeomorphisms whose action represents the evolution of the fluid trajectory by push-forward. An advected quantity $a_t$ satisfies an evolutionary partial differential equation (PDE) obtained from the time derivative of the pull-back relation $a_0 = g_t^*a_t$, as follows
\[
0 = {\partial_t}a_0 = {\partial_t}(g_t^*a_t )= g_t^*\big({\partial_t}a_t + \mathcal{L}_u a_t\big)
\quad\Longrightarrow\quad {\partial_t}a_t + \mathcal{L}_u a_t = 0 \,,
\]
where $\mathcal{L}_{u_t} (\,\cdot\,)$ denotes \emph{Lie derivative} by the vector field $u_t$ whose characteristic curves comprise the fluid trajectories. 
\end{definition}

\begin{definition}[Stochastic advection by Lie transport (SALT)\cite{Holm2015}]
In the setting of stochastic advection by Lie transport (SALT) the deterministic reconstruction equation in \eqref{eq:reconstructiondeterministic} is replaced by the semimartingale
\begin{equation}
{\sf d}g(X,t) = u(g_t(X),t)dt + \sum_{i=1}^M \xi_i(g_t(X))\circ dW_t^i,
\label{eq:reconstructionstochastic}
\end{equation}
where the symbol $\circ$ means that the stochastic integral is taken in the Stratonovich sense. The initial data is given by $g(X,0)=X$. The $W_t^i$ are independent, identically distributed Brownian motions, defined with respect to the standard stochastic basis $(\Omega,\mathcal{F},(\mathcal{F}_t)_{t\geq 0},\mathbb{P})$. The $\xi_i(\,\cdot\,)\in\mathfrak{X}$ are prescribed vector fields which are meant to represent uncertainty due to effects on advection of unknown rapid time dependence.

A stochastically advected quantity $a_t$ satisfies an evolutionary stochastic partial differential equation (SPDE) obtained as a  semimartingale relation via the pull-back relation $a_0 = g_t^*a_t$, as follows
\begin{equation}
0 = {\sf d} a_0 = {\sf d}(g_t^*a_t )= g_t^*\big({\sf d} a_t + \mathcal{L}_{{\sf d}x_t} a_t\big)
\quad\Longrightarrow\quad {\sf d} a_t + \mathcal{L}_{{\sf d}x_t} a_t = 0 \,,
\label{eq:KIWformula}
\end{equation}
where $\mathcal{L}_{u_t} (\,\cdot\,)$ denotes \emph{Lie derivative} by the vector field ${\sf d}x_t$ whose characteristic curves comprise the stochastic fluid trajectories in equation \eqref{eq:reconstructionstochastic}. 
\end{definition}

\begin{remark}[Kunita-It\^o-Wentzell (KIW) formula.]
Equation \eqref{eq:KIWformula} is called the \emph{KIW formula}, after its discovery by Kunita as an extension of the It\^o-Wentzell formula to define a stochastic Lie derivative for tensors and differential $k$-forms. For references and a discussion of its recent role in stochastic advection for fluid dynamics, see \cite{AdLHLT2020}.
\end{remark}

\begin{remark}[The deterministic limit.]
In what follows, any of the stochastic fluid equations derived from the Euler-Poincar\'e variational approach will reduce to the corresponding deterministic fluid equations by simply setting $\xi_i\to0$ in the reconstruction equation for the stochastic fluid trajectory in \eqref{eq:reconstructionstochastic}. 
\end{remark}

\begin{definition}[The diamond operator]
The \emph{diamond operator} is defined for $a\in V^*$, $u\in\mathfrak{X}$ and fixed $b\in V$ as
\begin{equation}
\langle b\diamond a, u\rangle_{\mathfrak{X}^{*}\times\mathfrak{X}} := -\langle b,\mathcal{L}_u a\rangle_{V^*\times V}.
\end{equation}
Here, $\langle\,\cdot\,,\,\cdot\,\rangle_{\mathfrak{X}^{*}\times\mathfrak{X}}$ and  $\langle\,\cdot\,,\,\cdot\,\rangle_{\mathfrak{X}^{*}\times\mathfrak{X}}$ denote the real-valued, non-degenerate, symmetric pairings between corresponding dual spaces, which can be defined on a case-by-case basis. 
The diamond operator provides a map dual to the Lie derivative, as $\mathcal{L}_{(\,\cdot\,)}b:\mathfrak{X}\to V$ and $b\diamond(\,\cdot\,):V^*\to\mathfrak{X}^{*}$. This duality is crucial in defining the  Euler-Poincar\'e variational principle.
\end{definition}

\begin{definition}[The variational derivative] The variational derivative of a functional $F:\mathcal{B}\to\mathbb{R}$, where $\mathcal{B}$ is a Banach space, is denoted $\delta F/\delta \rho$ with $\rho\in\mathcal{B}$. The variational derivative $\delta F/\delta \rho$ can be defined via the linearisation of the functional $F$ with respect to the following infinitesimal deformation 
\begin{equation}
\delta F[\rho]:= \frac{d}{d\epsilon}\Big|_{\epsilon=0} F[\rho+\epsilon \delta\rho] = \int \frac{\delta F}{\delta \rho}(x)\delta\rho(x)\,dx =: \left\langle\frac{\delta F}{\delta \rho},\delta \rho\right\rangle.
\end{equation}
In the definition above, $\epsilon\ll1\in\mathbb{R}$ is a parameter, $\delta\rho\in\mathcal{B}$ is an arbitrary function and the  variation can be understood as a Fr\'echet derivative. With the definition of the functional derivative in place, the following lemma can be formulated.
\end{definition}

\begin{theorem}[Stochastic Euler-Poincar\'e theorem]\label{thm:SEP}
With the notation as above, the following statements are equivalent.
\begin{enumerate}[i)]
\item The constrained variational principle
\begin{equation}
\delta\int_{t_1}^{t_2}\ell(u,a)\,dt = 0
\end{equation}
holds on $\mathfrak{X}\times V^*$, using variations $\delta u$ and $\delta a$ of the form
\begin{equation}
\delta u = {\sf d}w - [{\sf d}x_t,w], \qquad \delta a = -\mathcal{L}_w a,
\label{EPstoch-var}
\end{equation}
where $w(t)\in \mathfrak{X}$ is arbitrary and vanishes at the endpoints in time for arbitrary times $t_1,t_2$.
\item The stochastic Euler-Poincar\'e equations
\begin{equation}
{\sf d}\frac{\delta \ell}{\delta u} + \mathcal{L}_{{\sf d}x_t}\frac{\delta \ell}{\delta u} 
= \frac{\delta \ell}{\delta a}\diamond a\,dt\,,
\label{eq:stochep}
\end{equation}
hold on $\mathfrak{X}^*$ and the stochastic advection equations
\begin{equation}
{\sf d}a + \mathcal{L}_{{\sf d}x_t}a = 0\,,
\label{eq:stochadv}
\end{equation}
hold on $\times V^*$.
\end{enumerate}
\end{theorem} 

\begin{proof}
Using integration by parts and the endpoint conditions $w(t_1)=0=w(t_2)$, the variation can be computed to be
\begin{equation}
\begin{aligned}
\delta\int_{t_1}^{t_2}\ell(u,a)\,dt 
&= 
\int_{t_1}^{t_2}\left\langle\frac{\delta\ell}{\delta u},\delta u\right\rangle + \left\langle\frac{\delta\ell}{\delta a},\delta a\right\rangle\,dt\\
&= \int_{t_1}^{t_2}\left\langle\frac{\delta\ell}{\delta u},{\sf d}w-[{\sf d}x_t,w]\right\rangle + \left\langle\frac{\delta\ell}{\delta a}\,dt,-\mathcal{L}_w a\right\rangle\\
&= \int_{t_1}^{t_2}\left\langle -{\sf d}\frac{\delta\ell}{\delta u} - \mathcal{L}_{{\sf d}x_t}\frac{\delta\ell}{\delta u} + \frac{\delta\ell}{\delta a}\diamond a\,dt,w\right\rangle\\
&= 0\,.
\end{aligned}
\label{eq:StochEPeqns}
\end{equation}
Since the vector field $w$ is arbitrary, one obtains the stochastic Euler-Poincar\'e equations. Finally, the advection equation \eqref{eq:stochadv} follows by applying the KIW formula to $a(t)=g_{t*}a_0$.
\end{proof}

\begin{remark}
This version of the stochastic Euler-Poincar\'e theorem is equivalent to the version presented in \cite{Holm2015}, which uses stochastic Clebsch constraints. In \cite{Holm2015} one can also find an investigation of the It\^o formulation of the stochastic Euler-Poincar\'e equation. 
\end{remark}

\begin{theorem}[Stochastic Kelvin-Noether Theorem] \label{thm:Kelvin}
Let $c$ denote a compact embedded one-dimensional smooth submanifold of $M$ and denote $c_t=g_t(c)$ for all $t\in [0,T]$. If the mass density $D_0$ (a top-form) is initially non-vanishing, then
\[
{\sf d}\oint_{c_t}\frac{1}{D_t}\frac{{\delta} \ell}{{\delta} u}(u_t,a_t)
= \oint_{c_t}\frac{1}{D_t}\frac{{\delta} \ell}{{\delta} a}(u_t,a_t)\diamond a_t\,.
\]
The integrated form of this relation is 
\[
\oint_{c_t}\frac{1}{D_t}\frac{{\delta} \ell}{{\delta} u}(u_t,a_t)
=\oint_{c_0}\frac{1}{D_0}\frac{{\delta} \ell}{{\delta} u}(u_0,a_0)
+ \int_0^t\oint_{c_s}\frac{1}{D_s}\frac{{\delta} \ell}{{\delta} a}(u_s,a_s)\diamond a_s \rmd s.
\]
\end{theorem}

\begin{proof}
The KIW formula \eqref{eq:KIWformula} -- also known as the \emph{Lie chain rule} -- implies that 
\begin{align*}
{\sf d}\oint_{c_t}\frac{1}{D_t}\frac{{\delta} \ell}{{\delta} u}(u_t,a_t)
&=
\oint_{c_0}
{\sf d}
g_t^*\Big(\frac{1}{D_t} \frac{\delta\ell}{\delta u}\Big)\bigg)
=
\oint_{c_0}
g_t^*\bigg(\big({\sf d} + \mathcal{L}_{{\sf d}x_t} \big)\Big(\frac{1}{D_t} \frac{\delta\ell}{\delta u}\Big)\bigg)
\\&=
\oint_{c_t}
\big({\sf d} + \mathcal{L}_{{\sf d}x_t} \big)\Big(\frac{1}{D_t} \frac{\delta\ell}{\delta u}\Big)
= \oint_{c_t}\frac{1}{D_t}\frac{{\delta} \ell}{{\delta} a}(u_t,a_t)\diamond a_t\,,
\end{align*}
where the first step also uses the stochastic advection equation for $D_t$ in \eqref{eq:stochep} and 
the final step follows by substituting the stochastic Euler-Poincar\'e equations in \eqref{eq:stochadv}.

\end{proof}

\subsection{2D SALT Stochastic Atmospheric Model (SAM) }\label{Appendix-SALT}

\paragraph{Summary.}
This part of the appendix derives the SALT atmospheric model (SAM) for the isothermal ideal gas in equations \eqref{Atmos-SALT-Kel} and \eqref{Stoch_THETA_A}. \\ 

In the present notation, the Lagrangian for the deterministic 2D Compressible Atmospheric
Model (DAM) in Eulerian $(x,y)$ coordinates is,
\begin{align}
\ell\big[\bs{u}, D,\theta \big] = 
\int_{\Omega} \frac{D}{2} |\bs{u}|^2 
+ D \bs{u}\cdot\bs{R}(\bs{x})
- c_vD\theta\Pi  \,dx\,dy,
\label{CAM-Lag}
\end{align}
where $\bs{u}$ denotes 2D fluid velocity, $D$ is mass density, $\theta$ denotes the potential temperature, 
the function $\bs{R}(\bs{x})$ with ${\rm curl} \bs{R}=2\bs{\Omega}$ denotes the vector potential for the Coriolis parameter, 
$c_v$ is specific heat at constant volume, and $\Pi$ is the well-known Exner function, given by
\[
\Pi = \left(
\frac{p}{p_0}
\right)^{R/c_p},
\]
in which  $p_0$ is a reference pressure level, $c_p$ is specific heat at constant pressure and $R=c_p-c_v$ is the gas
constant. In these variables, the equation of state for an ideal gas in 2D with $n$ degrees of freedom is expressed as
\[
\Pi = \left(\frac{RD\theta}{p_0}\right)^{R/c_p} = \left(\frac{RD\theta}{p_0}\right)^{1-\gamma^{-1}}
= \left(\frac{RD\theta}{p_0}\right)^{2/(n+2)}
\,,
\]
since the specific heat ratio $\gamma=c_p/c_v = 1+2/n$ for ideal gases whose molecules possess $n$ degrees of freedom,   
comprising spatial translations, rotations and oscillations. Ideal gases of diatomic molecules in 3D have three translations, 
plus rotations and oscillations, so $n=5$ and $\gamma=7/5$ in that case.

Accordingly, the Lagrangian in \eqref{CAM-Lag} specialises for a ideal gas in 2D to 
\begin{align}
\ell\big[\bs{u}, D,\theta \big] = 
\int_{\Omega} \frac{D}{2} |\bs{u}|^2 
+ D \bs{u}\cdot\bs{R}(\bs{x})
- \kappa (D\theta)^\alpha  \,dx\,dy,
\label{2DAM-Lag}
\end{align}
where the constants $(\kappa,\alpha)$ take the values, 
\[
\kappa=c_v (R/p_0)^{2/(n+2)} 
\quad\hbox{and}\quad
\alpha = \frac{n+4}{n+2} = 1 + \frac{2}{n+2} = 2 -\gamma^{-1}
\,.
\]

We obtain the following variational derivatives of the Lagrangian in \eqref{CAM-Lag},
\begin{align}
\begin{split}
\frac1D\dede{\ell}{\bs{u}}  &=  \bs{u} + \bs{R}(\bs{x})
\,, \\
\dede{\ell}{D} & =  \frac{1}{2} |\bs{u}|^2 +  \bs{u}\cdot\bs{R}(\bs{x}) - \kappa\alpha(D\theta)^{\alpha-1}\theta \,, \\
\dede{\ell}{\theta}  &= -\, \kappa\alpha(D\theta)^{\alpha-1}D
\,.
\end{split}
\label{CAM-Lag-vars}
\end{align}

Substitution of the variational derivatives (\ref{CAM-Lag-vars}) of
the Lagrangian (\ref{CAM-Lag}) into the stochastic Euler-Poincar\'e equations
with SALT in \eqref{eq:stochep} gives the SAM system
 \begin{align}
\begin{split}
\left({\sf d} + \mathcal{L}_{{\sf d}x_t}\right) \left( \big( \bs{u} + \bs{R}(\bs{x}) \big) \cdot d\bs{x} \right)
 &=  \,{d\,}\left(\frac{1}{2} |\bs{u}|^2 +  \bs{u}\cdot\bs{R}(\bs{x}) -(\kappa/\gamma) (D\theta)^{\alpha} \right) dt
\,,\\
({\sf d} + \mathcal{L}_{{\sf d}x_t}) (D\theta \,dxdy) & =  0
\,.
\end{split}
\label{EPSD-geom-CAM}
\end{align}
Consequently, we recover the Kelvin circulation conservation law for the SAM in the compact form
\begin{align}
{\sf d}\oint_{c_t}\hspace{-2mm} 
 \big( \bs{u} + \bs{R}(\bs{x}) \big) \cdot d \bs{x} 
=
0 
\,,
\label{SAM-circons}
\end{align}
where $c_t=g_t(c_0)$ for all $t\in [0,T]$ denotes the push-forward by the SAM flow of the initial $c_0$,
a compact embedded one-dimensional smooth submanifold of $M$.

\begin{corollary}\label{CAM-PV}
  The system of SAM equations in \eqref{EPSD-geom-CAM} implies that
  potential vorticity $q:= \omega / (D\theta)$ is conserved along flow lines of the stochastic fluid
  trajectory ${\sf d}x_t$,
\begin{align}
{\sf d} q +{\sf d}\bs{x}_t \cdot\nabla q = 0
\quad\hbox{with potential vorticity }q:= \omega / (D\theta) 
\quad\hbox{and} \quad
\omega := \bs{\widehat{z}}\cdot{\rm curl}\big( \bs{u} + \bs{R}(\bs{x}) \big) 
\,.
\label{SAM-vorticitythm}
\end{align}
In turn, this formula implies that the following infinite family of integral quantities is conserved
\begin{align}
C_\Phi = \int_\Omega (D\theta)\Phi(q)\,dxdy
\,,
\label{SAM-enstrophy-thm}
\end{align}
for any differentiable function $\Phi$.
\end{corollary}
\begin{corollary}
The \emph{Deterministic} AM equations in (\ref{EPSD-geom-CAM}) with $\bs{\xi}_i\to0$ are Hamiltonian, with conserved energy\footnote{The energy $E$ in \eqref{SAM-erg} is not conserved for $\bs{\xi}_i\ne0$, though, because $\bs{\xi}_i\ne0$ injects the explicit time dependence of stochastic Lagrangian trajectories into the Euler-Poincar\'e variations \eqref{EPstoch-var} in Hamilton's principle.}
\begin{align}
E = \int_{\Omega} \frac{D}{2} |\bs{u}|^2 + \kappa (D\theta)^\alpha \,dxdy.
\label{SAM-erg}
\end{align}
The Lie-Poisson Hamiltonian structure of deterministic fluid equations is discussed in \cite{HMR1998} from the viewpoint of the Euler-Poincar\'e of Hamilton's principle for fluid dynamics. 
\end{corollary}
\begin{remark}
The system of SAM equations in \eqref{EPSD-geom-CAM} may also be written equivalently in \emph{standard} fluid dynamics notation as
\begin{align}
\begin{split}
{\sf d} \bs{u} + {\sf d}\bs{x}_t \cdot\nabla \bs{u} - {\sf d}\bs{x}_t  \times 2 \bs{\Omega}
 + \sum_{i=1}^M \Big( u_j  \nabla  \xi_i^j  +  \nabla \big(\bs{\xi}_i  \cdot \bs{R}\big)\Big) \circ dW_t^i 
&=  -\,(\kappa/\gamma) \nabla (D\theta)^\alpha  dt
\,, \\
{\sf d} (D\theta) + \nabla  \cdot \big( D\theta\, {\sf d}\bs{x}_t  \big) & =  0
\,.\end{split}
\label{Eady-EPSDeqns2}
\end{align}
\end{remark}
If one assumes low Mach number, so that $D\approx1$, and then adds viscosity and diffusion of heat, this set of equations will reproduce the SALT atmospheric model in equations \eqref{Atmos-SALT-Kel} and \eqref{Stoch_THETA_A} when one also sets $\alpha=1=\gamma$, which holds for the isothermal case of the ideal gas. 

%%%%%%%%%%%%%%%%%%%%%%%%%%%%%%%%%%%%%%%%%%%%%%%%%%%%%%%%%%%%%%%%%%%%%%%%%%%%%%%%%%%%%%%
%\begin{comment}
%%%%%%%%%%%%%%%%%%%%%%%%%%%%%%%%%%%%%%%%%%%%%%%%%%%%%%%%%%%%%%%%%%%%%%%%%%%%%%%%%%%%%%%
\subsection{2D LA-SALT Stochastic Atmospheric Model (LASAM) }\label{Appendix-LASALT}

\paragraph{Summary.}
Here we provide a geometric derivation of the Lagrangian Averaged Stochastic Advection by Lie Transport (LA-SALT) atmospheric model (LASAM) for the isothermal ideal gas studied in section \ref{LA-SALT-eqns-sec}. \\

The simplest way to derive the 2D LA-SALT Stochastic Atmospheric Model (LASAM) is to alter the \emph{Stratonovich} stochastic Lagrangian trajectory in equations \eqref{EPSD-geom-CAM} to the \emph{Stratonovich} LA-SALT stochastic path, which reads
 
\begin{equation}
dx_{t}^{a}\rightarrow \ {{\sf d}\mathbf{X}_{t}^{a}} 
:= {\mathbb{E}[\mathbf{u}^{a}]}(\bs{x},t)dt + \sum_{i}\xi _{i}^{a}(\bs{x})\circ dW_{i}(t).
\label{LA-SALTpath}
\end{equation}

The LASAM system then emerges in \emph{Stratonovich} stochastic geometric form as
 \begin{align}
\begin{split}
\left({\sf d} + \mathcal{L}_{{{\sf d}\mathbf{X}_{t}^{a}}} \right)  
\big( \bs{u} + \frac{1}{Ro^{a}}\bs{R}(\bs{x}) \big) \cdot d\bs{x}  
 &=  \,{d\,}\left(\frac{1}{2} |\bs{u}|^2 +  \frac{1}{Ro^{a}}\bs{u}\cdot\bs{R}(\bs{x}) 
 -\frac{1}{Ro^{a}}(\kappa/\gamma) (D\theta)^{\alpha} \right) dt
\,,\\
\left({\sf d} + \mathcal{L}_{{{\sf d}\mathbf{X}_{t}^{a}}}\right) (D\theta \,dxdy) & =  0
\,,
\end{split}
\label{EPSD-geom-LASAM-Strat}
\end{align}
where $\mathcal{L}_{{{\sf d}\mathbf{X}_{t}^{a}}}$ denotes Lie derivative with respect to \emph{Stratonovich} LA-SALT stochastic path in \eqref{LA-SALTpath}.
The stochastic geometric LASAM system in equation \eqref{EPSD-geom-LASAM-Strat} is given in its corresponding \emph{It\^o} form by
 \begin{align}
\begin{split}
\left({\sf d} + \mathcal{L}_{{{\sf d}\mathbf{\wh{X}}_{t}^{a}}}\right)
\Big( \big( \bs{u} + \frac{1}{Ro^{a}}\bs{R}(\bs{x}) \big) \cdot d\bs{x} \Big) 
& - \frac12  \sum_{i} \mathcal{L}_{\xi _{i}^{a}}\left(\mathcal{L}_{\xi _{i}^{a}}
\Big( \big( \bs{u} + \frac{1}{Ro^{a}}\bs{R}(\bs{x}) \big) \cdot d\bs{x} \Big) \right)
\\&=  \,{d\,}\left(\frac{1}{2} |\bs{u}|^2 +  \frac{1}{Ro^{a}}\bs{u}\cdot\bs{R}(\bs{x}) 
 -\frac{1}{Ro^{a}}(\kappa/\gamma) (D\theta)^{\alpha} \right) dt
\,,\\
\left({\sf d} + \mathcal{L}_{{{\sf d}\mathbf{\wh{X}}_{t}^{a}}}\right) (D\theta \,dxdy) 
& - \frac12  \sum_{i} \mathcal{L}_{\xi _{i}^{a}}\left(\mathcal{L}_{\xi _{i}^{a}} (D\theta \,dxdy) \right)
 =  0
\,,
\end{split}
\label{EPSD-geom-LASAM-Ito}
\end{align}
where the \emph{It\^o} Lagrangian trajectory for the LA-SALT, reads
 
\begin{equation*}
{{\sf d}\mathbf{\wh{X}}_{t}^{a}} 
:= {\mathbb{E}[\mathbf{u}^{a}]}(\bs{x},t)dt + \sum_{i}\xi _{i}^{a}(\bs{x})dW_{i}(t).
\end{equation*}

In standard notation for fluid dynamics and with $\alpha=1=\gamma$, the \emph{Stratonovich} LASAM equations in \eqref{EPSD-geom-LASAM-Strat} become

\begin{align}
d\mathbf{u}^{a}+({d\mathbf{X}_{t}}^{a}\cdot \nabla )\mathbf{u}
^{a}+\frac{1}{Ro^{a}}{d\mathbf{X}_{t}}^{a\bot }
& +{
\sum_{i}\Big(u_{j}^{a}\nabla \xi _{i}^{j}+\frac{1}{Ro^{a}}\nabla \Big(R_{j}
\mathbf{(x})\xi _{i}^{j}\Big)\Big)\circ dW_{t}^{i}} 
 \notag\\
& \hspace{-22mm}{+\,u_{j}^{a}\nabla \mathbb{E}[{u^{a}}^{j}]dt+
\frac{1}{Ro^{a}}\nabla (\mathbb{E}[{\mathbf{u}}^{a}]\cdot \mathbf{R})dt}+
\frac{1}{Ro^{a}}\nabla \theta ^{a}\,dt=\frac{1}{Re^{a}}\triangle \mathbf{u}^{a}\,dt\,,
 \notag\\
& d\theta ^{a}+{d\mathbf{X}_{t}}^{a}\cdot \nabla \theta
^{a}=-\gamma (\theta ^{o}-\theta ^{a})\,dt + \frac{1}{Pe^{a}}\triangle \theta ^{a}\,dt\,.
\label{COUPLED_SWE_T_A_STOCH_LA-Strat}
\end{align}
In the previous equation, we have used the continuity equation to eliminate the areal density $D$.

Likewise, in standard notation for fluid dynamics and with $\alpha=1=\gamma$, the \emph{It\^o} LASAM equations in \eqref{EPSD-geom-LASAM-Ito} become

\begin{align}
& d\mathbf{u}^{a}+({d\mathbf{\wh{X}}_{t}}^{a}\cdot \nabla )\mathbf{u}
^{a}+\frac{1}{Ro^{a}}{d\mathbf{\wh{X}}_{t}}^{a\bot }
 + {
\sum_{i}\Big(u_{j}^{a}\nabla \xi _{i}^{j}+\frac{1}{Ro^{a}}\nabla \Big(R_{j}
\mathbf{(x})\xi _{i}^{j}\Big)\Big) dW_{t}^{i}}  \nonumber \\
&+\frac12 \bigg[ \mathbf{\hat{z}}\times \xi \Big( {\rm div}\Big(\xi\,\big(\,\mathbf{\hat{z}}\cdot{\rm curl}\,(\,\mathbb{E}[{\mathbf{u}}^{a}] + \frac{1}{Ro^{a}}\mathbf{R}(\mathbf{x}) \big)\Big) \,\,\Big)  
- \nabla \bigg( \xi\cdot\nabla\Big(\xi \cdot\big(\mathbb{E}[{\mathbf{u}}^{a}] + \frac{1}{Ro^{a}}\mathbf{R}(\mathbf{x})\big) \Big)
\bigg)\bigg]dt
\nonumber \\& \hspace{22mm}{+\,u_{j}^{a}\nabla \mathbb{E}[{u^{a}}^{j}]dt+
\frac{1}{Ro^{a}}\nabla (\mathbb{E}[{\mathbf{u}}^{a}]\cdot \mathbf{R})dt}+
\frac{1}{Ro^{a}}\nabla \theta ^{a}\,dt = \frac{1}{Re^{a}}\triangle \mathbf{u}^{a}\,dt\,,
\label{COUPLED_SWE_VELOC_A_STOCH_LA-Ito} \\
& d\theta ^{a}+{d\mathbf{\wh{X}}_{t}}^{a}\cdot \nabla \theta
^{a} - \frac12 \Big(\xi\cdot\nabla(\xi\cdot\nabla \theta^a)  \Big)dt 
=-\gamma (\theta ^{o}-\theta ^{a})\,dt +  \frac{1}{Pe^{a}}\triangle \theta ^{a}\,dt
\,.
\label{COUPLED_SWE_T_A_STOCH_LA-Ito}
\end{align}
In the final equation, we have again used the continuity equation to eliminate the areal density $D$.

\begin{remark}[Expectation equations]
Taking the expectation of equations \eqref{COUPLED_SWE_VELOC_A_STOCH_LA-Ito} and  \eqref{COUPLED_SWE_T_A_STOCH_LA-Ito} yields a closed set of deterministic PDE for the expectations 
$\mathbb{E}[{\mathbf{u}}^{a}]$ and $\mathbb{E}[\theta ^{a}]$. Subtracting the expectations from equations \eqref{COUPLED_SWE_VELOC_A_STOCH_LA-Ito} and  \eqref{COUPLED_SWE_T_A_STOCH_LA-Ito} yields \emph{linear equations} for the differences,   
\begin{align*}
{\mathbf{u}}^{a'}:= {\mathbf{u}}^{a} - \mathbb{E}[{\mathbf{u}}^{a}]
\quad\hbox{and}\quad
\theta ^{a'} := \theta ^{a} - \mathbb{E}[\theta ^{a}]
\,.
\label{fluctuations_u_theta}
\end{align*}
Since ${\mathbf{u}}^{a'}$ and $\theta ^{a'}$ satisfy $\mathbb{E}[{\mathbf{u}}^{a'}]=0$ and $\mathbb{E}[\theta ^{a'}]=0$, one may regard these difference variables as fluctuations of ${\mathbf{u}}^{a}$ and $\theta ^{a}$ away from their expected values. 
\end{remark}

\end{document}